\documentclass[a4paper,12pt,reqno]{amsart}

\usepackage{amsmath,amsfonts,amssymb,amsthm}
\usepackage{graphics}
\usepackage{epsfig}
\usepackage{xcolor}
\usepackage{amsrefs}
\usepackage[margin=1.4in]{geometry}
\usepackage{multirow}
\usepackage{chngpage}

\usepackage{siunitx}

\newcommand{\ignore}[1]{}

\def\Xint#1{\mathchoice
{\XXint\displaystyle\textstyle{#1}}%
{\XXint\textstyle\scriptstyle{#1}}%
{\XXint\scriptstyle\scriptscriptstyle{#1}}%
{\XXint\scriptscriptstyle%
\scriptscriptstyle{#1}}%
\!\int}
\def\XXint#1#2#3{{\setbox0=\hbox{$#1{#2#3}{%
\int}$ }
\vcenter{\hbox{$#2#3$ }}\kern-.6\wd0}}

\def\dashint{\Xint-}
\newcommand{\av}{\dashint}

\newtheorem{definition}{Definition}
\newtheorem{proposition}{Proposition}
\newtheorem{theorem}{Theorem}

\newtheorem{lemma}{Lemma}
\newtheorem{corollary}{Corollary}

\newcommand{\omegay}{\omega_y}
\newcommand{\sigmaL}{\sigma^{(L)}}
\newcommand{\phiL}{\phi^{(L)}}
\newcommand{\qL}{q^{(L)}}
\newcommand{\aL}{a^{(L)}}
\newcommand{\GL}{G^{(L)}}
\newcommand{\rL}{r^{(L)}}
\newcommand{\uL}{u^{(L)}}
\newcommand{\tuL}{\tilde{u}^{(L)}}
\newcommand{\xiL}{\xi^{(L)}}
\newcommand{\gL}{g^{(L)}}

\newcommand{\abs}[1]{\lvert#1\rvert}

\title{Optimal artificial boundary condition for random elliptic media}

\author{Jianfeng Lu} \address{Mathematics Department, Duke University,
  Box 90320, Durham, NC 27708, USA} \email{jianfeng@math.duke.edu}

\author{Felix Otto} \address{Max Planck Institute for Mathematics in
  the Sciences, Inselstr. 22 04103 Leipzig, Germany} \email{otto@mis.mpg.de}

\date{\today}

\begin{document}



\maketitle 

\begin{abstract}
We are given a uniformly elliptic coefficient field 
that we regard as a realization of a stationary and finite-range (say, range unity) ensemble of coefficient fields.
Given a (deterministic) right-hand-side supported in a ball of size $\ell\gg 1$ and of vanishing average, 
we are interested in an algorithm to compute the (gradient of the) solution near the origin, 
just using the knowledge of the (given realization of the) coefficient field in some large box of size $L\gg\ell$. 
More precisely, we are interested in the most seamless (artificial) boundary condition
on the boundary of the computational domain of size $L$. 

\smallskip

Motivated by the recently introduced multipole expansion in random media, we propose an algorithm. 
We rigorously establish an error estimate (on the level of the gradient) in terms of $L\gg\ell\gg 1$, 
using recent results in quantitative stochastic homogenization. More precisely,
our error estimate has an a priori and an a posteriori aspect: With a priori overwhelming probability,
the (random) prefactor can be bounded by a constant that is computable without much further effort, 
on the basis of the given realization in the box of size $L$.

\smallskip

We also rigorously establish that the order of the error estimate in both $L$ and $\ell$ is optimal,
where in this paper we focus on the case of $d=2$.
This amounts to a lower bound on the variance of the quantity of interest when conditioned on
the coefficients inside the computational domain, 
and relies on the deterministic insight that a sensitivity analysis wrt a defect 
commutes with (stochastic) homogenization.
Finally, we carry out numerical experiments that show that this optimal convergence rate already sets in at only moderately large $L$,
and that more naive boundary conditions perform worse both in terms of rate and prefactor.

\end{abstract}


\section{Introduction and main results}

Let the dimension $d\ge 2$ and the ellipticity ratio $\lambda>0$ be fixed.
We will be considering symmetric tensor fields $a$ on $d$-dimensional space $\mathbb{R}^d$
that are uniformly elliptic:
\begin{align}\label{L03}
|\xi|^2\ge \xi\cdot a(x)\xi\ge\lambda|\xi|^2\quad\mbox{for all points $x$ and vectors $\xi$},
\end{align}
where wlog we've set the upper bound to unity. Symmetry is notationally convenient at a few places, but by no
means essential; while we use scalar notation and language, all results hold for systems.
For some localized rhs $g$, say near the origin, 
we are interested in the decaying (ie Lax-Milgram) {\it whole-space} solution $u$ of
\begin{align}\label{L98}
-\nabla\cdot a\nabla u=\nabla\cdot g.
\end{align}
More precisely, we are interested in $-\nabla u$ near the origin, say,
at the origin: $-\nabla u(0)$. In the language of electrostatics, we
are interested in the electric field generated by the neutral and
localized charge distribution $\nabla\cdot g$.  We pose the question
to which precision $-\nabla u(0)$ can be inferred without solving a
PDE in whole-space. Let us denote by $\ell$ the (linear) size of the
support of $g$.  In case of {\it constant} coefficients $a_h$, the
explicit fundamental function $G_h$ allows to reduce the determination
of $-\nabla u(0)$ to the evaluation of an integral over $B_\ell$, the
centered ball of radius $\ell$.  

\medskip

In our case of variable coefficients, we ask the
question of whether one can do better than solving a boundary value
problem with homogeneous boundary data, say the Dirichlet problem
\begin{align}\label{L88}
-\nabla\cdot a\nabla u_0=\nabla\cdot g\;\mbox{in}\;Q_L,\quad
u_0=0\;\mbox{on}\;\partial Q_L
\end{align}
on the centered cube $Q_L:=(-L,L)^d$ for some large scale $L\ge\ell$ (where we take cubes instead
of balls for computational convenience). Under the assumption
that $\ell$ is the only scale of $g$ in the sense that there exists a function $\hat g$ such that
\begin{align}\label{L97}
g(x)=\hat g(\frac{x}{\ell})\quad
\mbox{with}\;\hat g\;\mbox{supported in}\; B_1\;\mbox{and}\;
\hat\nabla\hat g\;\mbox{H\"older continuous},
\end{align}
one expects and we experimentally show in Section \ref{Num}
that the approximation \eqref{L88} is no better than what generically holds in the constant-coefficient case, namely
$\nabla(u_0-u)(0)=O((\frac{\ell}{L})^d)$.
More precisely, we ask the question whether one can do better {\it without} knowing the coefficients 
outside $Q_{2L}$, which is hopeless when $a$ has no further structure. In this paper, we thus consider the
case when $a$ comes from a stationary finite-range ensemble $\langle\cdot\rangle$ of uniformly elliptic coefficient fields;
wlog we assume the range to be unity. We recall that stationary means that $a$ and its shifted version $a(\cdot+z)$ for
any shift vector $z\in\mathbb{R}^d$ have the same distribution under $\langle\cdot\rangle$;
unit range means that for two subset $D,D'\subset\mathbb{R}^d$ with ${\rm dist}(D,D')>1$,
the restrictions $a_{|D}$ and $a_{|D'}$ of the coefficient field $a$ are independent under $\langle\cdot\rangle$.


\subsection{Approximation algorithm and error bound}

Loosely speaking, our first main result Theorem \ref{T1} states that there is an algorithm, the outcome of which is $u^{(L)}$,
that
\begin{itemize}
\item only involves knowing the realization $a$ restricted to $Q_{2L}$,
\item next to the solution of a Dirichlet problem on $Q_L$ only requires
the solution of $2$ (in $d=2$) respectively $9$ (in $d=3$) further Dirichlet problems on $Q_{2L}$, 
\item improves upon (\ref{L88}) by (almost) a factor $L^{-\frac{d}{2}}$, 
albeit with a random prefactor,
\item with overwhelming probability in $L\gg 1$, this prefactor is
  dominated by a constant that can be computed at the cost of further
  $2$ (for $d=2$) respectively $27$ (for $d=3$) (constant-coefficient) Dirichlet problems on
  $Q_{2L}$.
\end{itemize}
Note that Theorem \ref{T1} is a mixture of an a priori result, namely the
probabilistic estimate on when the approach is successful at all, and
an a posteriori result, namely the domination of the prefactor by the computable quantity $r_*^{(L)}$.
Loosely speaking, the second main result, Theorem \ref{T2}, states that in terms of
scaling (both in $L$ and $\ell$), there is no better algorithm since
for a relevant class of ensembles $\langle\cdot\rangle$, the square
root of the variance of $\nabla u(0)$ conditioned on $a_{|Q_{2L}}$ is
of order $(\frac{\ell}{L})^dL^{-\frac{d}{2}}$.  The argument shows
that the factor $L^{-\frac{d}{2}}$ is of CLT-type and loosely speaking
arises as the inverse of the square root of the volume of the
neighboring ``annulus'' $Q_{4L}-Q_{2L}$.

\medskip

This paper only discusses the algorithm that gives the (near) optimal result in case of $d=2$;
the optimal algorithm for the case of $d=3$ (and higher $d$ in general)
would require a refinement, namely the second-order corrector, but no new concepts. 
More precisely, the theory of dipoles developed in \cite{BellaGiuntiOttoPCMINotes}
would have to be replaced by its systematic generalization to multipoles (in particular
quadrupoles) developed in \cite{BellaGiuntiOttoarXiv}, relying on second-order correctors. 
It is the exponent $\beta$ 
on $\frac{1}{L}$ in (\ref{F01}), which can be taken arbitrarily close to $1$ ($=\frac{d}{2}$ for $d=2$), 
that provides the near-optimal CLT improvement over the homogeneous boundary value problem (\ref{L88}).

\medskip

Since we do not assume (H\"older-) continuity of the realization $a$, we do not have access
to the point evaluation (at the origin) of the gradient. Hence in all statements, pointwise control
is replaced by $L^2$-control over an order-one ball.

\medskip

{\bf Correctors $\phi_i$}. Not surprisingly, our algorithm makes use of the correctors $\phi_i$, which
for every coordinate direction $i=1,\cdots,d$ provide $a$-harmonic coordinates $x_i+\phi_i$ by satisfying
\begin{align}\label{L31}
-\nabla\cdot a(e_i+\nabla\phi_i)=0,
\end{align}
where $e_i$ denotes the unit vector in direction $i$.
According to the classical qualitative theory of stochastic homogenization (by ``qualitative'' theory 
we mean the one only relying on ergodicity and stationarity of the ensemble $\langle\cdot\rangle$), 
for almost every realization $a$, a corrector $\phi_i$ of sublinear growth, that is,
\begin{align}\label{L93}
\lim_{r\uparrow\infty}\frac{1}{r}\big(\av_{B_r}(\phi_i-\av_{B_r}\phi_i)^2\big)^\frac{1}{2}=0
\end{align}
can be constructed. Here and in the sequel, $\av_{B_r}$ denotes the average over the ball $B_r$ (of radius $r$ centered at the origin).
Moreover, again almost surely, the homogenized coefficients $a_h$ may be inferred from
\begin{align}\label{L91}
a_he_i=\lim_{L\uparrow\infty}\av_{B_L}q_i\quad\mbox{where}\quad q_i:=a(e_i+\nabla\phi_i).
\end{align}
Hence a naive guess would be to replace the approximation $u_0$ defined through (\ref{L88}) by
the approximation $u_I$ defined through
\begin{align}\label{L89}
-\nabla\cdot a\nabla u_I=\nabla\cdot g\;\mbox{in}\;Q_L,\quad
u_I=\tilde u_h\;\mbox{on}\;\partial Q_L,
\end{align}
where $\tilde u_h$ is the decaying solution of the homogenized problem in the whole space, that is,
\begin{align}\label{L24}
-\nabla\cdot a_h\nabla \tilde u_h=\nabla\cdot g.
\end{align}
This can be seen to generically yield no improved scaling of the error
in $L$ (ie at fixed $\ell$); it only improves the scaling of the error
when $\ell\gg 1$.  Incidentally, it would not fall in the class of the
algorithms we consider, since inferring the homogenized coefficient
$a_h$ requires solving the whole-space problem (\ref{L31}).  In the
context of multiscale method, this is the approach taken in
\cite{OdenVemaganti} (with additional steps to approximate $a_h$).

\medskip

In both periodic and random homogenization, it is known that the so-called two-scale expansion
\begin{align*}
(1+\phi_i\partial_i)\tilde u_h,
\end{align*}
where we use Einstein's convention of summation over repeated indices,
provides a better approximation to $u$ than $\tilde u_h$ itself; in particular,
this approximation is necessary to get closeness of the gradients (in the regime $\ell\gg 1$).
Hence a second attempt would be to replace the approximation $u_I$ defined through (\ref{L89}) by
the approximation $u_{II}$ defined through
\begin{align}\label{L90}
-\nabla\cdot a\nabla u_{II}=\nabla\cdot g\;\mbox{in}\;Q_L,\quad
u_{II}=(1+\phi_i\partial_i)\tilde u_h\;\mbox{on}\;\partial Q_L,
\end{align}
As our numerical experiments in Section \ref{Num} show, this generically yields no improved scaling of the error in $L$.

\medskip

{\bf Dipoles}. The problem with all three approaches, (\ref{L88}), (\ref{L89}), and (\ref{L90}),
is that as soon as $\ell\ll L$,
the far field of $\tilde u_h$ generically has the wrong dipole behavior. This phenomenon
was observed in \cite{BellaGiuntiOttoPCMINotes}, where the right-hand side $g$ in (\ref{L24})
was replaced by $g_i(e_i+\partial_i\phi)$ in order for the gradient of 
the two scale-expansion $\nabla(1+\phi_i\partial_i)\tilde u_h$ to be $O( (\tfrac{\ell}{L})^d (\tfrac{1}{L})^{\beta})$-close to $\nabla u$. 
This is the right strategy for concentrated $g$, ie for $\ell\sim 1$. In order
to also treat a more spread rhs, ie $\ell\gg 1$, we hold on to $\tilde u_h$ but correct it through
an $a_h$-dipole $\delta u_h$ coming from the first moments of $g_i\partial_i\phi$. In formulas,
we pass to 
\begin{align}\label{L18}
u_h:=\tilde u_h+(\int\nabla\phi_i\cdot g)\partial_iG_h,
\end{align}
where $G_h$ is the fundamental solution for $-\nabla\cdot a_h\nabla$.
Hence our more educated ansatz is to replace the approximation $u_{II}$ defined through (\ref{L90}) by
the approximation $u_{III}$ defined through
\begin{align}\label{L94}
-\nabla\cdot a\nabla u_{III}=\nabla\cdot g\;\mbox{in}\;Q_L,\quad
u_{III}=(1+\phi_i\partial_i)u_h\;\mbox{on}\;\partial Q_L.
\end{align}
Corollary \ref{L2} shows that this approximation indeed reduces the (generic) error of (\ref{L88})
by a factor with the desired $L$-scaling $\frac{1}{L^\beta}$. Corollary \ref{L2} relies on Lemma \ref{L1},
which is a minor modification of \cite{BellaGiuntiOttoPCMINotes}.

\medskip

{\bf Flux correctors $\sigma_i$}. As can be seen from Lemma \ref{L1}, the prefactor comes in form of $r_*^\beta$ for some 
length scale $r_*$. This length scale has the interpretation that for larger scales $r\ge r_*$,
the quantified sublinear growth of the correctors $\phi_i$ sets in, cf (\ref{L04}). 
The sublinear growth is quantified through the exponent $\beta$
(with $\beta$ close to $0$ meaning almost linear growth, and $\beta$ close to $1$ meaning
almost no growth). However, for quantitative results like Lemma \ref{L1}, which closely follows
\cite{BellaGiuntiOttoPCMINotes}, itself inspired by \cite{FischerOttoCPDE}, it is not sufficient
to monitor just $\phi_i$. In fact, the harmonic vector field $e_i+\nabla\phi_i$ (the electric field
in the language of electrostatics) is not just a closed $1$-form; but through the flux 
$q_i:=a(e_i+\nabla\phi_i)$ (the electric current in the language of electrostatics) it provides a 
closed $(d-1)$-form. Hence there is not just the $0$-form (a scalar potential) $x_i+\phi_i$, 
or rather its correction $\phi_i$,
but there naturally is also a $(d-2)$-form (a vector potential in the 3-d language, or a stream function
in the 2-d language), or rather its correction $\sigma_i$, which we can write as a skew-symmetric
tensor field $\sigma_i=\{\sigma_{ijk}\}_{j,k=1,\cdots,d}$. In view of (\ref{L91}),
this correction should satisfy
\begin{align}\label{L01}
q_i=a_he_i+\nabla\cdot\sigma_i,
\end{align}
where $(\nabla\cdot\sigma_i)_j:=\partial_k\sigma_{ijk}$. Note that by skew symmetry of $\sigma_i$ we have $\nabla\cdot\nabla\cdot\sigma_i=0$
so that (\ref{L01}) contains the familiar (\ref{L31}), as it implies $\nabla \cdot q_i = 0$.
Clearly, (\ref{L01}) determines $\sigma_i$ only up to a $(d-3)$-form, ie the freedom of the choice of a gauge.
A particularly simple choice of gauge is
\begin{align}\label{L30}
-\Delta\sigma_{ijk}=\partial_jq_{ik}-\partial_kq_{ij}.
\end{align}
This skew-symmetric field is not uncommon in periodic homogenization; in 
qualitative stochastic homogenization \cite[Lemma 1, Corollary 1]{GNOarXiv} 
it has been shown to almost surely exist with sublinear growth:
\begin{align}\label{L92}
\lim_{r\uparrow\infty}\frac{1}{r}\big(\av_{B_r}|\sigma_i-\av_{B_r}\sigma_i|^2\big)^\frac{1}{2}=0.
\end{align} 

\medskip

{\bf Radius $r_*$}. Loosely speaking, starting from the length scales $r$
at which the lhs expression (\ref{L93}) and (\ref{L92}) drop below a threshold only depending on $d$
and $\lambda$, the operator $-\nabla\cdot a\nabla$ inherits the
regularity theory of $-\nabla\cdot a_h\nabla$, both for Schauder
theory on the $C^{1,\alpha}$ level (then the threshold depends in
addition on $0<\alpha<1$) \cite[Corollary 3]{GNOarXiv}, and for the
Calderon-Zygmund theory on the $\dot H^{1,p}$ level (then the threshold
depends in addition on $1<p<\infty$) \cite[Corollary 4]{GNOarXiv}. This type of
theory had been developed by Avellaneda \& Lin \cite{AL} for the
periodic case; Armstrong \& Smart \cite{AS} were first to extent this
to the random case.  As mentioned, $r_*$ quantifies (\ref{L92}) and
(\ref{L93}); it is defined as the length scale starting from which
$\beta$-sublinear growth of the scalar and vector correctors kick in,
that is,
\begin{align}\label{L04}
\frac{1}{r}\big(\av_{B_r}|(\phi,\sigma)-\av_{B_r}(\phi,\sigma)|^2\big)^\frac{1}{2}\le(\frac{r_*}{r})^\beta
\quad\mbox{for all}\;r\ge r_*,
\end{align}
where $(\phi,\sigma)$ is the collection of all components $\{\phi_i,\sigma_{ijk}\}_{i,j,k}$.
While the form (\ref{L04}) is more natural, it is only seemingly weaker than
\begin{align}\label{L32}
\frac{1}{r}\big(\av_{B_r}|(\phi,\sigma)-\av_{B_{r_*}}(\phi,\sigma)|^2\big)^\frac{1}{2}\le
C(\beta)(\frac{r_*}{r})^\beta
\quad\mbox{for all}\;r\ge r_*,
\end{align}
as we shall show at the beginning of the proof of Lemma \ref{L1}. 
In \cite[Theorem 1 ii)]{FischerOttoSPDE} it is shown that $r_*$ satisfies (optimal) stretched exponential bounds
even under weak correlation decay of $\langle\cdot\rangle$.

\medskip

Clearly, this notion of $r_*$ singles out the origin (it can and will be defined for other bounds
$r_*(y)$, making it a stationary random field). The origin plays a special role in our analysis in 
two ways: It is where we want to monitor the error (on the level of the gradient) and where
the rhs $g$ is concentrated (on scale $\ell$). In view of the above-mentioned
$C^{1,\alpha}$-regularity theory that kicks in (only) from scales $r_*$ onwards, it is not
surprising that we can localize the error (on the level of the gradients) only to scales
$R\ge r_*$, see Proposition \ref{P1} (which in Theorem \ref{T1} is expressed in terms of the proxy $r_*^{(L)}$).
For the same reason, we need the condition $\ell\ge r_*$ in Proposition \ref{P1}. 

\medskip

{\bf Algorithm}. The ``algorithm'' (\ref{L94}) is not admissible, since it involves solving 
the $d$ whole-space problems (\ref{L31}).
The natural idea is to replace these whole-space problems by Dirichlet problems on $Q_{2L}$;
for reasons that will become clearer later, we do the same for the $d\times\frac{d(d-1)}{2}$
whole-space problems (\ref{L30}) (while having constant coefficients they feature an extended supported rhs). 
For any coordinate direction $i=1,\cdots,d$, let the function $\phiL_i$ 
and the skew-symmetric tensor field $\sigmaL_i=\{\sigmaL_{ijk}\}_{j,k=1,\cdots,d}$ 
be determined through
\begin{align}
-\nabla\cdot a(e_i+\nabla\phiL_i)=0\;\mbox{in}\;Q_{2L},&\quad
\phiL_i=0\;\mbox{on}\;\partial Q_{2L},\label{L95}\\
-\Delta\sigmaL_{ijk}=\partial_j\qL_{ik}-\partial_k\qL_{ij}\;\mbox{on}\;Q_{2L},&\quad
\sigmaL_{ijk}=0\;\mbox{on}\;\partial Q_{2L},\label{L96}
\end{align}
where we have set for abbreviation $\qL_i:=a(e_i+\nabla\phiL_i)$. While
an easy calculation shows that (\ref{L95}) \& (\ref{L96}) imply
$\Delta(q_i-\nabla\cdot\sigma_i)=0$, this does in general not yield
$\qL_i-\av_{Q_{2L}}\qL_i$ $=\nabla\cdot\sigmaL_i$. The latter would be
automatic in case of periodic boundary conditions, in which case we
would replace the homogenized coefficient $a_h$, cf (\ref{L91}),
by $\aL_h e_i$ $=\av_{Q_{2L}}\qL_i$. In our (more ambitious) case of
Dirichlet boundary conditions, we pick a mask $\hat\omega$ of an
averaging function with
\begin{align}\label{F09}
\omega(x)=\frac{1}{L^d}\hat \omega(\frac{x}{L})\quad
\mbox{with}\;\hat\omega\;\mbox{supported in}\; Q_1,\;\int\hat\omega=1,\;\text{and}\;
\hat\nabla\hat\omega\;\mbox{bounded},
\end{align}
and set
\begin{align}\label{L99}
\aL_he_i:=\int\omega \qL_i;
\end{align}
it is a consequence of (\ref{L47}) in Lemma \ref{L3} that $\aL_h$ is elliptic.
We now make the corresponding changes on the level of $\tilde u_h$ and $u_h$:
We substitute $\tilde u_h$ defined in (\ref{L24}) by the decaying solution $\tuL_h$ of
\begin{align}\label{L48}
-\nabla\cdot \aL_h\nabla \tuL_h=\nabla\cdot g
\end{align}
and $u_h$ defined in (\ref{L18}) by $\uL_h$ defined through
\begin{align}\label{L52}
\uL_h:=\tuL_h+(\int\nabla\phiL_i\cdot g)\partial_i\GL_h,
\end{align}
where $\GL_h$ denotes the fundamental solution of $-\nabla\cdot \aL_h\nabla$.

\medskip

As mentioned above, Theorem \ref{T1} is a mixture of an a priori and
an a posteriori result: Through the scale $L_0$, which only depends on the dimension $d<\infty$, 
the ellipticity ratio $\lambda>0$, the sublinear growth exponent $\beta<1$, 
and the stretched exponential exponent $s<2(1-\beta)$,
Theorem \ref{T1} provides an {\it a priori} estimate on the
probability that the random pick of a realization $a$ is so bad that
the algorithm fails. Through (\ref{F02}), which characterizes the
scale $\rL_*$ in a computable fashion, it provides an a posteriori
estimate on the constant $(\rL_*)^\beta$ in the
$\frac{1}{L^\beta}$-improvement of the error estimate (\ref{F01}).  We
think of (\ref{F01}) as an {\it a posteriori} estimate, since it
relies on an auxiliary computation based on the given realization $a$.


\begin{theorem}\label{T1}
Let $\langle\cdot\rangle$ be a stationary ensemble of uniformly elliptic coefficient fields $a$,
cf (\ref{L03}), that is of unit range. 
Then for any exponents $\beta\in(0,1)$ and $s\in(0,2(1-\beta))$, there exists a scale $L_0=L_0(d,\lambda,\beta,s)$
so that for any scale $L\ge L_0$, with probability $1-\exp(-(\frac{L}{L_0})^s)$ 
a realization $a$ has the following property: 

\smallskip

Let $\{\phiL_i\}_{i=1,\cdots,d}$, $\sigmaL_i=\{\sigmaL_{ijk}\}_{i,j,k=1,\cdots,d}$,
and $\aL_h$ be defined through (\ref{L95}), (\ref{L96}), and (\ref{L99}).
Suppose that $(\phiL,\sigmaL)$ is quantitatively sublinear 
in the sense of that there exists a scale $1\le \rL_*\le L$ with
\begin{align}\label{F02}
\frac{1}{r}\big(\av_{Q_{r}}|(\phiL,\sigmaL)|^2\big)^\frac{1}{2}
\le(\frac{\rL_*}{r})^\beta\quad\mbox{for}\;2L\ge r\ge \rL_*.
\end{align}

\smallskip

Let $g$ be of the form (\ref{L97}) for some scale $\ell\in[\rL_*,L]$ and mask $\hat g$,
and let $u$ be the decaying solution of (\ref{L98}).
Consider the solution $\uL$ of
\begin{align}\label{L56}
-\nabla\cdot a\nabla \uL=\nabla\cdot g\;\mbox{in}\;Q_L,
\quad \uL=(1+\phiL_i\partial_i)\uL_h\;\mbox{on}\;\partial Q_L,
\end{align}
where $\uL_h$ is defined through  (\ref{L48}) and (\ref{L52}).
Then we have 
\begin{align}\label{F01}
\big(\av_{B_R}|\nabla(\uL-u)|^2\Big)^\frac{1}{2}
\le C(\frac{\ell}{L})^d(\frac{\rL_*}{L})^{\beta}\quad\mbox{for}\;L\ge R\ge \rL_*
\end{align}
with a constant $C$ that only depends on $d$, $\lambda$, $\beta$, and the H\"older norm
of $\hat\nabla\hat g$ as well as the supremum norm of $\hat\nabla\hat\omega$.
\end{theorem}


Theorem \ref{T1} has three ingredients, the deterministic a priori
error estimate provided by Proposition \ref{P1}, the stochastic
ingredient Lemma \ref{L6}, and the deterministic Lemma \ref{L5}, that
allows to pass from an a priori to an a posteriori error estimate. We
say that Proposition \ref{P1} provides an {\it a priori} and {\it
  deterministic} error estimate since it is formulated in terms of
$r_*$ characterizing the sublinear growth of the augmented corrector
$(\phi,\sigma)$, cf (\ref{L04}).  In addition, it starts from a given
uniformly elliptic coefficient field $a$, which might but does not have
to be a realization under $\langle\cdot\rangle$.  The only assumption
on the Dirichlet proxy $(\phiL,\sigmaL)$ is that it is well-behaved on
the large scale $2L$, cf (\ref{L70}), but not necessarily on smaller
scales as in (\ref{L32}).

\begin{proposition}\label{P1}
Consider a given $\lambda$-uniformly elliptic coefficient field $a$ on $\mathbb{R}^d$,
cf (\ref{L03}). Suppose that there exists a tensor $a_h$ and, for $i=1,\cdots,d$, a scalar field 
$\phi_i$ and a skew-symmetric tensor field $\sigma_i$ such that (\ref{L01}) holds.
Suppose that for given $\beta\in(0,1)$ there exists a radius $r_*$ such that (\ref{L04}) holds.

\smallskip

For given $L\ge r_*$, let $\phiL_i$, $\sigmaL_i$, and $\aL_h$ be defined through
(\ref{L96}), (\ref{L95}), and (\ref{L99}). We assume that
\begin{align}\label{L70}
\frac{1}{2L}\big(\av_{Q_{2L}}|(\phiL,\sigmaL)|^2\big)^\frac{1}{2}\le(\frac{r_*}{2L})^\beta.
\end{align}

\smallskip

Given $\hat g$ and $\ell\in[r_*,L]$ let $u$ and $\uL$ be defined as in Theorem \ref{T1}.
Then we have 
\begin{align}\label{L57}
\big(\av_{B_R}|\nabla(\uL-u)|^2\Big)^\frac{1}{2}
\le C(\frac{\ell}{L})^d(\frac{r_*}{L})^{\beta}\quad\mbox{for}\;L\ge R\ge r_*,
\end{align}
where the constant $C$ is of the same type as in Theorem \ref{T1}.
\end{proposition}


The following lemma is the key ingredient for Proposition \ref{P1},
it shows that indeed the multipole has to be corrected in the sense of
(\ref{L18}). Its proof essentially follows \cite[Theorem 0.2]{BellaGiuntiOttoPCMINotes}.

\begin{lemma}\label{L1}
Let $a$, $a_h$, $\phi_i$, $\sigma_i$, $\beta$, and $r_*$ be as in Proposition \ref{P1}.
Let $g$ be of the form (\ref{L97}) for some $\ell\ge r_*$ and $\hat g$,
and let $u$ be the decaying solution of (\ref{L98}).
Let $u_h$ be defined through  (\ref{L24}) and (\ref{L18}).
Then we have
\begin{align}\label{L17}
\big(\av_{B_R^c}|\nabla(u-(1+\phi_i\partial_i)u_h)|^2\Big)^\frac{1}{2}
\le C(\frac{\ell}{R})^d(\frac{r_*}{R})^{\beta}\quad\mbox{for}\;R\ge r_*,
\end{align}
where $\phi_i$ is normalized through
\begin{align}\label{L62}
\av_{B_{r_*}}\phi_i=0.
\end{align}
Here $C$ is a constant of the same type as in Theorem \ref{T1}. Furthermore
$\av_{B_R^c}$ is the abbreviation of $\frac{1}{R^d}\int_{B_R^c}$, where $B_R^c$
is the complement of $B_R$.

\end{lemma}


Equipped with Lemma \ref{L1}, we may assess the effect of a computational domain $Q_{L}$
endowed with the Dirichlet conditions given by $(1+\phi_i\partial_i)u_h$, cf (\ref{L94}).

\begin{corollary}\label{L2}
Let $a$, $a_h$, $\phi_i$, $\sigma_{i}$, $\beta$, and $r_*$ be as in Proposition \ref{P1}.
Let $g$ be of the form (\ref{L97}) for some $\ell\ge r_*$ and $\hat g$,
and let $u$ be the decaying solution of (\ref{L98}).

\smallskip

For $L\ge\ell$ let $u_{III}$ be the solution of the Dirichlet problem (\ref{L94}).
Then we have
\begin{align}\label{L26}
\big(\av_{B_R}|\nabla(u_{III}-u)|^2\Big)^\frac{1}{2}
\le C(\frac{\ell}{L})^d(\frac{r_*}{L})^{\beta}\quad\mbox{for}\;L\ge R\ge r_*,
\end{align}
where $C$ is of the type as in Theorem \ref{T1}.
\end{corollary}


In order to pass from Corollary \ref{L2} to Proposition \ref{P1}, that is, from $u_{III}$ to $u^{(L)}$,
we need to replace $(\phi,\sigma)$ by the computable $(\phiL,\sigmaL)$
in order to pass from $u_{III}$ to $\uL$.

\begin{lemma}\label{L3}
Let $a$, $a_h$, $\phi_i$, and $\sigma_i$ be as in Proposition \ref{P1}.
Suppose that for given $\beta\in(0,1)$ there exists a radius $r_*$ such that (\ref{L04}) holds in the weaker
form of
\begin{align}\label{F19}
\frac{1}{2L}\big(\av_{Q_{2L}}|(\phi,\sigma)-\av_{Q_{2L}}(\phi,\sigma)|^2\big)^\frac{1}{2}\le(\frac{r_*}{2L})^\beta,
\end{align}
for given $L\ge r_*$. Let $\phiL_i$, $\sigmaL_i$, and $\aL_h$ be defined through
(\ref{L95}), (\ref{L96}), and (\ref{L99}). Suppose that (\ref{L70}) holds.
Then we have 
\begin{align}
\big(\av_{Q_{\frac{3}{2}L}}|\nabla(\phiL-\phi)|^2\big)^\frac{1}{2}&\le C(\frac{r_*}{L})^{\beta},\label{L44}\\
\big(\av_{Q_{L}}|\nabla(\sigmaL-\sigma)|^2\big)^\frac{1}{2}&\le C(\frac{r_*}{L})^{\beta},\label{F20}\\
\mbox{and}\quad|\aL_h-a_h|&\le C(\frac{r_*}{L})^{\beta},\label{L47}
\end{align}
where $C$ is of the same type as in Theorem \ref{T1}.
\end{lemma}


As mentioned, the following Lemma \ref{L5} allows to pass from the deterministic {\it a priori} estimate
of Proposition \ref{P1} to the deterministic {\it a posteriori} estimate of Theorem \ref{T1}.
It shows that sublinear growth of $(\phi,\sigma)$ on scales $R\ge L$ with (some) pre-factor $(\rL_*)^\beta$, cf (\ref{L80}),
and sublinear growth of the Dirichlet proxy $(\phiL,\sigmaL)$ on scales $2L\ge r\ge \rL_*$, cf (\ref{L73}),
implies sublinear growth of $(\phi,\sigma)$ on all scales $r\ge r_*^{(L)}$, that is, $r_*\lesssim \rL_*$.

\begin{lemma}\label{L5}
Let $a$, $a_h$, $\phi_i$, and $\sigma_i$ be as in Proposition \ref{P1}.
For given $L\ge r_*$, let $\aL_h$, $\phiL_i$, and $\sigmaL_i$ be defined through
(\ref{L99}), (\ref{L95}), and (\ref{L96}). We assume that for a given $\beta\in(0,1)$
there exists a scale $\rL_*\le L$ such that we have on the one side
\begin{align}\label{L80}
\frac{1}{R}\big(\av_{B_R}|(\phi,\sigma)-\av_{B_R}(\phi,\sigma)|^2\big)^\frac{1}{2}\le(\frac{\rL_*}{R})^\beta
\quad\mbox{for}\;R\ge L
\end{align}
and on the other side
\begin{align}\label{L73}
\frac{1}{r}\big(\av_{B_r}|(\phiL,\sigmaL)|^2\big)^\frac{1}{2}\le(\frac{\rL_*}{r})^\beta
\quad\mbox{for}\; 2L\ge r\ge \rL_*.
\end{align}
Then we have
\begin{align}\label{L85}
\frac{1}{r}\big(\av_{B_r}|(\phi,\sigma)-\av_{B_r}(\phi,\sigma)|^2\big)^\frac{1}{2}\le C(\frac{\rL_*}{r})^\beta
\quad\mbox{for}\;r\ge \rL_*
\end{align}
with the constant $C$ being of the same type as in Theorem \ref{T1}.
\end{lemma}


The only stochastic ingredient, which we ``take from the shelf'',
for Theorem \ref{T1} is part ii) of the following lemma;
part i) is needed in the argument for the lower bound, Theorem \ref{T2} below. Lemma \ref{L6} provides
stochastic bounds on the (augmented) corrector $(\phi,\sigma)$, and essentially amounts to
saying that it is bounded with overwhelming probability in $d>2$ and almost so in $d=2$
(the form of the statement of Lemma \ref{L6} is marginally weakened by not distinguishing
the cases $d=2$ and $d>2$).
By now, there are several approaches to such a result: The (historically) first result of this type is
\cite[Proposition 2.1]{GloriaOttoAP},
and is based on functional inequalities (at first in case of a discrete medium;
see \cite[Proposition 1]{GloriaOttoJEMS} for an extension to the continuum case)
and thus (indirectly) relies on an underlying product structure of the ensemble $\langle\cdot\rangle$
(eg a Gaussian field or Poisson point process like in Definition \ref{D1} below).
Here, we follow a more recent approach that is based on a finite, say unit, range assumption.
There are two possible references for this second approach: \cite{AKM} based
on the variational approach of \cite{AS}, and \cite{GloriaOttofiniterange}, based
on a semi-group approach. For reasons of familiarity, we opt for the second reference,
which also has the advantage of treating $\sigma$ next to $\phi$.

\begin{lemma}\label{L6}
Let $\langle\cdot\rangle$ be a stationary ensemble of uniformly elliptic coefficient fields, cf (\ref{L03}),
that has unit range of dependence. Then for every $i=1,\cdots,d$, there exist a (random) scalar field $\phi_i$
and a (random) skew symmetric tensor field $\sigma_i$ such the gradient fields $\nabla\phi_i$ and $\nabla\sigma_i$
are stationary, of finite second moments, and of vanishing expectation, and such that (\ref{L01}) (and thus (\ref{L31}))
and (\ref{L30}) hold. Moreover,
\begin{itemize}
\item[i)] For every exponent $0<\beta<1$ there exists a (random) radius $r_*$ such that (\ref{L04}) holds
and which satisfies
\begin{align}\label{io06}
\langle\exp(r_*^s)\rangle\le C
\end{align}
for any exponent $0<s<2$ with $C=C(d,\lambda,\beta,s)$.
\item[ii)] For every exponents $0<\beta<1$ and $0<s<2(1-\beta)$, there exists
a scale $L_0=L_0(d,\lambda,\beta,s)$ such that for every $L\ge L_0$ the statement
\begin{align*}
\frac{1}{R}\big(\av_{B_R}|(\phi,\sigma)-\av_{B_R}(\phi,\sigma)|^2\big)^\frac{1}{2}\le(\frac{1}{R})^\beta
\quad\mbox{for all}\;R\ge L
\end{align*}
fails with probability $\le \exp(-(\frac{L}{L_0})^s)$.
\end{itemize}
\end{lemma}


\subsection{The lower bound}

Our second main result states that the fluctuations of $\nabla u(0)$, when conditioned on the
coefficient field inside the ball $B_L$ are at least of the order of $(\frac{\ell}{L})^d(\frac{1}{L})^\frac{d}{2}$. 
The first factor is a (deterministic) consequence of fact that $\nabla u(0)$, in view of the rhs supported in $B_\ell$,
is not very sensitive in the coefficients in $B_L^c$. 
The second factor scales as the inverse of the square root of the volume 
of the annulus $B_{2L}-B_L$, and thus has a CLT-flavor to it. More precisely, instead of $\nabla u(0)$ we monitor
smooth averages of $\nabla u$ on a sufficiently large but order-one scale $R$ near the origin. 
Similarly to (\ref{F09}) we fix a (universal) averaging function on $B_{R}$ through
\begin{align}\label{g18}
\omega(x)=\frac{1}{R^d}\hat \omega(\frac{x}{R})\quad\mbox{with}\quad
\hat\omega=0\;\mbox{in}\;B_1^c,\quad\int\hat\omega=1,\quad
\hat\nabla\hat\omega\;\mbox{bounded}
\end{align}
and consider $\int\omega\nabla u$.
We establish this lower bound on the fluctuations under convenient assumptions on the ensemble: 

\begin{definition}\label{D1}
Let $\langle\cdot\rangle$ denote the distribution of the Poisson point process $X$ in $\mathbb{R}^d$ of
unit intensity. We assume that there exists a measurable map $A$ from the space of point configurations into
the space of $\lambda$-uniformly elliptic coefficient fields; in other words, we
consider the ensemble of such coefficient fields given by
\begin{align*}
a(x)=A(x,\cdot)
\end{align*}
and thus make the following assumptions on $A$:
For all points $x\in\mathbb{R}^d$, point configurations $X,X'\subset\mathbb{R}^d$, shift vectors $z\in\mathbb{R}^d$
we impose
\begin{itemize}
\item shift-invariance, that is,
\begin{align}\label{l22}
A(z+x,z+X)=A(x,X),
\end{align}
\item locality, that is,
\begin{align}\label{l06}
A(x,X)=A(x,X')\quad\mbox{provided}\;X\cap B_1(x)=X'\cap B_1(x),
\end{align}
\item monotonicity, that is,
\begin{align}\label{l18}
A(x,\emptyset)<a_h\quad\mbox{as inequality between symmetric matrices}.
\end{align}
\end{itemize}
In particular, $\langle\cdot\rangle$ is stationary and of range unity, so that by Lemma \ref{L6} i), for every $\beta<1$
and $s<2$, there exists $r_*$ with (\ref{L04}) and such that
\begin{align}\label{l13}
\langle \exp(r_*^s)\rangle\le C(d,\lambda,\beta,s).
\end{align}
\end{definition}


\begin{theorem}\label{T2}
Let the ensemble be as in Definition \ref{D1}.
Consider the solution $u$ of (\ref{L98}) with rhs of the form (\ref{L97}) for a given $\hat g$.
Then there exists a radius $R$ such for all scales $L$, $\ell$ with $C\le\ell\le\frac{1}{C}L$ we have
\begin{align}\label{j17}
\big\langle\big|\int\omega\nabla u-\langle\int\omega\nabla u|B_L\rangle\big|^2\big\rangle^\frac{1}{2}
\ge\frac{1}{C}(\frac{\ell}{L})^d(\frac{1}{L})^\frac{d}{2}|\int\hat g|,
\end{align}
where $\omega$ is defined as in (\ref{g18}). Here the radius $R$ and the constant $C$ depend on the ensemble,
on the sup norm of $\hat\nabla\hat\omega$,
and on the H\"older norm of $\hat\nabla\hat g$.
\end{theorem}


Like Theorem \ref{T1}, Theorem \ref{T2} relies on a purely deterministic result, namely Proposition \ref{P2} which is
of independent interest. Proposition \ref{P2} monitors the effect of a ``defect'' in the medium $a$; 
more precisely, we consider the medium $a'$ given by
\begin{align}\label{g10}
a'=a_0\quad\mbox{in}\;B_R(y)\quad\mbox{and}\quad a'=a\quad\mbox{outside of}\;B_R(y),
\end{align}
where $a_0$ is some other $\lambda$-uniformly elliptic coefficient field, cf (\ref{L03}),
and $y$ is some point and $R$ some radius. We are interested in the effect on the solution
$u$ of our whole-space problem (\ref{L98}) with localized rhs $g$.
Hence we compare $u$ to $u'$ given by the decaying solution of
\begin{align}\label{g11}
-\nabla\cdot a'\nabla u'=\nabla\cdot g.
\end{align}
Here, we think of $B_R(y)$ as being far from the origin where $g$ is localized, cf (\ref{L97}).
Proposition \ref{P2} states that to leading order, the effect of the inclusion is captured by its
effect on the level of the homogenized coefficients. 
More precisely, we consider the coefficient field $a_h'$ given by
\begin{align}\label{g14}
a'_h=a_0\quad\mbox{in}\;B_R(y)\quad\mbox{and}\quad a'_h=a_h\quad\mbox{outside of}\;B_R(y)
\end{align}
and the decaying solution $u_h'$ of
\begin{align}\label{g12}
-\nabla\cdot a'_h\nabla u'_h=-\nabla\cdot a_h\nabla u_h,
\end{align}
where $u_h$ is given by (\ref{L18}).
Proposition \ref{P2} states that indeed $u'-u\approx u_h'-u_h$.
As for the other results, Proposition \ref{P2} does so on the level
of the gradient (and thus involves the two-scale expansion) and in a localized way.
It is the positivity of the exponent $\alpha$ in (\ref{g15}) that ensures that
indeed to leading order, $\nabla(u'-u)$ behaves like $\partial_i(u_h'-u_h)(e_i+\nabla\phi_i)$, 
since, as a classical argument shows, $\nabla(u'_h-u_h)(0)$ generically scales as
$(\frac{\ell}{|y|})^d(\frac{R}{|y|})^d$. Loosely speaking, Proposition \ref{P2}
states that {\it sensitivity analysis} (the dependence of a solution on the coefficient field) {\it and
homogenization commute}.

\begin{proposition}\label{P2}
Consider a given $\lambda$-uniformly elliptic coefficient field $a$ on $\mathbb{R}^d$,
cf (\ref{L03}). Suppose that there exists a tensor $a_h$ and, for $i=1,\cdots,d$, a scalar field
$\phi_i$ and a skew-symmetric tensor field $\sigma_i$ such that (\ref{L01}) holds.
Suppose that for given $\beta\in(0,1)$ there exists a radius $r_*(0)$ such that (\ref{L04}) holds
and a radius $r_*(y)$ such that (\ref{L04}) holds with the origin replaced by some point $y$.

\smallskip

We are given another $\lambda$-uniformly elliptic coefficient field $a_0$, cf (\ref{L03}),
and consider $a'$ and $a'_h$ given through (\ref{g10}) and (\ref{g14}).
Consider $u$, $u'$, $u_h$, and $u_h'$ defined through (\ref{L98}), (\ref{g11}), (\ref{L18}), and (\ref{g12}).

\smallskip

Under the assumption that $\ell\ge r_*(0)$ and provided the radius $R$ satisfies $\frac{1}{2}|y|\ge R\ge r_*(0),r_*(y)$
we have 
\begin{align}\label{g15}
\big(\av_{B_R}|\nabla\big((u'-u)&-(1+\phi_i\partial_i)(u_h'-u_h)\big)|^2\big)^\frac{1}{2}\nonumber\\
&\le C(\frac{\ell}{|y|})^d(\frac{R}{|y|})^d\big((\frac{R}{|y|})^\beta+(\frac{r_*(y)}{R})^\alpha\big)
\end{align}
with an exponent $\alpha>0$ that only depends on $d$, $\lambda$ and $\beta$
and a constant $C<\infty$ that in addition depends on the H\"older norm of $\hat g$.
\end{proposition}


Proposition \ref{P2} relies on two ingredients: The first ingredient
is Lemma \ref{L9} which to leading order characterizes $u'-u$ in terms
of the homogenized solution $u_h$ and a tensor $\delta a_{ij}$,
defined in (\ref{g55}), that captures how the corrector $\phi$ is
affected by the defect. Since we are interested in $u'-u$ near the
origin, the Green's function $G_h$ of the homogenized operator is also
involved. Since we characterize $u'-u$ on the level of gradients, the
estimate involves the correction $\nabla\phi_k$.  The second
ingredient is Lemma \ref{L8} that establishes a version of Proposition
\ref{P2} with the general solution $u$ of an equation involving
$-\nabla\cdot a\nabla$ replaced by a specific one, namely the
corrector $\phi_i$. With help of this lemma, we establish Corollary
\ref{Le10} that characterizes the tensor $\delta a_{ij}$ in terms of
its counterpart $\delta a_h$.

\medskip

How is the solution of the corrector equation (\ref{L31}) affected by replacing $a$ by $a'$, cf (\ref{g10})?
We denote by $\phi'_i$ the solution of
\begin{align}\label{g16}
-\nabla\cdot a'(e_i+\nabla\phi'_i)=0
\end{align}
that behaves as $\phi_i$ at infinity (in the sense that $\nabla(\phi_i'-\phi_i)$ is square integrable).
Lemma \ref{L8} compares $\phi_i'$ to $\phi_{ih}'$, the decaying solution of
\begin{align}\label{g17}
-\nabla\cdot a'_h(e_i+\nabla\phi_{ih}')=0.
\end{align}
More precisely, Lemma \ref{L8} compares the two solutions on the level of the potentials and of flux averages.
Similarly to (\ref{F09}) we fix a (universal) averaging function on $B_{2R}(y)$ through 
%
\begin{align}\label{g18bis}
\lefteqn{\omegay(x)=\frac{1}{R^d}\hat \omega(\frac{x-y}{R})\quad\mbox{with}}\nonumber\\
&\hat\omega=\frac{1}{2|B_1|}\;\mbox{in}\;B_1,\quad\hat\omega=0\;\mbox{outside}\;B_2,\quad\int\hat\omega=1,\quad
\hat\nabla\hat\omega\;\mbox{bounded}.
\end{align}
%
In terms of arguments, Lemma \ref{L8} takes inspiration from \cite[Proposition 1]{GNOarXiv}.

\begin{lemma}\label{L8}
We make the same assumptions on $a$, $a_h$, $\phi_i$, $\sigma_i$, $r_*(y)$ as in Proposition \ref{P2}.
Let $a_0$, $a'$, and $a_h'$ be as in Proposition \ref{P2}. For $i=1,\cdots,d$ we consider the solutions
$\phi_i'$ and $\phi_{ih}'$ of (\ref{g16}) and (\ref{g17}). Then for $R\ge r_*(y)$ 
and $\omegay$, cf (\ref{g18bis}), we have
\begin{align}
\frac{1}{R}\big(\av_{B_{2R}(y)}(\phi_i'-\phi_{ih}'-\av_{B_{2R}(y)}(\phi_i'-\phi_{ih}'))^2\big)^\frac{1}{2}\le C(\frac{r_*(y)}{R})^\alpha,\label{g36}\\
\big|\int\omegay a'(e_i+\nabla\phi_i')-\int\omegay a_h'(e_i+\nabla\phi_{ih}')\big|
\le C(\frac{r_*(y)}{R})^\alpha,\label{g35}
\end{align}
where $\alpha$ and $C$ only depend on $d$, $\lambda$ and $\beta$.
\end{lemma}

\medskip

As we shall see in Lemma \ref{L9}, it is the tensor 
\begin{align}\label{g55}
\delta a_{ij}:=\int(e_j+\nabla\phi_j')\cdot(a'-a)(e_i+\nabla\phi_i)
\end{align}
that to leading order governs the sensitivity of the solution operator of $-\nabla\cdot a\nabla$ under changing
the coefficients from $a$ to $a'$.  It is an easy consequence of Lemma \ref{L8} that to leading order, we may
replace the medium described by $a$ (in both $a$ itself and in its perturbation $a'$, and the corresponding
correctors) by its homogenized version $a_h$ in the expression (\ref{g55}).

\begin{corollary}\label{Le10}
Under the assumptions of Lemma \ref{L8} we have for $R\ge r_*(y)$
\begin{align}\label{g41}
\lefteqn{\big|\frac{1}{R^d}\int(e_j+\nabla\phi_j')\cdot(a'-a)(e_i+\nabla\phi_i)}\nonumber\\
&-\frac{1}{R^d}\int(e_j+\nabla\phi_{jh}')\cdot(a'_h-a_h)e_i\big|\le C(\frac{r_*(y)}{R})^\alpha,
\end{align}
where $\alpha$ and $C$ are as Lemma \ref{L8}.
\end{corollary}


The following lemma heavily relies on \cite{BellaGiuntiOttoPCMINotes}, indirectly through
expanding on Lemma \ref{L1}, but also more directly for localization.

\begin{lemma}\label{L9}
We make the same assumptions on $a$, $a_h$, $\phi_i$, $\sigma_i$, $r_*$, $r_*(y)$ as in Proposition \ref{P2}.
Let $a_0$, $a'$, $a_h'$, $u$, and $u'$ be as in Proposition \ref{P2}. Then for $R\ge r_*, r_*(y)$ we have
\begin{align}\label{g85}
\big(\av_{B_R}|\nabla(u'-u)&
-\partial_iu_h(y)\delta a_{ij}\partial_j\partial_kG_h(-y)(e_k+\nabla\phi_k)|^2\big)^\frac{1}{2}\nonumber\\
&\le C(\frac{\ell}{|y|})^d(\frac{R}{|y|})^{d+\beta},
\end{align}
where $C$ depends on $d$, $\lambda$, $\beta$ and the H\"older norm of $\hat g$.
\end{lemma}


Equipped with Lemma \ref{L9} and Corollary \ref{Le10}, we now may show that $\nabla u(0)$, or rather
a smooth average $\int\omega\nabla u$ of $\nabla u$ near the origin, substantially reacts to a change
in the medium at $B_R(y)$. This reaction is characterized in terms of the tensor 
\begin{align}\label{g96}
\delta a_{hij}:=\int(e_j+\nabla\phi_{jh}')\cdot(a_h'-a_h)e_i
\end{align}
appearing in (\ref{g41}).
In fact, Theorem \ref{T2} is not inferred from Proposition \ref{P2},
but rather directly from Corollary \ref{Le10} and Lemma \ref{L9},
which we combine for that purpose to

\begin{corollary}\label{Le11}
We make the same assumptions on $a$, $a_h$, $\phi_i$, $\sigma_i$, $r_*$, $r_*(y)$ as in Proposition \ref{P2}.
Let $a_0$, $a'$, $a_h'$, $u$, and $u'$ be as in Proposition \ref{P2}. Then provided $R\le\ell\le|y|$ we have
\begin{align}\label{j05bis}
\lefteqn{\big|\int\omega\nabla(u'-u)
-\int g_m\partial_m\partial_iG_h(y)\delta a_{hij}\partial_j\partial_kG_h(-y)e_k\big|\le C
(\frac{\ell}{|y|})^d(\frac{R}{|y|})^{d}}\nonumber\\
&\times\Big((\frac{R}{|y|})^{\beta}+(\frac{r_*}{R})^\beta+(\frac{r_*(y)}{R})^\alpha
+\frac{\ell}{|y|}+(\frac{r_*}{\ell})^\beta\Big)\big(1+\frac{r_*+r_*(y)}{R}\big)^{\frac{3}{2}d+\beta}.
\end{align}
where $\omega$ is defined in (\ref{g18}).
Here $C$ depends on $d$, $\lambda$, $\beta$, and the $L^2$-norm of $\hat\nabla \hat g$.
\end{corollary}

The first rhs summand in (\ref{j05bis}) comes from Lemma \ref{L9}, the second one
from eliminating $\nabla\phi_k$ in (\ref{g85}), the third one comes from
Corollary \ref{Le10}, the fourth comes from the dipole expansion of $u_h$,
and the fifth from simplifying the dipole moment. The last factor on the rhs of (\ref{j05bis}) arises because
we avoid any smallness condition on $\max\{r_*,r_*(y)\}$ in terms of $R$ or $\ell$.

\medskip

Under our assumptions on the ensemble, cf Definition \ref{D1}, and
equipped with the deterministic result of Corollary \ref{Le11}, we obtain
a lower bound on the variance of the expectation of
$\int\omega\nabla u$ conditioned on the Poisson process restricted to
$B_R(y)$.  Note that we first average over the Poisson process in the
complement $B_R^c(y)$, and then consider the variance wrt the Poisson process on $B_R(y)$
 --- as opposed to the opposite order, 
which would amount to a weaker result, in
particular too weak for the purpose of Theorem \ref{T2}.  This lower
bounds holds provided the order-one radius $R$, which also governs the
average through $\omega$, cf (\ref{g18}), is sufficiently large.  The lower bound is
optimal in terms of the scaling in the ratios of the length scales $1$
(the range of dependence of the random coefficient field $a$), $\ell$ (the
scale of the source $g$), and $|y|$ (the distance between the source
and the site at which the variance is probed).

\begin{lemma}\label{Le12}
Under the assumptions of Theorem \ref{T2} there exists an $R$ 
so that provided $C\le \ell\le\frac{1}{C}|y|$ we have 
\begin{align}\label{l05}
\big\langle\big|\langle\int\omega\nabla u|B_R(y)\rangle-\langle\int\omega\nabla u\rangle\big|^2\big\rangle^\frac{1}{2}
&\ge\frac{1}{C}(\frac{\ell}{|y|})^d(\frac{1}{|y|})^d|\int\hat g|,
\end{align}
where $R$ and $C$ are as in Theorem \ref{T2}.
\end{lemma}
\medskip

The following auxiliary lemma will be used in the proof of Lemma \ref{L9}
and states that the random (scalar) field $r_*=r_*(a,x)$ is not very sensitive to local changes in the coefficient field $a$.
Here we think of $r_*$ as being defined as the minimal radius with the property (\ref{L04});
we denote by $r_*'$ the radius for the medium $a'$, and by $r_*(y), r_*'(y)$
the corresponding radii with the origin replaced by the point $y$ (in the neighborhood of which
$a'$ differs from $a$). Here, and in analogy to $\phi_i'$, we think of $\sigma_{ijk}'$ as the solution
of $-\Delta\sigma_{ijk}'=\partial_jq'_{ik}-\partial_kq'_{ij}$, where $q_i':=a'(e_i+\nabla\phi_i')$,
that behaves at $\infty$ as $\sigma_{ijk}$ (ie $\int|\nabla(\sigma_i'-\sigma_i)|^2<\infty$). We note that
by construction we have $\Delta(q_i'-\nabla\cdot\sigma_i')=0$, so that by the decay of $q_i'-q_i$
and $\nabla\sigma_i'-\nabla\sigma_i$, the relation (\ref{L01}) is preserved, that is,
$q_i'=a_he_i+\nabla\cdot\sigma'_i$.

\begin{lemma}\label{L7}
We have 
\begin{align}\label{g01}
r_*'(y)+r_*'\le C(d,\beta)(r_*(y)+r_*+R).
\end{align}
\end{lemma}


\section{Numerical results}\label{Num}

In our numerical tests, we consider a random ensemble according to
Definition~\ref{D1}. Let $X$ be the Poisson point process on
$\mathbb{R}^2$ with unit intensity, the coefficient field $a$ is given by 
\begin{equation}
  a(x) = \max_{ \xi_i \in X } \bigl\{ \varphi (x - \xi_i) \bigr\} + \frac{1}{2},
\end{equation}
where $\varphi$ is a bump function supported in $B_{1/2}$
\begin{equation}
  \varphi(\hat{x}) = \exp\biggl(- \frac{1}{1 - 4 \abs{\hat{x}}^2} \biggr). 
\end{equation}

The elliptic PDE is 
\begin{equation}
  - \nabla \cdot a \nabla u = f
\end{equation}
where the rhs is given by 
\begin{equation}
  f(x) = x_2 \exp\biggl(- \frac{5}{5 - \abs{x}^2} \biggr)
\end{equation}
and thus $f$ is compactly supported inside the ball $B_{\sqrt{5}}$
with average $0$. It is clear that we can rewrite $f = \nabla \cdot g$ with $g$ supported in $B_{\sqrt{5}}$.

To numerically approximate the solution of the Dirichlet problems
\eqref{L95}, \eqref{L96}, and \eqref{L56}, we use a standard
second-order centered finite-difference scheme with mesh size
$\delta x = 0.1$. 

\begin{figure}[ht]
  \centering
  \includegraphics[width = .7\textwidth]{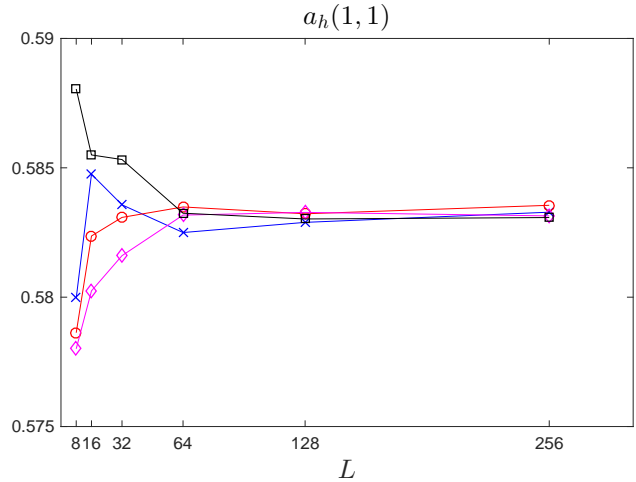}
  \caption{Finite domain approximation to the homogenized
    coefficient. Different lines correspond to different independent
    realization of the coefficient field $a$. } \label{fig:ahom}
\end{figure}

In Figure~\ref{fig:ahom}, we plot the numerically obtained $\aL_h$ as
a function of $L$ for $4$ independent realizations of the coefficient
field $a$. It shows that the Dirichlet approximation $\aL_h$ converges as $L$ increases, validating \eqref{L47} in Lemma~\ref{L3}.

\begin{figure}[ht]
  \centering
  \includegraphics[width = .7\textwidth]{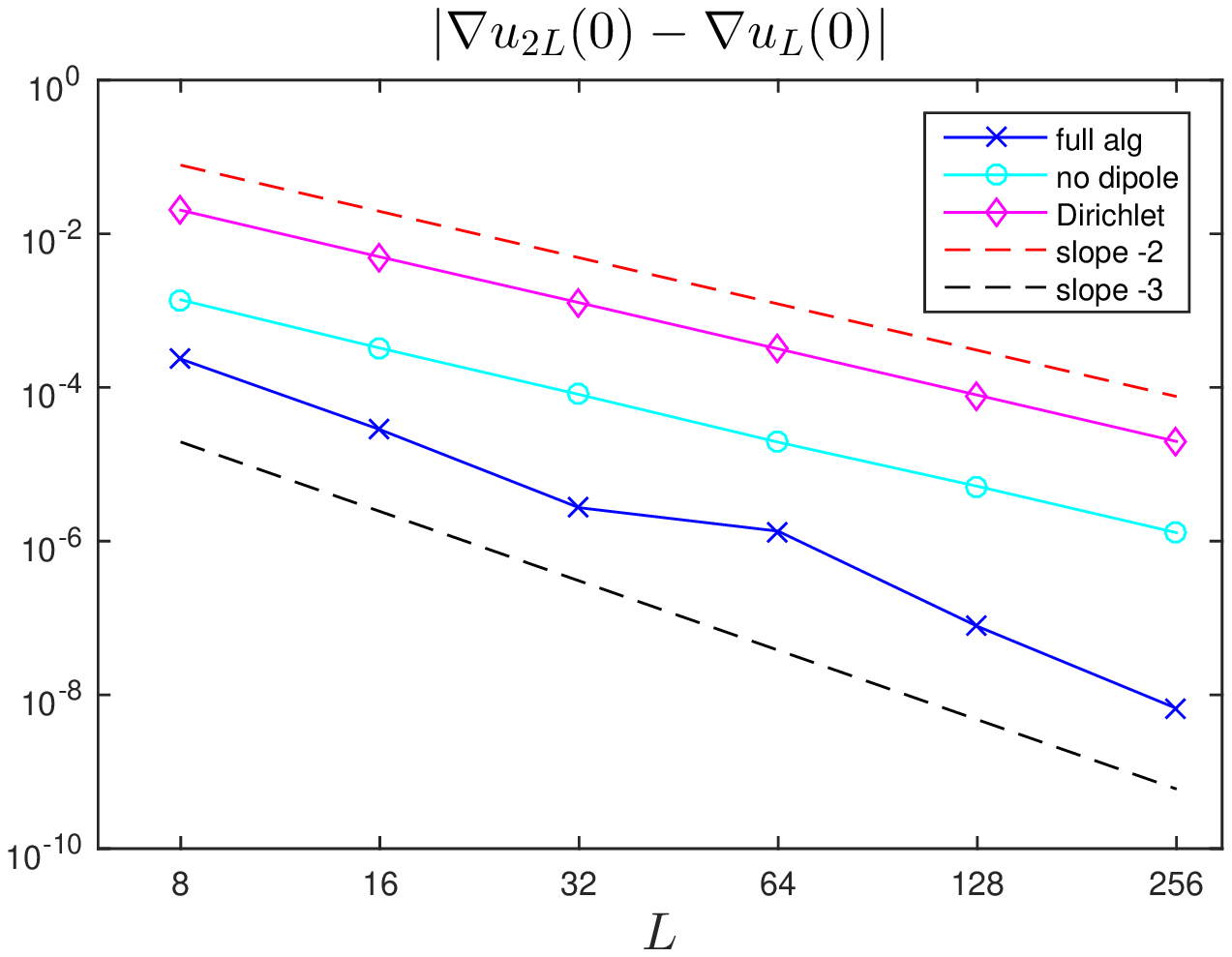}\\
  \includegraphics[width = .7\textwidth]{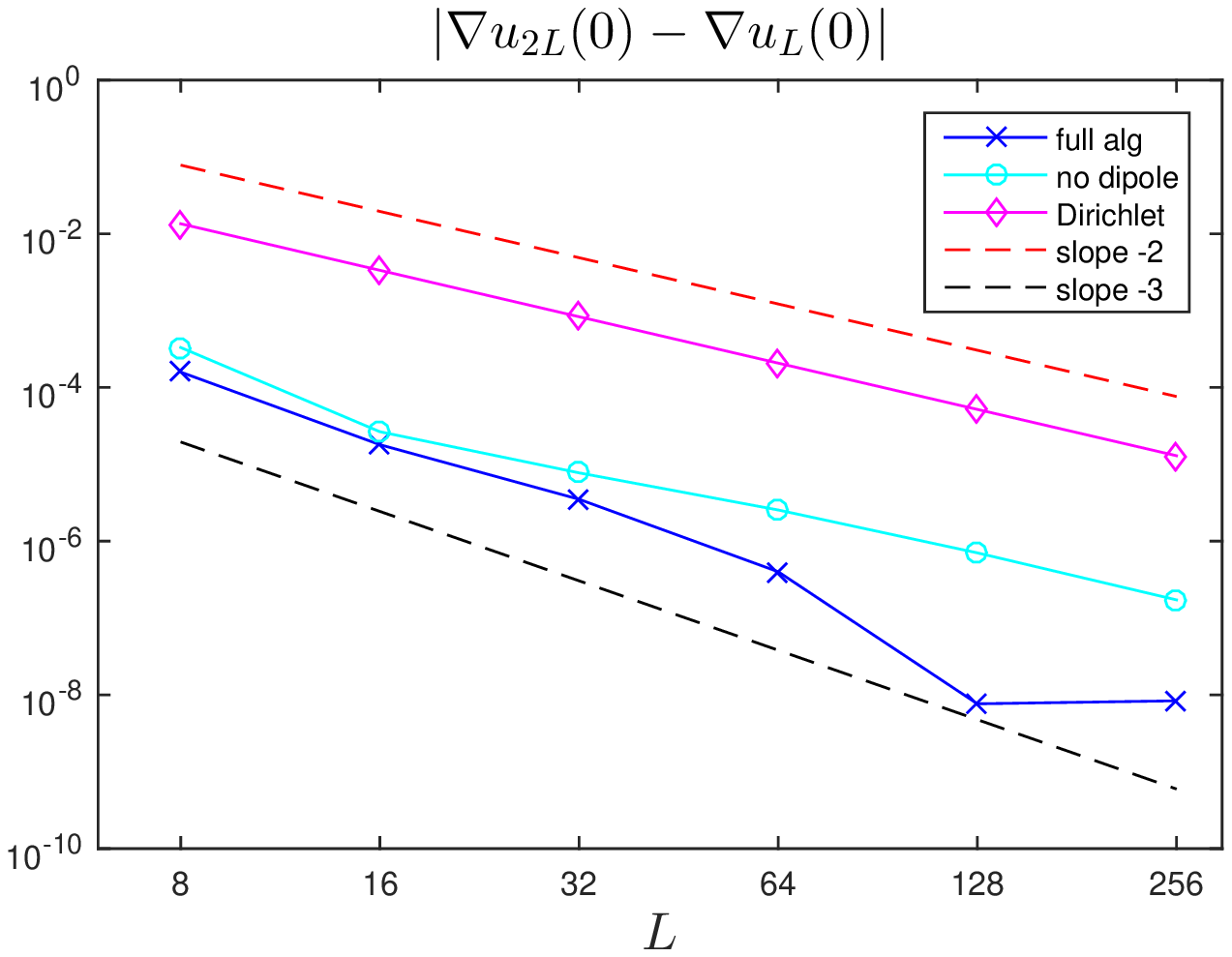}
  \caption{Numerical convergence rate of $\nabla u(0)$ for the
    proposed algorithm, an approximation without the dipole
    correction, and the Dirichlet approximation. The two figures
    correspond to two independent realization of the random
    medium.} \label{fig:convgradu}
\end{figure}

We now consider our algorithm to approximate $\nabla u(0)$ based on
the $a \vert_{Q_{2L}}$. To validate the algorithm, we check the
numerical convergence rate and compare it with two approximations that
with a slower convergence rate. Recall that our algorithms consists of
the steps:
\begin{itemize}
\item Solve the equation \eqref{L95} for approximate correctors $\phiL_i$;
\item calculate the homogenized coefficient $\aL_h$ as in \eqref{L99};
\item obtain $\uL_h$ from the homogenized equation \eqref{L48} and the dipole correction \eqref{L52}; 
\item solve the equation \eqref{L56} for $\uL$. 
\end{itemize}
In comparison, we consider two other algorithms: 1) Solving the
equation with homogeneous Dirichlet boundary condition, ie
\eqref{L88}, referred as ``Dirichlet algorithm''; 2) Dropping the
dipole correction, ie, instead of \eqref{L56}, one solves \eqref{L90}
using approximate homogenized coefficient and correctors
\begin{equation}
  - \nabla \cdot a \nabla \uL_{II}  = \nabla \cdot g \ \text{in}\ Q_L, \quad
  \uL_{II} = ( 1 + \phiL_i \partial_i) \tuL_h \ \text{on}\ \partial Q_L. 
\end{equation}
Thus in the boundary condition, $\uL_h$ is replaced by $\tuL_h$ so the
dipole correction is dropped. This will be referred as the ``no dipole
algorithm''. In Figure~\ref{fig:convgradu}, we compare the numerical
convergence for the three algorithms with two different realization of
the random media. The difference
$\nabla u^{(2L)}(0) - \nabla u^{(L)}(0)$ is plotted for various $L$
for the three algorithms. As the plot is in loglog scale, a straight
line shows algebraic convergence with the slope indicating the
convergence rate. We observe from the numerical results that the
proposed algorithm achieves almost $O(L^{-3})$ convergence, which is
consistent with the theoretical rate
$O((\tfrac{\ell}{L})^2 (\tfrac{1}{L})^{\beta})$ with $\beta$ close to
$1$. The wiggled line is due to the randomness of the coefficients
added in the annulus $Q_{4L} - Q_{2L}$ when $L$ is doubled.  We also
observe that the other two algorithms have the slower rate of
$O(L^{-2})$, the main reason being that they do not capture the
correct dipole in the far field. It can be also seen that by using the
information from homogenization, the ``no dipole algorithm'' indeed
reduces the error compared with the simple minded Dirichlet
approximation.

Finally, we study the error committed by the proposed algorithm with
the uncertainty of $\nabla u(0)$ due to the unknown coefficient
outside $Q_{2L}$. In the second row of Table~\ref{tab:lower}, we
compare the approximation $\nabla \uL(0)$ with the reference solution
obtained by a calculation with a large domain ($L = 768$ is
chosen). In the remaining rows, we show how $\nabla u(0)$ varies 
when we re-sample the random medium outside $Q_{2L}$ (where again, we take $\nabla u^{(768)}(0)$ as a proxy for $\nabla u(0)$). These thus show the sensitivity of $\nabla u(0)$ with 
respect to the change of the media. We observe from the numerical
results that while the error made by our algorithm is larger than the
typical fluctuations  of $\nabla u(0)$, the decay rate matches and the error is
also comparable to the sensitivity for various $L$.

\begin{table}[ht]
  \begin{adjustwidth}{-.5in}{-.5in}  
  \begin{tabular}{c|cccc}
    \hline
    & $L = 32$ & $L =  64$ & $L = 128$ & $L = 256$ \\
    \hline
    $\abs{\nabla \uL (0) - \nabla u^{(768)}(0)}$ & \num{2.7709e-6} & \num{1.2843e-6} & \num{8.3910e-8} & \num{5.8996e-9} \\
    \hline 
    \multirow{3}{0.35\textwidth}{\centering $\nabla u(0)$ difference with re-sampled media in $Q_{2L}^c$} 
    & \num{4.8919e-8} & \num{7.7095e-9} & \num{4.9692e-9} & \num{7.8225e-10} \\
    & \num{5.3976e-8} & \num{1.9546e-8} & \num{8.6449e-9} & \num{6.7321e-10} \\
    & \num{1.8771e-7} & \num{3.1458e-8} & \num{1.0232e-8} & \num{5.9754e-10} \\
    \hline
   \end{tabular}
    \smallskip
   \caption{The error of the proposed numerical scheme compared with sensitivity of $\nabla u(0)$ when the media is re-sampled outside $Q_{2L}$. Reference value $\abs{\nabla u^{(768)}(0)} = \num{3.5859e-1}$.} \label{tab:lower}
   \end{adjustwidth}
  
\end{table}


\section{Proofs}

\begin{proof}[\sc Proof of Lemma \ref{L1}]\

We start with the remark that (\ref{L04}) actually implies (\ref{L32}).
%
%
Indeed, what separates (\ref{L32}) from (\ref{L04}) is 
\begin{align}\label{g58}
|\av_{B_r}(\phi,\sigma)-\av_{B_{r'}}(\phi,\sigma)|
\lesssim r(\frac{r_*}{r})^\beta
\quad\mbox{for all}\;r\ge r'\ge r_*.
\end{align}
Since $\beta<1$, we may apply dyadic decomposition to reduce this to
\begin{align*}
|\av_{B_r}(\phi,\sigma)-\av_{B_{r'}}(\phi,\sigma)|
\lesssim r(\frac{r_*}{r})^\beta
\quad\mbox{for all}\;2r'\ge r\ge r'\ge r_*,
\end{align*}
which by the triangle inequality is a consequence of (\ref{L04}). In view of (\ref{L32}),
we may change $(\phi,\sigma)$ by an additive constant, and thus without affecting (\ref{L01})
to upgrade (\ref{L04}) to
\begin{align}\label{L33}
\frac{1}{r}\big(\av_{B_r}|(\phi,\sigma)|^2\big)^\frac{1}{2}
\lesssim(\frac{r_*}{r})^\beta
\quad\mbox{for all}\;r\ge r_*.
\end{align}
This is in line with (\ref{L62}).

\medskip

For later reference we also note that
\begin{align}\label{L55}
\big(\av_{B_R}|e_i+\nabla\phi_i|^2\big)^\frac{1}{2}\lesssim 1\quad\mbox{for}\;R\ge r_*,
\end{align}
which is a consequence of Caccioppoli's estimate on the $a$-harmonic function $x_i+\phi_i$, cf (\ref{L31}),
and (\ref{L33}) in form of $\big(\av_{B_R}(x_i+\phi_i)^2\big)^\frac{1}{2}$
$\lesssim R(\frac{r_*}{R})^\beta$.

\medskip

We now follow the steps of the proof of \cite[Theorem 0.2]{BellaGiuntiOttoPCMINotes} and like there may wlog assume that $r_*=1$.
In analogy to Step 2 in the proof of \cite[Theorem 0.2]{BellaGiuntiOttoPCMINotes} we claim
\begin{align}\label{io11}
\big(\int|\nabla u|^2\big)^\frac{1}{2}\lesssim\ell^{\frac{d}{2}},
\end{align}
and
\begin{align}\label{L06}
|x|^d|\nabla u_h(x)|+|x|^{d+1}|\nabla^2u_h(x)|\lesssim\ell^d\quad\mbox{for all}\;x\not=0.
\end{align}
We furthermore claim that the $u$ and $u_h$ have not only vanishing ``constant invariant''
\begin{align}\label{io14}
\int\nabla\eta_r\cdot a\nabla u=\int\nabla\eta_r\cdot a_h\nabla u_h=0\quad\mbox{for}\;r\ge\ell,
\end{align}
but also identical ``linear invariants'', that is,
\begin{align}\label{io13}
\lefteqn{\int\nabla\eta_r\cdot((x_i+\phi_i)a\nabla u-ua(e_i+\nabla\phi_i))}\nonumber\\
&=\int\nabla\eta_r\cdot(x_ia_h\nabla u_h-u_ha_he_i)\quad\mbox{for}\;r\ge\ell,
\end{align}
where $\eta_r(x)=\eta(\frac{x}{r})$, and $\eta$ is a cut-off function for $B_1$ in $B_2$.
Indeed, (\ref{io11}) follows from the energy estimate for (\ref{L98}) and the form (\ref{L97})
of the rhs. By the triangle inequality, we split (\ref{L06}) into the corresponding estimate
for $\tilde u_h$ and for $\delta u_h:=\xi_i\partial_iG_h$ with $\xi_i:=\int\nabla\phi_i\cdot g$.
We first turn to $\tilde u_h$ defined through (\ref{L24}) which in view of (\ref{L97}) is of
the form $\nabla \tilde u_h(x)=\hat\nabla\hat u_h(\frac{x}{\ell})$ with $-\hat\nabla\cdot a_h\hat\nabla \hat u_h$
$=\hat\nabla\cdot\hat g$, so that the desired estimate follows from
\begin{align*}
|\hat\nabla \hat u_h|\lesssim\frac{1}{|\hat x|^d+1}\quad\mbox{and}\quad
|\hat\nabla^2 \hat u_h|\lesssim\frac{1}{|\hat x|^{d+1}+1},
\end{align*}
which in turn are a consequence of standard Schauder theory based on the H\"older continuity of $\hat\nabla\cdot\hat g$
for the near-field, and the Green's function representation for the far-field. 
We now turn to $\delta u_h$ and note that the estimates
\begin{align*}
|\nabla \delta u_h|\lesssim\frac{\ell^d}{|x|^d}\quad\mbox{and}\quad
|\nabla^2 \delta u_h|\lesssim\frac{\ell^d}{|x|^{d+1}}
\end{align*}
follow from the homogeneity of $\nabla^2 G_h$ of degree $-d$ and the fact that $|\xi_i|\lesssim\int_{B_\ell}|\nabla\phi_i|$, 
Jensen's inequality together with (\ref{L55}).
We now turn to (\ref{io14}). The identity for $u$ follows directly from (\ref{L98}) and the fact that $g$ is
supported in $B_\ell$, cf (\ref{L97}). For the same reason, the contribution to the constant invariant
of $u_h$ coming from $\tilde u_h$ vanishes; the one coming from $\delta u_h$, which is constant in $r>0$ by $a_h$-harmonicity
(thus the name ``invariant''),
must vanish by its homogeneity of order $1-d$.
We finally turn to (\ref{io13}). By (\ref{L98}) and (\ref{L31}) we learn from an integration by parts and the support condition
on $g$ that for $r\ge\ell$,
\begin{align}\label{io14bis}
\int\nabla\eta_r\cdot((x_i+\phi_i)a\nabla u-ua(e_i+\nabla\phi_i))
=\int(e_i+\nabla\phi_i)\cdot g.
\end{align}
Likewise, we obtain from (\ref{L24}) that
\begin{align*}
\int\nabla\eta_r\cdot(x_ia_h\nabla \tilde u_h-\tilde u_ha_he_i)=\int e_i\cdot g.
\end{align*}
Finally, by definition of the Green's function $G_h$ we have
\begin{align*}
\int\nabla\eta_r\cdot(x_ia_h\nabla \partial_j G_h-\partial_j G_ha_he_i)=\delta_{ij}.
\end{align*}
By $u_h=\tilde u_h+\xi_j\partial_jG_h$ and the definition of $\xi$, the last two identities combine to
\begin{align*}
\int\nabla\eta_r\cdot((x_ia_h\nabla u_h-u_ha_he_i)=\int(e_i+\nabla\phi_i)\cdot g.
\end{align*}
This and (\ref{io14bis}) yield the identity (\ref{io13}).

\medskip

Following Step 3 in the proof of \cite[Theorem 0.2]{BellaGiuntiOttoPCMINotes} (with $\eta=0$, the additional multiplicative factor $\ell^d$,
and $\beta$ playing the role of $1-\alpha$)
we consider the error in the two-scale expansion $w:=u-(1+\phi_i\partial_i)u_h$ we claim that
\begin{align}\label{io12}
\big(\int_{B_1^c}|\nabla w|^2\big)^\frac{1}{2}\lesssim\ell^{d}
\end{align}
and that
\begin{align}\label{L14}
-\nabla\cdot a\nabla w=\nabla\cdot h\quad\mbox{with}\quad h:=(\phi_ia-\sigma_i)\nabla\partial_i u_h,
\end{align}
where 
\begin{align}\label{L06bis}
\big(\frac{1}{R^d}\int_{B_r^c}|h|^2\big)^\frac{1}{2}\lesssim\ell^d(\frac{1}{r})^{d+\beta}\quad\mbox{for}\;r\ge 1.
\end{align}
We first turn to (\ref{io12}) and note
$\nabla w=\nabla u-\partial_iu_h(e_i+\nabla\phi_i)-\phi_i\nabla\partial_iu_h$. Hence for the $\nabla u$-contribution
the desired estimate follows from (\ref{io11}) (and $\ell\ge 1$). According to (\ref{L06}), the contribution
from $\partial_iu_h(e_i+\nabla\phi_i)$ is
estimated by $\ell^d\big(\int_{B_1^c}|x|^{-2d}|e_i+\nabla\phi_i|^2\big)^\frac{1}{2}$; dividing the integral into
dyadic annuli and using (\ref{L55}), we see that it is $\lesssim 1$. Still according to (\ref{L06}),
the contribution from $\phi_i\nabla\partial_iu_h$ is estimated by 
$\ell^d\big(\int_{B_1^c}|x|^{-2(d+1)}\phi_i^2\big)^\frac{1}{2}$; dividing also this integral into
dyadic annuli and using (\ref{L33}), we see that it is $\lesssim 1$.
As $-\nabla\cdot a\nabla u=\nabla\cdot g=\nabla\cdot a_h\nabla u_h$ (in view of (\ref{L98}), (\ref{L24}),
and the $a_h$-harmonicity of $G_h$ and thus $\delta u_h$ away from the origin), 
the identity (\ref{L14}) follows by a straightforward calculation
as in Step 3 of \cite{BellaGiuntiOttoPCMINotes}. It is the merit of $\sigma_i$ that the
rhs can be written in divergence-form $\nabla\cdot h$. 
We finally turn to the estimate (\ref{L06bis}); by (\ref{L06}) and (\ref{L03}) we have
$|h|\lesssim\ell^d|x|^{-(d+1)}|(\phi,\sigma)|$. Hence the estimate follows once more from (\ref{L33}) 
and division into dyadic annuli.

\medskip

Following Step 4 and Step 5 of \cite[Theorem 0.2]{BellaGiuntiOttoPCMINotes}
(still with $\alpha$ replaced by $1-\beta$ and the additional factor of $\ell^d$)
we obtain
\begin{align*}
\big(\frac{1}{R^d}\int_{B_r^c}|\nabla w|^2\big)^\frac{1}{2}\lesssim\ell^d(\frac{1}{r})^{d+\beta}\quad\mbox{for}\;r\ge 2,
\end{align*}
which (in view of $r_*=1$) turns into (\ref{L17}). Note that the outcome of Step 5 of \cite[Theorem 0.2]{BellaGiuntiOttoPCMINotes}
is worse by a $\log R$; which however can be easily avoided, cf \cite[Lemma 3]{BellaGiuntiOttoarXiv}
with $k=1$ and $r=r_*=1$.


\end{proof}

\medskip


\begin{proof}[\sc Proof of Corollary \ref{L2}]\ 

Let us introduce a tool we will often use, namely the $C^{0,1}$-estimate for an $a$-harmonic function $w$
in $B_R$, the crucial role of which was recognized in \cite{AS}. As is obvious from its characterization (\ref{L04}), $r_*$ 
dominates, up to a multiplicative constant $C=C(d,\lambda,\beta)$, the ``minimal radius'' introduced in \cite{GNOarXiv}
and from which on the ``mean-value property'', as an estimate, holds for $|\nabla w|^2$, see \cite[Theorem 1]{GNOarXiv} for
the proof. We record this
\begin{align}\label{io10}
\av_{B_r}|\nabla w|^2\le C(d,\lambda,\beta)\av_{B_R}|\nabla w|^2\quad\mbox{for}\;r_*\le r\le R.
\end{align}

\medskip

We note that the combination of (\ref{L24}), (\ref{L18}) and (\ref{L94}) may be rephrased in terms of $w:=u_{III}-u$ 
and $w_D:=u-(1+\phi_i\partial_i)u_h$ as
\begin{align}\label{L28}
-\nabla\cdot a\nabla w=0\;\mbox{in}\;Q_L\quad\mbox{and}\quad w=w_D\;\mbox{on}\;\partial Q_L.
\end{align}
Once we show that this implies
\begin{align}\label{L27}
\big(\int_{Q_L}|\nabla w|^2\big)^\frac{1}{2}\lesssim
\big(\int_{Q_{2L}-Q_L}|\nabla w_D|^2\big)^\frac{1}{2},
\end{align}
we see that (\ref{L26}) follows from (\ref{L17}) for $R=L$, and a subsequent application of the $C^{0,1}$-estimate
for $w$ to get from $Q_L$ to $B_R$ for $R\ge r_*$. 

\medskip

We now turn to the argument that (\ref{L28}) implies (\ref{L27}) under
the mere assumption of uniform ellipticity, cf (\ref{L03}). Hence by rescaling, we may wlog
assume $L=1$. Since (\ref{L28}) and (\ref{L27}) are oblivious to additive
constant we may wlog assume $\int_{Q_{2}-Q_1}w_D=0$, which allows to construct an extension
$\bar w_D$ of $w_D$ on $Q_1$ such that
\begin{align}\label{L29}
\big(\int_{Q_2}|\nabla\bar w_D|^2\big)^\frac{1}{2}\lesssim
\big(\int_{Q_2-Q_1}|\nabla w_D|^2\big)^\frac{1}{2}.
\end{align}
This extension allows us to reformulate (\ref{L28}) as
\begin{align*}
-\nabla\cdot a\nabla(w-\bar w_D)=\nabla\cdot a\nabla\bar w_D\;\mbox{in}\;Q_1,
\quad w-\bar w_D=0\;\mbox{on}\;\partial Q_1.
\end{align*}
Hence (\ref{L27}) follows by the energy estimate on the latter, followed by a triangle
inequality in $L^2$ yielding $\big(\int_{Q_1}|\nabla w|^2\big)^\frac{1}{2}$
$\lesssim \big(\int_{Q_1}|\nabla \bar w_D|^2\big)^\frac{1}{2}$, into which we insert
(\ref{L29}).
\end{proof}

\medskip

\begin{proof}[\sc Proof of Lemma \ref{L3}]\

We first turn to (\ref{L44}) and, for notational simplicity, drop the index $i=1,\cdots,d$.
By (\ref{L31}), which is as mentioned a consequence of (\ref{L01}), and (\ref{L95}) 
we have $-\nabla\cdot a\nabla(\phiL-\phi)$ $=0$ in $Q_{2L}$,
so that by Caccioppoli's estimate and the triangle inequality in $L^2$ we have
\begin{align*}
\big(\av_{Q_{\frac{3}{2}L}}&|\nabla(\phiL-\phi)|^2\big)^\frac{1}{2}
\lesssim\frac{1}{L}\big(\av_{Q_{2L}}((\phiL-\phi)-\av_{Q_{2L}}\phi)^2\big)^\frac{1}{2}\nonumber\\
&\le\frac{1}{L}\big(\av_{Q_{2L}}(\phiL)^2\big)^\frac{1}{2}
+\frac{1}{L}\big(\av_{Q_{2L}}(\phi-\av_{Q_{2L}}\phi)^2\big)^\frac{1}{2}.
\end{align*}
Hence (\ref{L44}) follows from the $\phi$-parts in (\ref{L70}) and (\ref{F19}) .

\medskip

For (\ref{F20}), we argue in a similar way: According to (\ref{L01}), which yields (\ref{L30}), and (\ref{L96})
we have $-\Delta(\sigmaL_{jk}-\sigma_{jk})$ $=\partial_j(\qL_k-q_k)-\partial_k(\qL_j-q_j)$ in $Q_{2L}$, 
where we recall $q=a(e+\nabla\phi)$, $\qL=a(e+\nabla\phiL)$ and thus $\qL-q =a\nabla(\phiL-\phi)$.
By Caccioppoli's estimate and the upper bound on $a$ provided by (\ref{L03}) we therefore obtain
\begin{align*}
\lefteqn{\big(\av_{Q_{L}}|\nabla(\sigmaL-\sigma)|^2\big)^\frac{1}{2}}\nonumber\\
&\lesssim\big(\av_{Q_{\frac{3}{2}L}}|\nabla(\phiL-\phi)|^2\big)^\frac{1}{2}
+\frac{1}{L}\big(\av_{Q_{2L}}|\sigmaL|^2\big)^\frac{1}{2}
+\frac{1}{L}\big(\av_{Q_{2L}}|\sigma-\av_{Q_{2L}}\sigma|^2\big)^\frac{1}{2}.
\end{align*}
Hence (\ref{F20}) follows from (\ref{L44}) and the $\sigma$-parts of both (\ref{L70}) and (\ref{F19}).

\medskip

We now turn to (\ref{L47}) and recall our averaging function $\omega$, cf (\ref{F09}).
We first claim that also for the whole-space case,
the effective coefficient can approximately be recovered from averaging the flux wrt to $\omega$:
\begin{align}
|a_he-\int\omega a(e+\nabla\phi)|\lesssim(\frac{r_*}{L})^\beta.\label{L45}
\end{align}
%
The argument for this is based on (\ref{L01}), which yields the identity
$a_he-\int\omega a(e+\nabla\phi)$ $=\int(\sigma-\av_{B_L}\sigma)$ $\nabla\omega$, and thus by (\ref{F09}) the estimate 
$|a_he$ $-\int\omega a(e+\nabla\phi)|$ $\lesssim\frac{1}{L}\av_{B_L}|\sigma-\av_{B_L}\sigma|$. Hence we obtain
(\ref{L45}) from (\ref{F19}) and Jensen's inequality. 

\medskip

We now may conclude: By the upper bound provided in (\ref{L03}), we obtain from
(\ref{L45}) and the definition (\ref{L99}) by the triangle inequality
\begin{align*}
|(\aL_h-a_h)e|\lesssim\av_{Q_L}|\nabla(\phiL-\phi)|+(\frac{r_*}{L})^\beta,
\end{align*}
so that (\ref{L47}) follows from (\ref{L44}) and Jensen's inequality.
\end{proof}

\medskip

\begin{proof}[\sc Proof of Proposition \ref{P1}]\

We first compare the two solutions $\tuL_h$ and $\tilde u_h$ of (\ref{L48}) and (\ref{L24}), respectively,
and claim that
\begin{align}\label{L51}
|\nabla(\tuL_h-\tilde u_h)|+L|\nabla^2(\tuL_h-\tilde u_h)|
\lesssim (\frac{r_*}{L})^{\beta}(\frac{\ell}{L})^d\quad\mbox{on}\;Q_L^c.
\end{align}
Indeed, the difference $w:=\tuL_h-\tilde u_h$ satisfies
\begin{align*}
-\nabla\cdot \aL_h\nabla w=\nabla\cdot (\aL_h-a_h)\nabla\tilde u_h.
\end{align*}
The form (\ref{L97}) of the rhs $g$ in (\ref{L24}) transmits to the solution $\tilde u_h$:
$\nabla\tilde u_h(x)$ $=\hat\nabla\widehat{\tilde u_h}(\frac{x}{\ell})$, and then also to $w$:
$\nabla w(x)$ $=\hat\nabla\hat w(\frac{x}{\ell})$. Recall that $\hat g$ is compactly supported
in $B_1$ with H\"older continuous derivatives; in particular, $\hat g$ is in the class of vector fields 
decaying as $O(|x|^{-d})$ with derivatives decaying as $O(|x|^{-d-1})$ and differences of the derivatives
decaying as $O(|x|^{-d-1-\alpha}|x-y|^\alpha)$, in line with $\alpha$-H\"older continuity for some $\alpha\in(0,1)$.
This class is preserved under the constant-coefficient Helmholtz projection (recall that $\aL_h$ satisfies
(\ref{L03})). Hence we obtain in particular
\begin{align*}
|\hat\nabla\hat w|\lesssim|\aL_h-a_h|\frac{1}{(|\hat x|+1)^d},&\quad
|\hat\nabla^2\hat w|\lesssim|\aL_h-a_h|\frac{1}{(|\hat x|+1)^{d+1}}.
\end{align*}
Translating back to the microscopic variables and using (\ref{L47}) in Lemma \ref{L3} we have
\begin{align*}
|\nabla w|\lesssim (\frac{r_*}{L})^{\beta}(\frac{\ell}{|x|+\ell})^d,&\quad
\ell|\nabla^2w|\lesssim(\frac{r_*}{L})^{\beta}(\frac{\ell}{|x|+\ell})^{d+1},
\end{align*}
from which we extract (\ref{L51}).

\medskip

We now compare $\uL_h$ and $u_h$ defined in (\ref{L52}) and (\ref{L18}) and claim that
\begin{align}\label{L53}
|\nabla(\uL_h-u_h)|+L|\nabla^2(\uL_h-u_h)|
\lesssim (\frac{r_*}{L})^{\beta}(\frac{\ell}{L})^d\quad\mbox{on}\;Q_L^c.
\end{align}
For later purpose, we also record
\begin{align}\label{L61}
|\nabla \uL_h|+L|\nabla^2 \uL_h|
\lesssim(\frac{\ell}{L})^d\quad\mbox{on}\;Q_L^c.
\end{align}
In order to pass from (\ref{L51}) to (\ref{L53}), and from (\ref{L06}) and (\ref{L53}) to (\ref{L61}), it remains to control
the dipole contributions:
\begin{align*}
|\nabla(\xiL_i\partial_i\GL_h-\xi_i\partial_iG_h)|&+L|\nabla^2(\xiL_i\partial_i\GL_h-\xi_i\partial_iG_h)|
\lesssim (\frac{r_*}{L})^{\beta}(\frac{\ell}{L})^d\nonumber\\
\mbox{and}\quad
|\nabla\xi_i\partial_i\GL_h|&+L|\nabla^2\xi_i\partial_i\GL_h|
\lesssim(\frac{\ell}{L})^d\quad\mbox{both on}\;Q_L^c,
\end{align*}
where we have set for abbreviation $\xiL_i$ $:=\int\nabla\phiL_i\cdot g$ and (as above) $\xi_i =\int\nabla\phi_i\cdot g$.
Because of the obvious estimates on the constant-coefficient (and thus homogeneous) fundamental solution
\begin{align*}
|\nabla\partial_i\GL_h|+L|\nabla^2\partial_i\GL_h|\lesssim\frac{1}{L^d}\quad\mbox{on}\;Q_L^c
\end{align*}
and
\begin{align*}
|\nabla\partial_i(\GL_h-G_h)|+L|\nabla^2\partial_i(\GL_h-G_h)|&\lesssim|\aL_h-a_h|\frac{1}{L^d}\nonumber\\
&\lesssim (\frac{r_*}{L})^{\beta}\frac{1}{L^d}\quad\mbox{on}\;Q_L^c,
\end{align*}
it suffices to show
\begin{align*}
|\xi|\lesssim\ell^d\quad\mbox{and}\quad|\xiL-\xi|\lesssim (\frac{r_*}{L})^{\beta}\ell^d.
\end{align*}
By definition of $\xiL,\xi$, the form (\ref{L97}) of $g$, and Jensen's inequality the latter two estimates follow from
\begin{align}\label{L64}
\big(\av_{B_\ell}|\nabla\phi|^2\big)^\frac{1}{2}\lesssim 1
\quad\mbox{and}\quad
\big(\av_{B_\ell}|\nabla(\phiL-\phi)|^2\big)^\frac{1}{2}\lesssim (\frac{r_*}{L})^{\beta}.
\end{align}
The first estimate in (\ref{L64}) follows from (\ref{L55}) since $\ell\ge r_*$. The second estimate 
in (\ref{L64}) follows from (\ref{L44}) in Lemma \ref{L3} and the $C^{0,1}$-estimate for an $a$-harmonic function
$w$ in $Q_L$, cf (\ref{io10}),
which requires $L\ge \ell\ge r_*$.

\medskip

We finally compare $\uL$ defined through (\ref{L56}) with $u_{III}$ defined through (\ref{L94}) and claim that
\begin{align*}
\big(\av_{B_R}|\nabla(\uL-u_{III})|^2\big)^\frac{1}{2}\lesssim (\frac{\ell}{L})^d(\frac{r_*}{L})^{\beta}
\quad\mbox{for}\;L\ge R\ge r_*.
\end{align*}
This allows us to pass from (\ref{L26}) in Corollary \ref{L2} to this proposition's statement (\ref{L57}).
We note that because $w:=\uL-u_{III}$ satisfies
\begin{align*}
-\nabla\cdot a\nabla w=0\quad\mbox{in}\;Q_L,\quad w=w_D\quad\mbox{on}\;\partial Q_L,
\end{align*}
where $w_D:=(1+\phiL_i\partial_i)\uL_h-(1+\phi_i\partial_i)u_h$. Hence we have by (a slight adaptation of)
the argument at the beginning of the proof of Corollary \ref{L2} that
\begin{align*}
\big(\av_{Q_L}|\nabla w|^2\big)^\frac{1}{2}\lesssim\big(\av_{Q_{\frac{3}{2}L}-Q_L}|\nabla w_D|^2\big)^\frac{1}{2}.
\end{align*}
We combine this with the $C^{0,1}$-estimate, cf (\ref{io10}), for $a$-harmonic functions,
%
%
so that it remains to show
\begin{align}\label{L58}
\big(\av_{Q_{\frac{3}{2}L}-Q_L}|\nabla w_D|^2\big)^\frac{1}{2}
\lesssim (\frac{\ell}{L})^d(\frac{r_*}{L})^{\beta}.
\end{align}
We appeal to the triangle inequality in $L^2$ to split $\nabla w_D$ into the four contributions
\begin{align*}
\nabla w_D&=\partial_i(\uL_h-u_h)(e_i+\nabla\phi_i)+\phi_i\nabla\partial_i(\uL_h-u_h)\nonumber\\
&+\partial_i\uL_h\nabla(\phiL_i-\phi_i)+(\phiL_i-\phi_i)\nabla\partial_i\uL_h.
\end{align*}
For the first contribution, (\ref{L58}) follows from the first part of (\ref{L53}) and (\ref{L55}). 
For the second contribution, (\ref{L58}) follows from the second part of (\ref{L53}) and (\ref{L33}), recall
the normalization of $\phi_i$ through $\av_{B_{r_*}}\phi_i=0$, cf (\ref{L62}).
For the third contribution, we appeal to the first part of (\ref{L61}) and (\ref{L44}).
For the fourth contribution we appeal once more to (\ref{L33}), and now also to (\ref{L70}),
which we combine to $\big(\av_{Q_{\frac{3}{2}L}}(\phiL_i-\phi_i)^2\big)^\frac{1}{2}$ $\lesssim L(\frac{r_*}{L})^\beta$. 
In connection with the second part of (\ref{L61}), we see that also that this last contribution is controlled as stated in (\ref{L58}).
\end{proof}

\medskip


\begin{proof}[\sc Proof of Lemma \ref{L5}]\

As mentioned in the introduction, (\ref{L95}) and (\ref{L96}) do not ensure the analogue of (\ref{L01}),
that is, $\qL_i$ $:=a(e_i+\nabla\phiL_i)$ $=\aL_he_i+\nabla\cdot\sigmaL_i$.
However, we will need the $C^{1,\alpha}$-theory for $a$ from \cite[Theorem 1]{GNOarXiv} on scales $\frac{L}{2}\ge r\ge \rL_*$, 
which relies on this identity next to the controlled sublinearity (\ref{F02}). Hence
for arbitrary but fixed $i=1,\cdots,d$,
we need to modify $\sigmaL_i$ to (a still skew symmetric) $\tilde\sigma_i$ with
\begin{align}\label{F10}
\qL_i=\aL_h e_i+\nabla\cdot\tilde\sigma_i\quad\mbox{on}\;Q_L,
\end{align}
while retaining the $\sigma$-part of (\ref{F02}) in (the seemingly weaker) form of
\begin{align}\label{F11}
\frac{1}{r}\big(\av_{Q_r}|\tilde\sigma_i-\av_{Q_r}\tilde\sigma_i|^2\big)^\frac{1}{2}\lesssim (\frac{\rL_*}{r})^\beta
\quad\mbox{for}\;L\ge r\ge \rL_*.
\end{align}
To this purpose, we first argue that the ``defect'' $\gL_i:=\qL_i-\aL_h e_i-\nabla\cdot\sigmaL_i$, 
which as mentioned in the introduction is component-wise harmonic on $Q_{2L}$ by (\ref{L95}) and (\ref{L96}), satisfies
\begin{align}\label{F16}
\big(\av_{Q_L}|\gL_i|^2\big)^\frac{1}{2}\lesssim(\frac{\rL_*}{L})^\beta.
\end{align}
Indeed, this follows immediately from the identity (\ref{L01}) in conjunction with the three estimates
(\ref{L44}), (\ref{F20}), and (\ref{L47}) in Lemma \ref{L3} (applied with $r_*$ replaced by $\rL_*$).

\medskip

We now turn to (\ref{F10}) and (\ref{F11}). In view of the definition of the defect $\gL_i$, Poincar\'e's inequality, 
and (\ref{F02}), it is enough to construct $\delta\tilde\sigma_i$ such that
\begin{align}\label{F19bis}
\gL_i=\nabla\cdot\delta\tilde\sigma_i\quad\mbox{on}\;Q_L
\end{align}
and
\begin{align}\label{F12}
\big(\av_{Q_r}|\nabla\delta\tilde\sigma_i|^2\big)^\frac{1}{2}\lesssim (\frac{\rL_*}{L})^\beta
\quad\mbox{for}\;L\ge r\ge \rL_*.
\end{align}
To this purpose, we extend the restriction of $\gL_i$ to $Q_{L}$ periodically (without changing the notation)
and let $\delta\sigmaL_i$ be the $Q_L$-periodic solution of
\begin{align}\label{F14}
-\Delta\delta\sigmaL_{ijk}=\partial_j\gL_{ik}-\partial_k\gL_{ij}\quad\mbox{with}\quad\av_{Q_L}\delta\sigmaL_{ijk}=0.
\end{align}
By a remark in the introduction, periodic boundary conditions ensure
\begin{align}\label{F18}
\gL_i-\av_{Q_L}\gL_i=\nabla\cdot\delta\sigmaL_i.
\end{align}
We claim that (\ref{F12}) holds for $\delta\sigmaL_i$:
\begin{align}\label{F15}
\big(\av_{Q_r}|\nabla\delta\sigmaL_i|^2\big)^\frac{1}{2}\lesssim (\frac{\rL_*}{L})^\beta
\quad\mbox{for}\;L\ge r\ge \rL_*.
\end{align}
Indeed, since by the energy estimate for the periodic problem (\ref{F14}), that is,
\begin{align*}
\big(\av_{Q_L}|\nabla\delta\sigmaL_i|^2\big)^\frac{1}{2}\lesssim \big(\av_{Q_L}|\gL_i|^2\big)^\frac{1}{2},
\end{align*}
for (\ref{F15}) it suffices to establish the stronger estimate
\begin{align}\label{F17}
\sup_{Q_{\frac{L}{4}}}|\nabla\delta\sigmaL_i|\lesssim \big(\av_{Q_L}|\gL_i|^2\big)^\frac{1}{2}
+\big(\av_{Q_L}|\nabla\delta\sigmaL_i|^2\big)^\frac{1}{2}
\end{align}
and to appeal to (\ref{F16}). Estimate (\ref{F17}) is a consequence of Sobolev's estimate
\begin{align*}
\sup_{Q_{\frac{L}{4}}}|\nabla\delta\sigmaL_i|\lesssim
\big(\av_{Q_{\frac{L}{4}}}|\nabla^{n+1}\delta\sigmaL_i|^2\big)^\frac{1}{2}
+\big(\av_{Q_{\frac{L}{4}}}|\nabla\delta\sigmaL_i|^2\big)^\frac{1}{2},
\end{align*}
where $n$ is an integer larger than $\frac{d}{2}$, of a localized energy estimate for 
the constant-coefficient equation (\ref{F14}) in form of
\begin{align*}
\big(\av_{Q_{\frac{L}{4}}}|\nabla^{n+1}\delta\sigmaL_i|^2\big)^\frac{1}{2}
+\big(\av_{Q_{\frac{L}{4}}}|\nabla\delta\sigmaL_i|^2\big)^\frac{1}{2}
\lesssim
\big(\av_{Q_{\frac{L}{2}}}|\nabla^{n}\gL_i|^2\big)^\frac{1}{2}
+\big(\av_{Q_{\frac{L}{2}}}|\nabla\delta\sigmaL_i|^2\big)^\frac{1}{2},
\end{align*}
and of a localized energy estimate for $-\Delta \gL_i=0$
\begin{align*}
\big(\av_{Q_{\frac{L}{2}}}|\nabla^{n}\gL_i|^2\big)^\frac{1}{2}
\lesssim \big(\av_{Q_L}|\gL_i|^2\big)^\frac{1}{2}.
\end{align*}
In order to pass from (\ref{F18}) to (\ref{F19bis}), 
it is enough to (explicitly) construct an {\it affine} $\delta\tilde\sigma_{ijk}$
with $\av_{Q_{L}}\gL_i$ $=\nabla\cdot\delta\tilde\sigma_{ijk}$ and such that $|\nabla\delta\tilde\sigma_{ijk}|$
$\lesssim|\av_{Q_{L}}\gL_i|$.
 
\medskip

We now note that we control the {\it gradient} of the proxy $(\phiL,\sigmaL)$ down
to scales $\rL_*$:
\begin{align}\label{L75}
\big(\av_{B_r}|\nabla(\phiL,\sigmaL)|^2\big)^\frac{1}{2}\lesssim 1\quad\mbox{for}\;\frac{L}{4}\ge r\ge \rL_*.
\end{align}
Indeed, since by (\ref{L95}), $x_i+\phiL_i$ is $a$-harmonic in $B_L$, we have by Caccioppoli's estimate
\begin{align*}
\big(\av_{B_{\frac{L}{2}}}|e_i+\nabla\phiL_i|^2\big)^\frac{1}{2}
&\lesssim\frac{1}{L}\big(\av_{B_{L}}(x_i+\phiL_i-\av_{B_L}(x_i+\phiL_i))^2\big)^\frac{1}{2}\nonumber\\
&\lesssim 1+\frac{1}{L}\big(\av_{B_{L}}(\phiL_i-\av_{B_L}\phiL_i)^2\big)^\frac{1}{2},
\end{align*}
so that by (\ref{L73}) for $r=L$, we get (\ref{L75}) for $r=\frac{L}{2}$. The remaining range
of $\frac{L}{2}\ge r\ge \rL_*$ follows since there, 
according to our hypotheses (\ref{L73}) for $\phi^{(L)}$ and to (\ref{F10}) and (\ref{F11}) for $\tilde\sigma$, 
the medium $a$ is well-behaved and thus admits the $C^{0,1}$-estimate, cf (\ref{io10}), 
which applied to $x_i+\phiL_i$ yields
\begin{align*}
\big(\av_{B_{r}}|e_i+\nabla\phiL_i|^2\big)^\frac{1}{2}
\lesssim \big(\av_{B_{\frac{L}{2}}}|e_i+\nabla\phiL_i|^2\big)^\frac{1}{2}
\quad\mbox{for}\;\frac{L}{2}\ge r\ge \rL_*.
\end{align*}
This establishes the $\phiL$-contribution to (\ref{L75}), which in particular implies for the flux
$\qL_i=a(e_i+\nabla\phiL_i)$
\begin{align}\label{L76}
\big(\av_{B_r}|\qL_i|^2\big)^\frac{1}{2}\lesssim 1\quad\mbox{for}\;\frac{L}{2}\ge r\ge \rL_*.
\end{align}
By the equation for $\sigmaL_{i}$ with rhs given by the curl of $\qL_i$, cf (\ref{L96}), we obtain from
Caccioppoli's estimate
\begin{align*}
\big(\av_{B_{r}}|\nabla\sigmaL_{ijk}|^2\big)^\frac{1}{2}
\lesssim\frac{1}{r}\big(\av_{B_{2r}}(\sigmaL_{ijk})^2\big)^\frac{1}{2}
+\big(\av_{B_{2r}}|\qL_i|^2\big)^\frac{1}{2},
\end{align*}
so that we obtain the $\sigmaL$-part of (\ref{L75}) from the $\sigmaL$-part of (\ref{L73}) and
from (\ref{L76}).

\medskip

\ignore{
Next we argue that, again on the level of the gradients, the difference between proxy and
true corrector are small on scale $L$:
\begin{align}\label{L79}
\big(\av_{B_\frac{L}{2}}|\nabla(\phiL-\phi)|^2\big)^\frac{1}{2}
+\big(\av_{B_\frac{L}{4}}|\nabla(\sigmaL-\sigma)|^2\big)^\frac{1}{2}\lesssim(\frac{\rL_*}{L})^\beta.
\end{align}
Indeed, by (\ref{L31}) and (\ref{L95}), $\phiL-\phi$ is (component-wise) $a$-harmonic in $B_L$.
Hence by the argument at the beginning of the proof of Lemma \ref{L3} we obtain the $\phi$-part
of (\ref{L79}) from (\ref{L80}) for $R=L$ and (\ref{L73}) for $r=L$. For the $\sigma$-part we note
that by (\ref{L30}) and (\ref{L96}), $\sigmaL_{i}-\sigma_{i}$ solves
a (tensorial) Poisson equation with a rhs given by the curl of $a\nabla(\phiL_i-\phi_i)$. Hence by Caccioppoli
and the triangle inequality in $L^2$
\begin{align*}
\lefteqn{\big(\av_{B_\frac{L}{4}}|\nabla(\sigmaL-\sigma)|^2\big)^\frac{1}{2}}\nonumber\\
&\lesssim\big(\av_{B_L}|\sigmaL|^2\big)^\frac{1}{2}
+\big(\av_{B_L}|\sigma-\av_{B_L}\sigma|^2\big)^\frac{1}{2}
+\big(\av_{B_\frac{L}{2}}|\nabla(\phiL-\phi)|^2\big)^\frac{1}{2},
\end{align*}
so that the $\sigma$-part of (\ref{L79}) follows from the $\phi$-part of it,
together with the $\sigma$-parts of both (\ref{L80}) and (\ref{L73}) for
$R=r=L$.

\medskip
}

We now come to the central piece, namely that the differences between proxy and true corrector are
small {\it down to scales} $\rL_*$:
\begin{align}\label{L81}
\big(\av_{B_r}|\nabla(\phiL-\phi,\sigmaL-\sigma)|^2\big)^\frac{1}{2}
\lesssim(\frac{\rL_*}{L})^\beta\quad\mbox{for}\;\frac{L}{4}\ge r\ge \rL_*.
\end{align}
As noticed above, $\phiL_i-\phi_i$ is $a$-harmonic in $B_L$ so that passing from $r=\frac{L}{2}$,
cf (\ref{L44}) and (\ref{F20}) in Lemma \ref{L3} with $r_*$ replaced by $\rL_*$, 
to $\frac{L}{2}\ge r\ge \rL_*$ follows from the $C^{0,1}$-estimate already used
earlier, (\ref{io10}). This settles the $\phi$-part of (\ref{L81}); in order to deal with the
$\sigma$-part, we need the full $C^{1,\alpha}$-theory for $a$-harmonic functions from \cite[Theorem 1]{GNOarXiv}
(for an $\alpha\in(0,1)$ fixed, say $\alpha=\frac{1}{2}$), which holds for radii $L\ge r\ge \rL_*$ since, as already
remarked above, the medium is well behaved there in the sense that there exist a tensor $a_h^{(L)}$, scalar fields  $\phi_i^{(L)}$,
and skew symmetric vector fields $\tilde\sigma_i$, related by (\ref{F10}) and satisfying the estimates
(\ref{L73}) \& (\ref{F11}).
Applied to the $a$-harmonic function $\phi_i-\phiL_i$ in $B_L$ this yields a vector $\xi$ (which should carry an index $i$) 
such that 
\begin{align*}
\big(\av_{B_r}|\nabla(\phi_i-\phiL_i-\xi_j(x_j+\phiL_j))|^2\big)^\frac{1}{2}
&\lesssim(\frac{r}{L})^\alpha\big(\av_{B_\frac{L}{2}}|\nabla(\phi_i-\phiL_i)|^2\big)^\frac{1}{2}\nonumber\\
\mbox{for all}\;\frac{L}{2}\ge r\ge \rL_*\quad\mbox{and}\;|\xi|&\lesssim \big(\av_{B_\frac{L}{2}}|\nabla(\phi_i-\phiL_i)|^2\big)^\frac{1}{2}.\nonumber
\end{align*}
Recalling the definition of the fluxes $q_i=a(e_i+\nabla\phi_i)$ and $\qL_i=a(e_i+\nabla\phiL_i)$
and inserting (\ref{L80}), we obtain
\begin{align}\label{L83}
\big(\av_{B_r}|q_i-\qL_i-\xi_j\qL_j|^2\big)^\frac{1}{2}
\lesssim(\frac{r}{L})^\alpha(\frac{\rL_*}{L})^\beta\quad\mbox{for}\;\frac{L}{2}\ge r\ge \rL_*,
\quad|\xi|\lesssim(\frac{\rL_*}{L})^\beta.
\end{align}
By (\ref{L30}) and (\ref{L96}) we have that $w:=\sigma_{i}-\sigmaL_{i}-\xi_j\sigmaL_{j}$ solves
a Poisson equation with the curl of $q_i-\qL_i-\xi_j\qL_j$ as rhs. Hence the first part of (\ref{L83})
translates into
\begin{align}\label{L86}
-\Delta w=\nabla\cdot h\quad\mbox{with}\quad\big(\av_{B_r}|h|^2\big)^\frac{1}{2}
\lesssim(\frac{r}{L})^\alpha(\frac{\rL_*}{L})^\beta\quad\mbox{for}\;\frac{L}{2}\ge r\ge \rL_*.
\end{align}
In the next paragraph, we shall argue that thanks to $\alpha>0$, this implies
\begin{align}\label{L87}
(\av_{B_r}|\nabla w|^2\big)^\frac{1}{2}\lesssim(\frac{\rL_*}{L})^\beta+(\av_{B_\frac{L}{4}}|\nabla w|^2\big)^\frac{1}{2}
\quad\mbox{for}\;\frac{L}{4}\ge r\ge \rL_*.
\end{align}
By definition of $w$ and the triangle inequality, this yields
\begin{align*}
\lefteqn{(\av_{B_r}|\nabla(\sigma_i-\sigmaL_i)|^2\big)^\frac{1}{2}
\lesssim(\frac{\rL_*}{L})^\beta+(\av_{B_\frac{L}{4}}|\nabla(\sigma_i-\sigmaL_i)|^2\big)^\frac{1}{2}}\nonumber\\
&+|\xi_j|\Big(\big(\av_{B_r}|\nabla\sigmaL_j|^2\big)^\frac{1}{2}
+\big(\av_{B_\frac{L}{4}}|\nabla\sigmaL_j|^2\big)^\frac{1}{2}\Big)
\quad\mbox{for}\;\frac{L}{4}\ge r\ge \rL_*.
\end{align*}
Inserting the estimate on $|\xi|$ from (\ref{L83}), the localized estimates on
$\sigmaL$ from (\ref{L75}), and the large-scale estimate on $\sigma-\sigmaL$ from (\ref{F20}),
we get the localized estimate on $\sigma-\sigmaL$ stated in (\ref{L81}).

\medskip

It remains to argue that (\ref{L86}) implies (\ref{L87}). Like in the proof of Lemma \ref{L1},
we resort to a construction via a decomposition into dyadic annuli. For any dyadic multiple $r$
of $\rL_*$ with $r\le \frac{L}{4}$ we consider the Lax-Milgram solution $w_r$ of
\begin{align*}
-\Delta w_r=\nabla\cdot(I(B_{2r}-B_r)h);
\end{align*}
the solution of the Poisson equation with rhs $\nabla\cdot(I(B_{\rL_*})h)$ is denoted by $w_{\frac{\rL_*}{2}}$
for notational consistency.
By the energy estimate for the Poisson equation and (\ref{L86}) we have for all dyadic 
$\frac{\rL_*}{2}\le r\le\frac{L}{4}$
\begin{align*}
\big(\frac{1}{r^d}\int|\nabla w_r|^2\big)^\frac{1}{2}\lesssim(\frac{r}{L})^\alpha(\frac{\rL_*}{L})^\beta.
\end{align*}
From this and the (standard) mean-value property (note that unless $r=\frac{\rL_*}{2}$, $w_r$ is harmonic in $B_r$) 
we obtain for every radius $R\ge \rL_*$ (the lower bound arises because of $w_{\frac{\rL_*}{2}}$)
\begin{align*}
\big(\av_{B_R}|\nabla w_r|^2\big)^\frac{1}{2}\lesssim(\frac{r}{L})^\alpha(\frac{\rL_*}{L})^\beta.
\end{align*}
Because of $\alpha>0$ we obtain for $\tilde w:=\sum_{\frac{\rL_*}{2}\le r\le\frac{L}{4}}w_r$ that
\begin{align*}
\big(\av_{B_R}|\nabla\tilde w|^2\big)^\frac{1}{2}\lesssim(\frac{\rL_*}{L})^\beta
\quad\mbox{for}\;R\ge \rL_*.
\end{align*}
Since by construction, $w-\tilde w$ is harmonic in $B_\frac{L}{4}$, we obtain by the (standard) mean-value
property
\begin{align*}
\big(\av_{B_R}|\nabla(w-\tilde w)|^2\big)^\frac{1}{2}
\lesssim\big(\av_{B_\frac{L}{4}}|\nabla(w-\tilde w)|^2\big)^\frac{1}{2}\quad
\mbox{for}\;R\le\frac{L}{4}.
\end{align*}
By the triangle inequality in $L^2$, the two last inequalities yield (\ref{L87}).

\medskip

Equipped with the localized estimates on $\nabla(\phi-\phiL,\sigma-\sigmaL)$ from (\ref{L81}),
it is easy to conclude: Because of (\ref{L80}), it is sufficient to establish (\ref{L85}) for
$\frac{L}{4}\ge r\ge \rL_*$. In fact, we shall show that (\ref{L73}) entails (\ref{L85}) in this range.
This follows instantly from Poincar\'e's inequality in form of
\begin{align*}
\frac{1}{r}\big(\av_{B_r}|(\phi,\sigma)&-\av_{B_r}(\phi,\sigma)|^2\big)^\frac{1}{2}
-\frac{1}{r}\big(\av_{B_r}|(\phiL,\sigmaL)|^2\big)^\frac{1}{2}\nonumber\\
&\lesssim\big(\av_{B_r}|\nabla(\phiL-\phi,\sigmaL-\sigma)|^2\big)^\frac{1}{2},
\end{align*}
into which we plug (\ref{L81}).
\end{proof}

\medskip


\begin{proof}[\sc Proof of Lemma \ref{L6}]\

Since as a finite-range ensemble, $\langle\cdot\rangle$ is in particular ergodic, \cite[Lemma 1]{GNOarXiv}
applies, and yields the existence of $\phi_i$ and $\sigma_i$ with the stated properties.
We set for abbreviation $F_R:=\big(\av_{B_R}|(\phi,\sigma)-\av_{B_R}(\phi,\sigma)|^2\big)^\frac{1}{2}$.

\medskip

For given $1\le s'<2$, we start by extracting the stochastic bound
\begin{align}\label{io01}
\big\langle\exp\big(\frac{F_R}{C\log^\frac{1}{2}R}\big)^{s'}\big\rangle\le C\quad\mbox{for all}\;R\ge 2
\end{align}
from \cite[Corollary 2]{GloriaOttofiniterange}.
Here $C=C(d,\lambda,s')$ denotes a generic constant the value of which may change from line to line.
We focus on the case of $d=2$ (the result is stronger for $d>2$),
in which case the statement of \cite[Corollary 2]{GloriaOttofiniterange} takes the form of
\begin{align*}
\log\big\langle\exp\Big(\big(\frac{1}{C\log(2+|x|)}\int G_1(\cdot-x)\big|(\phi,\sigma)
-\int G_1(\phi,\sigma)\big|^2\big)^\frac{s'}{2}\Big)\big\rangle\le 1
\end{align*}
for all $x\in\mathbb{R}^d$,
where the Gaussian $G_1(x):=(2\pi)^{-\frac{d}{2}}\exp(-\frac{|x|^2}{2})$ plays the role of a spatial averaging function.
We rewrite this 
\begin{align*}
\big\langle\exp\Big(\big(\frac{1}{C\log(2+|x|)}\int G_1(\cdot-x)\big|(\phi,\sigma)
-\int G_1(\phi,\sigma)\big|^2+1\big)^\frac{s'}{2}\Big)\big\rangle\lesssim 1,
\end{align*}
where $\lesssim$ stands for $\le C(d,\lambda,s')$. Since $[0,\infty)\in F\mapsto\exp((F+1)^\frac{s'}{2})$ is convex
for $s'\ge 1$, we obtain for arbitrary $R\ge 2$ from averaging over $x\in B_R$
\begin{align*}
\big\langle\exp\big((\frac{1}{C\log R}\av_{B_R}\big|(\phi,\sigma)
-\int G_1(\phi,\sigma)\big|^2+1)^\frac{s'}{2}\big)\big\rangle\lesssim 1,
\end{align*}
where now $\int G_1(\phi,\sigma)$ may be replaced by $\av_{B_R}(\phi,\sigma)$, so that this turns into (\ref{io01}).

\medskip

We now turn to part i) of the lemma, ie (\ref{io01}), and fix $\beta<1$ and let the random radius
$r_*\in[1,+\infty]$ be minimal with (\ref{L04}).
Because of $\frac{1}{r'}(\av_{B_{r'}}\cdot)^\frac{1}{2}$
$\le(\frac{r}{r'})^{\frac{d}{2}+1}\frac{1}{r}(\av_{B_{r}}\cdot)^\frac{1}{2}$,
$r_*$ is dominated by the smallest {\it dyadic} radius, which for simplicity we call
again $r_*$, with the property that
\begin{align*}
\frac{1}{R}F_R\le 2^{-\frac{d}{2}-1}(\frac{r_*}{R})^\beta\quad\mbox{for all dyadic}\;R>r_*.
\end{align*}
Hence we have for any (deterministic) threshold $r\ge 2$
\begin{align}\label{io02}
\langle I(r_*\ge r)\rangle\le\sum_{R>r\;\mbox{dyadic}}\langle I\big(F_R\ge 2^{-\frac{d}{2}-1} r^{\beta}R^{1-\beta}\big)\rangle,
\end{align}
and obtain from (\ref{io01}) for $s<s'<2$ (say $s'=\frac{s+2}{2}$) by Chebyshev's inequality
\begin{align}\label{io03}
\langle I\big(F_R\ge 2^{-\frac{d}{2}-1} r^{\beta}R^{1-\beta}\big)\rangle
\le \exp\big(-(\frac{r^\beta R^{1-\beta}}{C\log^\frac{1}{2}r})^{s'}+1\big),
\end{align}
where $C=C(d,\lambda,s,\beta)$ denotes a generic constant the value of which may change line by line.
We now claim that $\beta<1$ implies that the first term in the dyadic sum dominates, ie,
\begin{align}\label{io04}
\sum_{R>r\;\mbox{dyadic}}\exp\big(-(\frac{r^\beta R^{1-\beta}}{C\log^\frac{1}{2}r})^{s'}\big)
\lesssim \exp\big(-(\frac{r}{C\log^\frac{1}{2}r})^{s'}\big),
\end{align}
where $\lesssim$ stands short for $\le$ up to a generic multiplicative constant $C=C(d,\lambda,s,\beta)$.
Indeed, setting for abbreviation $C(r):=(\frac{r}{C\log^\frac{1}{2}r})^{s'}$,
this amounts to showing $\sum_{R>r}\exp(-C(r)(\frac{R}{r})^{s'(1-\beta)})$ $\lesssim \exp(-C(r))$,
which because of $C(r)\ge 1$ (for $r\gg 1$, which means that $r$ is larger than some constant only depending
on $d,\lambda,s,\beta$) reduces to the elementary $\sum_{n=1}^\infty\exp(-((2^n)^{s'(1-\beta)}-1))\lesssim 1$,
where we used $\beta<1$ and $s'>0$.
The combination of (\ref{io02}), (\ref{io03}), and (\ref{io04}) yields
\begin{align*}
\langle I(r_*\ge r)\rangle\lesssim\exp\big(-(\frac{r}{C\log^\frac{1}{2}r})^{s'}\big)\quad\mbox{for}\;r\gg 1,
\end{align*}
which implies (\ref{io06}) thanks to $s<s'$.

\medskip

We finally turn to part ii) of the lemma and fix $\beta<1$ and $0<s<2(1-\beta)$. We have to show
the existence of $L_0$ such that
\begin{align*}
\langle I(\exists R\ge L\quad F_R>R^{1-\beta})\rangle\le\exp(-(\frac{L}{L_0})^s)\quad\mbox{for all}\;L\ge L_0;
\end{align*}
by the same argument as for part i) this reduces to
\begin{align}\label{io07}
\lefteqn{\sum_{R\ge L\;\mbox{dyadic}}\langle I(F_R>2^{-\frac{d}{2}-1}R^{1-\beta})\rangle}\nonumber\\
&\le\exp(-(\frac{L}{L_0})^s)\quad\mbox{for all dyadic}\;L\ge L_0.
\end{align}
From (\ref{io01}) for $\frac{s}{1-\beta}<s'<2$ (say $s':=\frac{1}{2}(\frac{s}{1-\beta}+2)$)
we obtain for each summand with $R\ge 2$
\begin{align*}
\langle I\big(F_R\ge 2^{-\frac{d}{2}-1} R^{1-\beta}\big)\rangle
\le \exp\big(-(\frac{R^{1-\beta}}{C\log^\frac{1}{2}R})^{s'}+1\big),
\end{align*}
and as for (\ref{io04}) we find for the sum
\begin{align*}
\sum_{R>L\;\mbox{dyadic}}\exp\big(-(\frac{R^{1-\beta}}{C\log^\frac{1}{2}R})^{s'}\big)
\lesssim \exp\big(-(\frac{L^{1-\beta}}{C\log^\frac{1}{2}L})^{s'}\big),
\end{align*}
Since $s<s'(1-\beta)$ the last two statements imply (\ref{io07}) for some $L_0\lesssim 1$.

\end{proof}

\medskip


\begin{proof}[\sc Proof of Theorem \ref{T1}]\

According to Lemma \ref{L6} ii) and with probability $1-\exp(-(\frac{L}{L_0})^s)$,
the hypothesis (\ref{L80}) of Lemma \ref{L5} is satisfied with $r_*^{(L)}=1$; the second hypothesis
(\ref{L73}) is satisfied by assumption (\ref{F02}) of the theorem. Hence by (\ref{L85}) in Lemma \ref{L5}, (\ref{L04}) holds with
$r_*\lesssim\max\{r_*^{(L)},1\}=r_*^{(L)}$, so that we may apply Proposition \ref{P1} with $r_*^{(L)}$
playing the role of $r_*$. Hence (\ref{L57}) turns into the desired (\ref{F01}).

\end{proof}

\medskip


\begin{proof}[\sc Proof of Lemma \ref{L8}]\

For notational simplicity, we drop the index $i$.
We start by collecting some estimates on $\phi_h'$. Rewriting (\ref{g17}) as
$-\nabla\cdot a_h'\nabla\phi_h'$ $=\nabla\cdot(a_h'-a_h)e$ and noting that the rhs $(a_h'-a_h)e$
is bounded and supported in $B_R(y)$, cf (\ref{g14}), we obtain from the energy estimate 
\begin{align}\label{g20}
\big(\frac{1}{R^d}\int|\nabla\phi_h'|^2\big)^\frac{1}{2}\lesssim 1.
\end{align}
It is convenient to introduce
\begin{align}\label{g21}
v'_h:=x\cdot e+\phi'_h\quad\mbox{so that}\quad\nabla\cdot a_h'\nabla v_h'\stackrel{(\ref{g17})}{=}0,
\end{align} 
and to reformulate (\ref{g20}) as
\begin{align}\label{g22}
\big(\av_{B_{3R}(y)}|\nabla v_h'|^2\big)^\frac{1}{2}\lesssim 1,
\end{align}
which in view of (\ref{g21}), by Meyer's estimate, upgrades to
\begin{align}\label{g23}
\big(\av_{B_{2R}(y)}|\nabla v_h'|^{p}\big)^\frac{1}{p}\lesssim 1\quad\mbox{for some}\;p=p(d,\lambda)>2.
\end{align}
Rewriting (\ref{g17}) once more, this time as $-\nabla\cdot a_h\nabla\phi_h'$ $=\nabla\cdot(a_h-a_h')\nabla v_h'$,
and noting that $a_h$ is a constant coefficient and that $(a_h-a_h')\nabla v_h'$ is supported in $B_R(y)$, cf (\ref{g14}),
we see that $\phi_h'$ decays like the {\it gradient} of the fundamental solution $G_h$. Hence, from (\ref{g20})
we obtain the estimate
\begin{align}\label{g26}
R\sup_{x\in B_{2R}^c(y)}(\frac{|x-y|}{R})^{d+1}|\nabla^2 v_h'|\stackrel{(\ref{g21})}{=}R\sup_{x\in B_{2R}^c(y)}(\frac{|x-y|}{R})^{d+1}|\nabla^2\phi_h'|\lesssim 1.
\end{align}
We may even get closer to the boundary of $B_R^c(y)$ at the expense of a bad constant: For any boundary layer width $0<\rho\le R$ we
obtain from (\ref{g22}) and the fact that $v_h'$ is constant-coefficient harmonic in $B_R^c(y)$, cf (\ref{g21}),
\begin{align}\label{g27}
\sup_{B_{2R}(y)-B_{R+\rho}(y)}\big(\rho|\nabla^2 v_h'|+|\nabla v_h'|\big)\lesssim (\frac{R}{\rho})^\frac{d}{2}.
\end{align}

\medskip

We now turn to the comparison of $\phi'$ and $\phi_h'$, at first in the strong topology on the level of gradients.
We carry this out in terms of the harmonic functions
\begin{align}\label{g37}
v':=e\cdot x+\phi'\quad\mbox{so that}\quad\nabla\cdot a'\nabla v'\stackrel{(\ref{g16})}{=}0
\end{align}
and $v'_h$, by monitoring the error in the two-scale expansion, that is,
\begin{align}\label{g38}
w:=v'-(1+\zeta\phi_i\partial_i)v_h'.
\end{align}
Here $\zeta$ denotes a cut-off function with
\begin{align}\label{g24}
\zeta=0\;\mbox{in}\;B_{R+\rho},\quad\zeta=1\;\mbox{outside}\;B_{R+2\rho},\quad|\nabla\zeta|\lesssim\rho^{-1}
\end{align}
for a boundary layer width $0<\rho\le R$ to be optimized at the end of the proof.
Based on the equations (\ref{g21}), (\ref{g37}), and on (\ref{L01}), we obtain the following formula
\begin{align}\label{g29}
&-\nabla\cdot a'\nabla w=\nabla\cdot(h_{far}+h_{near})\quad\mbox{with}\nonumber\\
&h_{far}:=(\phi_ia-\sigma_i)\nabla(\zeta\partial_iv_h')
\quad\mbox{and}\quad h_{near}:=(1-\zeta)(a'-a_h')\nabla v_h';
\end{align}
which is the same (elementary) calculation as in Step 2 of the proof of \cite[Proposition 1]{GNOarXiv}
and a slight variation of (\ref{L14}).
For the proof of this lemma, we normalize $(\phi,\sigma)$ by $\av_{B_{r_*(y)}(y)}(\phi,\sigma)$ $=0$
so that by (\ref{L32}) (with the origin replaced by $y$) we have
\begin{align}\label{g28}
\frac{1}{r}\big(\av_{B_r(y)}|(\phi,\sigma)|^2\big)^\frac{1}{2}\lesssim(\frac{r_*(y)}{r})^\beta
\quad\mbox{for all}\;r\ge r_*(y).
\end{align}

\medskip

The near-field rhs $h_{near}$ is supported in the thin annulus $B_{R+2\rho}-B_R$ of thickness $2\rho$, 
as a consequence of (\ref{g24}), (\ref{g10}) and (\ref{g14}). With help of Meyer's estimate (\ref{g23}),
we may capitalize on this by H\"older's inequality:
\begin{align}\label{g30}
\big(\frac{1}{R^d}\int|h_{near}|^{2}\big)^\frac{1}{2}\lesssim (\frac{\rho}{R})^{\frac{1}{2}-\frac{1}{p}}.
\end{align}
In the complement of $B_{2R}(y)$, in view of (\ref{g24}) the far-field term assumes the simpler form 
$h_{far}$ $=(\phi_ia-\sigma_i)\nabla\partial_iv_h'$ and thus is easily estimated:
\begin{align}\label{g31}
\lefteqn{\big(\frac{1}{R^d}\int_{B_{2R}^c(y)}|h_{far}|^{2}\big)^\frac{1}{2}}\nonumber\\
&\stackrel{(\ref{g26})}{\lesssim} 
\frac{1}{R}\big(\frac{1}{R^d}\int_{B_{2R}^c(y)}\big((\frac{R}{|x-y|})^{d+1}|(\phi,\sigma)|\big)^{2}\big)^\frac{1}{2}
\stackrel{(\ref{g28})}{\lesssim}(\frac{r_*(y)}{R})^\beta;
\end{align}
the last estimate can be seen by a decomposition into dyadic annuli; here we also use our assumption that $R\ge r_*(y)$.
In $B_{2R}(y)$ by (\ref{g24}) the far-field term is supported in $B_{R+\rho}^c(y)$ and estimated by
$|h_{far}|$ $\lesssim |(\phi,\sigma)|(|\nabla^2 v_h'|+\frac{1}{\rho}|\nabla v_h'|)$. Hence we obtain from (\ref{g28})
and (\ref{g27}) 
\begin{align}\label{g32}
\big(\frac{1}{R^d}\int_{B_{2R}(y)}|h_{far}|^{2}\big)^\frac{1}{2}
\lesssim(\frac{R}{\rho})^{1+\frac{d}{2}}(\frac{r_*(y)}{R})^\beta.
\end{align}
Using (\ref{g30}), (\ref{g31}), and (\ref{g32}) in the energy estimate for (\ref{g29}) we obtain
\begin{align}\label{g33}
\big(\frac{1}{R^d}\int|\nabla w|^{2}\big)^\frac{1}{2}
\lesssim (\frac{\rho}{R})^{\frac{1}{2}-\frac{1}{p}}+(\frac{R}{\rho})^{1+\frac{d}{2}}(\frac{r_*(y)}{R})^\beta.
\end{align}

\medskip

We are interested in the difference of the potentials, cf (\ref{g36}), and in the difference of flux averages, 
cf (\ref{g35}), and thus need to post-process (\ref{g33}).
We first turn to the difference of the potentials and note that by definitions
(\ref{g21}), (\ref{g37}), and (\ref{g38}) of $v_h'$, $v'$ and $w$ we have
$\phi'-\phi_h'$ $=w+\zeta\phi_i\partial_iv_h'$,
so that by Poincar\'e's inequality and the support condition on $\zeta$, cf (\ref{g24}),
\begin{align*}
\lefteqn{\frac{1}{R}\big(\av_{B_{2R}(y)}(\phi'-\phi_h'-\av_{B_{2R}(y)}(\phi'-\phi_h'))^2\big)^\frac{1}{2}}\nonumber\\
&\lesssim \big(\frac{1}{R^d}\int|\nabla w|^{2}\big)^\frac{1}{2}
+\frac{1}{R}\big(\av_{B_{2R}(y)}|\phi|^2\big)^\frac{1}{2}\sup_{B_{2R}(y)-B_{R+\rho}(y)}|\nabla v_h'|,
\end{align*}
so that from (\ref{g33}) for the first rhs term and (\ref{g28}) \& (\ref{g27}) for the second rhs term
we obtain
\begin{align}\label{g39}
\frac{1}{R}\big(\av_{B_{2R}(y)}(\phi'-\phi_h'&-\av_{B_R(y)}(\phi'-\phi_h'))^2\big)^\frac{1}{2}\nonumber\\
&\lesssim (\frac{\rho}{R})^{\frac{1}{2}-\frac{1}{p}}+(\frac{R}{\rho})^{1+\frac{d}{2}}(\frac{r_*(y)}{R})^\beta.
\end{align}

\medskip

We now turn to the differences of flux averages.
The first post-processing step is based on the formula
\begin{align*}
a'\nabla v'-a_h'\nabla v_h'+\nabla\cdot(\zeta\partial_iv_h'\sigma_i)=a'\nabla w+h_{far}+h_{near},
\end{align*}
which in fact is the basis for the formula in (\ref{g29}).
Hence from (\ref{g30}), (\ref{g31}), (\ref{g32}), and (\ref{g33}) we obtain
\begin{align}\label{g34}
\big(\frac{1}{R^d}\int|a'\nabla v'-a_h'\nabla v_h'&+\nabla\cdot(\zeta\partial_iv_h'\sigma_i)|^{2}\big)^\frac{1}{2}\nonumber\\
\lesssim& (\frac{\rho}{R})^{\frac{1}{2}-\frac{1}{p}}+(\frac{R}{\rho})^{1+\frac{d}{2}}(\frac{r_*(y)}{R})^\beta.
\end{align}
The second post-processing step consists in noting that the additional term
$\nabla\cdot(\zeta\partial_iv_h'\sigma_i)$ in the integrand has small average: From integration by parts,
\begin{align*}
\int\omegay\nabla\cdot(\zeta\partial_iv_h'\sigma_i)=-\int\nabla\omegay\cdot(\zeta\partial_iv_h'\sigma_i),
\end{align*}
and since $\omegay$ is supported in $B_{2R}(y)$, cf (\ref{g18bis}), this rhs is bounded by
\begin{align*}
\sup|\nabla\omegay|&\big(\int_{B_{2R}(y)}|\nabla v_h'|^2\big)^\frac{1}{2}\big(\int_{B_{2R}(y)}|\sigma|^2\big)^\frac{1}{2}\nonumber\\
&\stackrel{(\ref{g18bis}),(\ref{g22})}{\lesssim}\frac{1}{R}\big(\av_{B_{2R}(y)}|\sigma|^2\big)^\frac{1}{2}
\stackrel{(\ref{g28})}{\lesssim}(\frac{r_*(y)}{R})^\beta.
\end{align*}
Combining this with (\ref{g34}) yields 
\begin{align}\label{g40}
\big|\int\omegay(a'\nabla v'-a_h'\nabla v_h')\big|
\lesssim (\frac{\rho}{R})^{\frac{1}{2}-\frac{1}{p}}+(\frac{R}{\rho})^{1+\frac{d}{2}}(\frac{r_*(y)}{R})^\beta.
\end{align}

\medskip

The last task is to upgrade (\ref{g39}) and (\ref{g40}) by 
optimizing in the boundary layer thickness $\rho$. Indeed, choosing $\frac{\rho}{R}$ $=(\frac{r_*(y)}{R})^\frac{\beta}{\frac{1}{2}-\frac{1}{p}+1+\frac{d}{2}}$ 
we obtain (\ref{g36}) and (\ref{g35}) with the potentially small but positive exponent
$\alpha=\beta\frac{\frac{1}{2}-\frac{1}{p}}{\frac{1}{2}-\frac{1}{p}+1+\frac{d}{2}}$. 
\end{proof}

\medskip


\begin{proof}[\sc Proof of Corollary \ref{Le10}]\

We start by noting that
\begin{align}\label{g47}
\big(\av_{B_{2R}(y)}|e_j+\nabla\phi_j'|^2\big)^\frac{1}{2}\lesssim 1.
\end{align}
Indeed, by (\ref{L55}) in form of
\begin{align}\label{g45}
\big(\av_{B_{2R}(y)}|e_j+\nabla\phi_j|^2\big)^\frac{1}{2}\lesssim 1
\end{align}
and the triangle inequality in $L^2$, it suffices to show
\begin{align}\label{g46}
\big(\av_{B_{2R}(y)}|\nabla(\phi_j'-\phi_j)|^2\big)^\frac{1}{2}\lesssim 1.
\end{align}
Since $-\nabla\cdot a'\nabla(\phi_j'-\phi_j)=\nabla\cdot(a'-a)(e_j+\nabla\phi_j)$, cf (\ref{L31}) and (\ref{g16}),
and since $a'-a$ is supported in $B_R(y)$, cf (\ref{g10}), we have by the energy estimate
\begin{align*}
\big(\frac{1}{R^d}\int|\nabla(\phi_j'-\phi_j)|^2\big)^\frac{1}{2}\lesssim \big(\av_{B_{2R}(y)}|e_j+\nabla\phi_j|^2\big)^\frac{1}{2},
\end{align*}
so that (\ref{g46}) is a consequence of (\ref{g45}).

\medskip

Since $a'-a$ and $a_h'-a_h$ are supported in $B_R(y)$, cf (\ref{g10}) and (\ref{g14}),
we may smuggle in the averaging function $\omegay$, cf (\ref{g18bis}), into the
lhs of the desired (\ref{g41}) so that it is enough to show
\begin{align}\label{g42}
\lefteqn{\big|\int\omegay(e_j+\nabla\phi_j')\cdot(a'-a)(e_i+\nabla\phi_i)}\nonumber\\
&-\int\omegay(e_j+\nabla\phi_{jh}')\cdot(a'_h-a_h)e_i\big|\lesssim(\frac{r_*(y)}{R})^\alpha.
\end{align}
By integration by parts, and using (\ref{L31}) and (\ref{g16}), we obtain for the first integral in (\ref{g42})
\begin{align}\label{g43}
\lefteqn{\int\omegay(e_j+\nabla\phi_j')\cdot(a'-a)(e_i+\nabla\phi_i)}\nonumber\\
&=\int\omegay(e_j+\nabla\phi_j')\cdot a'e_i
-\int\omegay e_j\cdot a(e_i+\nabla\phi_i)\nonumber\\
&-\int(\phi_i-\av_{B_{2R}(y)}\phi_i)(e_j+\nabla\phi_j')\cdot a'\nabla\omegay\nonumber\\
&+\int(\phi_j'-\av_{B_{2R}(y)}\phi_j')\nabla\omegay\cdot a(e_i+\nabla\phi_i).
\end{align}
On the first rhs term in (\ref{g43}) we apply (\ref{g35}) of Lemma \ref{L8}:
\begin{align}\label{g50}
\big|\int\omegay(e_j+\nabla\phi_j')\cdot a'e_i-\int\omegay(e_j+\nabla\phi_{jh}')\cdot a'_he_i\big|\lesssim(\frac{r_*(y)}{R})^\alpha.
\end{align}
Using (\ref{L01}), the second rhs term in (\ref{g43}) (without the minus sign) can be rewritten as
\begin{align*}
\int\omegay e_j\cdot a_h e_i-\int e_j\cdot(\sigma_i-\av_{B_{2R}(y)}\sigma_i)\nabla\omegay,
\end{align*}
where the second contribution is estimated as follows
\begin{align*}
\big|\int e_j\cdot(\sigma_i-\av_{B_{2R}(y)}\sigma_i)\nabla\omegay\big|
\stackrel{(\ref{g18bis})}{\lesssim}\frac{1}{R}\big(\av_{B_{2R}(y)}|\sigma_i-\av_{B_{2R}(y)}\sigma_i|^2\big)^\frac{1}{2}
\stackrel{(\ref{L04})}{\lesssim}(\frac{r_*(y)}{R})^{\beta},
\end{align*}
so that we obtain
\begin{align}\label{g51}
\big|\int\omegay e_j\cdot a(e_i+\nabla\phi_i)-\int\omegay e_j\cdot a_h e_i\big|\lesssim(\frac{r_*(y)}{R})^\beta.
\end{align}
The third rhs term in (\ref{g43}) is estimated as follows
\begin{align}\label{g52}
\lefteqn{\big|\int(\phi_i-\av_{B_{2R}(y)}\phi_i)(e_j+\nabla\phi_j')\cdot a'\nabla\omegay\big|}\nonumber\\
&\stackrel{(\ref{g18bis})}{\lesssim}\frac{1}{R}\big(\av_{B_{2R}(y)}(\phi_i-\av_{B_{2R}(y)}\phi_i)^2\big)^\frac{1}{2}
\big(\av_{B_{2R}(y)}|e_j+\nabla\phi_j'|^2\big)^\frac{1}{2}
&\stackrel{(\ref{L04}),(\ref{g47})}{\lesssim}(\frac{r_*(y)}{R})^{\beta}.
\end{align}
We now turn to the last rhs term in (\ref{g43}). We first apply (\ref{g36}) in Lemma \ref{L8}
to the effect of
\begin{align}\label{g48}
\lefteqn{\big|\int(\phi_j'-\av_{B_{2R}(y)}\phi_j')\nabla\omegay\cdot a(e_i+\nabla\phi_i)
-\int(\phi_{jh}'-\av_{B_{2R}(y)}\phi_{jh}')\nabla\omegay\cdot a(e_i+\nabla\phi_i)\big|}\nonumber\\
&\stackrel{(\ref{g18bis})}{\lesssim}\frac{1}{R}
\Big(\av_{B_{2R}(y)}\big(\phi_j'-\phi_{jh}'-\av_{B_{2R}(y)}(\phi_j'-\phi_{jh}')\big)^2\Big)^\frac{1}{2}
\big(\av_{B_{2R}(y)}|e_i+\nabla\phi_i|^2\big)^\frac{1}{2}\nonumber\\
&\stackrel{(\ref{g36}),(\ref{L55})}{\lesssim}(\frac{r_*(y)}{R})^{\alpha}.
\end{align}
We then note that by (\ref{L01}), the skew symmetry of $\sigma_i$ and two integration by parts
\begin{align*}
\lefteqn{\int(\phi_{jh}'-\av_{B_{2R}(y)}\phi_{jh}')\nabla\omegay\cdot a(e_i+\nabla\phi_i)}\nonumber\\
&=-\int\omegay\nabla\phi_{jh}'\cdot a_he_i-\int\nabla\omegay\cdot(\sigma_i-\av_{B_{2R}(y)}\sigma_i)\nabla\phi_{jh}',
\end{align*}
so that
\begin{align}\label{g49}
\lefteqn{\big|\int(\phi_{jh}'-\av_{B_{2R}(y)}\phi'_{jh})\nabla\omegay\cdot a(e_i+\nabla\phi_i)
+\int\omegay\nabla\phi_{jh}'\cdot a_he_i\big|}\nonumber\\
&\stackrel{(\ref{g18bis})}{\lesssim}\frac{1}{R}\big(\av_{B_{2R}(y)}|\sigma_i-\av_{B_{2R}(y)}\sigma_i|^2\big)^\frac{1}{2}
\big(\av_{B_{2R}(y)}|\nabla\phi_{jh}'|^2\big)^\frac{1}{2}
\stackrel{(\ref{L04}),(\ref{g20})}{\lesssim}(\frac{r_*(y)}{R})^{\beta}.
\end{align}
The combination of (\ref{g48}) and (\ref{g49}) yields for the last term in (\ref{g43}):

\begin{align}\label{g53}
\big|\int(\phi_j'-\av_{B_{2R}(y)}\phi_{j}')
\nabla\omegay\cdot a(e_i+\nabla\phi_i)+\int\omegay\nabla\phi_{jh}'\cdot a_he_i\big|
\lesssim(\frac{r_*(y)}{R})^{\alpha}.
\end{align}
Inserting the four estimates (\ref{g50}), (\ref{g51}), (\ref{g52}) and (\ref{g53}) into (\ref{g43}),
we obtain (\ref{g42}).
\end{proof}

\medskip


\begin{proof}[\sc Proof of Lemma \ref{L9}]\

Starting point is Lemma \ref{L1}, more precisely (\ref{L17}) for $L:=\frac{1}{2}|y|$ and (\ref{L62}) in form of
\begin{align}\label{g57}
\big(\av_{B_{L}^c}|\nabla(u-(1+(\phi_i-\av_{B_{r_*}}\phi_i)\partial_i)u_h)|^2\big)^\frac{1}{2}
\lesssim (\frac{\ell}{|y|})^d (\frac{r_*}{|y|})^{\beta}.
\end{align}
The first post-processing step is to replace $\av_{B_{r_*}}\phi_i$ by $\av_{B_{r_*}(y)}\phi_i$
in (\ref{g57}). Indeed, by (\ref{g58}) (for $r=L=\frac{1}{2}|y|$ and $r'=r_*$) we have
\begin{align}\label{g59}
 \frac{1}{L}\big|\av_{B_L   }(\phi,\sigma)&-\av_{B_{r_*   }   }(\phi,\sigma)\big|
+\frac{1}{L}\big|\av_{B_L(y)}(\phi,\sigma) -\av_{B_{r_*(y)}(y)}(\phi,\sigma)\big|\nonumber\\
&\lesssim (\frac{\max\{r_*,r_*(y)\}}{|y|})^\beta.
\end{align}
As a direct consequence of (\ref{L04}) (for $r=3L$) and $B_L\cup B_L(y)\subset B_{3L}$ we have 
\begin{align*}
\frac{1}{L}\big|\av_{B_L}(\phi,\sigma)-\av_{B_L(y)}(\phi,\sigma)\big|
\lesssim(\frac{r_*}{|y|})^\beta,
\end{align*}
so that (\ref{g59}) implies
\begin{align*}
 \frac{1}{L}\big|\av_{B_{r_*}}(\phi,\sigma)-\av_{B_{r_*(y)}(y)}(\phi,\sigma)\big|
\lesssim (\frac{\max\{r_*,r_*(y)\}}{|y|})^\beta.
\end{align*}
We combine this with (\ref{L06}) in form of
\begin{align}\label{g65}
\sup_{B_L^c}|\nabla u_h|+\sup_{B_L^c}L|\nabla^2u_h|\lesssim(\frac{\ell}{|y|})^d
\end{align}
to the desired
\begin{align*}
\big(\av_{B_{L}^c}|\nabla(u-(1+(\phi_i-\av_{B_{r_*(y)}(y)}\phi_i)\partial_i)u_h)|^2\big)^\frac{1}{2}
\lesssim (\frac{\ell}{|y|})^d (\frac{\max\{r_*,r_*(y)\}}{|y|})^{\beta}.
\end{align*}
Hence wlog we assume $\av_{B_{r_*(y)}(y)}\phi_i=0$ so that on the one hand, the above simplifies to
\begin{align}\label{g60}
\big(\av_{B_{L}^c}|\nabla(u-(1+\phi_i\partial_i)u_h)|^2\big)^\frac{1}{2}
\lesssim (\frac{\ell}{|y|})^d (\frac{\max\{r_*,r_*(y)\}}{|y|})^{\beta}
\end{align}
and (\ref{L32}), with the origin replaced by $y$, assumes the form of
\begin{align}\label{g61}
\frac{1}{r}\big(\av_{B_r(y)}|(\phi,\sigma)|^2\big)^\frac{1}{2}\lesssim(\frac{r_*(y)}{r})^\beta
\quad\mbox{for all}\;r\ge r_*(y).
\end{align}

\medskip

Following Step 6 in the proof of \cite[Theorem 0.2]{BellaGiuntiOttoPCMINotes}, we now localize (\ref{g60}) around $y$, 
making use of (\ref{g61}).
To this purpose, we appeal once more to the formula (\ref{L14}) for the error in the two-scale convergence
\begin{align*}
-\nabla\cdot a\nabla(u-(1+\phi_i\partial_i)u_h)=\nabla\cdot h\quad\mbox{where}\quad
h:=(\phi_ia-\sigma_i)\nabla\partial_iu_h.
\end{align*}
The combination of (\ref{g65}) and (\ref{g61}) yields for the rhs
\begin{align}\label{g63}
\big(\av_{B_r(y)}|h|^2\big)^\frac{1}{2}\lesssim(\frac{\ell}{|y|})^d\frac{r}{L}(\frac{r_*(y)}{r})^\beta
\quad\mbox{for all}\;L\ge r\ge r_*(y).
\end{align}
We now argue, starting from the large-scale anchoring of (\ref{g60}) in form of
\begin{align}\label{g69}
\big(\av_{B_L(y)}|\nabla(u-(1+\phi_i\partial_i)u_h)|^2\big)^\frac{1}{2}
\lesssim (\frac{\ell}{|y|})^d (\frac{\max\{r_*,r_*(y)\}}{|y|})^{\beta},
\end{align}
that (\ref{g63}) allows for the desired localization to our scale of interest $R\ge r_*(y)$:
\begin{align}\label{g70}
\big(\av_{B_R(y)}|\nabla(u-(1+\phi_i\partial_i)u_h)|^2\big)^\frac{1}{2}
\lesssim(\frac{\ell}{|y|})^d(\frac{\max\{r_*,r_*(y)\}}{|y|})^\beta.
\end{align}
To this purpose, for every radius $r=R,2R,\cdots$ with $r\le L$  
we consider the Lax-Milgram solution of
\begin{align*}
-\nabla\cdot a\nabla u_r=\nabla\cdot(I(B_{r}(y)-B_\frac{r}{2}(y))h)
\end{align*}
with the understanding that for $r=R$, the rhs is given by $I(B_R(y))h$. From the
energy estimate and (\ref{g63}) we obtain 
\begin{align}\label{g66}
\big(\frac{1}{r^d}\int|\nabla u_r|^2\big)^\frac{1}{2}
\lesssim(\frac{\ell}{|y|})^d\frac{r}{L}(\frac{r_*(y)}{r})^\beta.
\end{align}
Since for $r\ge 2R$, $u_r$ is $a$-harmonic in $B_\frac{r}{2}(y)$, we may apply the $C^{0,1}$-estimate, cf (\ref{io10}),
to localize the above to
\begin{align*}
\big(\av_{B_{R}(y)}|\nabla u_r|^2\big)^\frac{1}{2}
\lesssim(\frac{\ell}{|y|})^d\frac{r}{L}(\frac{r_*(y)}{r})^\beta.
\end{align*}
Since the exponent $1-\beta$ of $r$ is positive, $\tilde u:=\sum_{L\ge r\ge R}u_r$ satisfies
\begin{align}\label{g69bis}
\big(\av_{B_{R}(y)}|\nabla\tilde u|^2\big)^\frac{1}{2}
\lesssim(\frac{\ell}{|y|})^d(\frac{r_*(y)}{|y|})^\beta.
\end{align}
Likewise, we obtain directly from (\ref{g66})
\begin{align}\label{g67}
\big(\frac{1}{L^d}\int|\nabla\tilde u|^2\big)^\frac{1}{2}
\lesssim(\frac{\ell}{|y|})^d(\frac{r_*(y)}{|y|})^\beta.
\end{align}
Since by construction, $u-(1+\phi_i\partial_i)u_h-\tilde u$ 
is harmonic in $B_\frac{L}{2}(y)$, we may apply the $C^{0,1}$-estimate
to the effect of
\begin{align}\label{g68}
\lefteqn{\big(\av_{B_R(y)}|\nabla(u-(1+\phi_i\partial_i)u_h-\tilde u)|^2\big)^\frac{1}{2}}\nonumber\\
&\lesssim\big(\av_{B_L(y)}|\nabla(u-(1+\phi_i\partial_i)u_h-\tilde u)|^2\big)^\frac{1}{2}.
\end{align}
By the triangle inequality in $L^2$, we obtain (\ref{g70}) from 
(\ref{g68}), (\ref{g69bis}), (\ref{g67}), and (\ref{g69}).

\medskip

We now post-process (\ref{g70}), making use of $r_*,r_*(y)\le R$, to
\begin{align}\label{g77}
\big(\av_{B_R(y)}|\nabla u-\partial_i u_h(y)(e_i+\nabla\phi_i)|^2\big)^\frac{1}{2}
\lesssim(\frac{\ell}{|y|})^d(\frac{R}{|y|})^\beta
\end{align}
and also note for later purpose
\begin{align}\label{g78}
\big(\av_{B_R(y)}|\nabla u|^2\big)^\frac{1}{2},
\big(\av_{B_R(y)}|\partial_i u_h(y)(e_i+\nabla\phi_i)|^2\big)^\frac{1}{2}
\lesssim(\frac{\ell}{|y|})^d.
\end{align}
In deriving (\ref{g77}) from (\ref{g70}), 
we first replace $\nabla(1+\phi_i\partial_i)u_h$ by $\partial_iu_h(e_i+\nabla\phi_h)$,
which means that by the triangle inequality in $L^2$, we have to estimate $\phi_i\nabla\partial_iu_h$.
We have
\begin{align*}
\big(\av_{B_R(y)}|\phi_i\nabla\partial_iu_h|^2\big)^\frac{1}{2}
\le\sup_{B_R(y)}|\nabla^2u_h|\big(\av_{B_R(y)}|\phi|^2\big)^\frac{1}{2};
\end{align*}
by $B_R(y)\subset B_L^c$ and (\ref{g65}), the first factor is estimated by $\frac{1}{|y|}(\frac{\ell}{|y|})^d$.
By (\ref{g61}), the second factor is estimated by $R(\frac{r_*(y)}{R})^\beta$.
Hence because of $R\le|y|$, this contribution is contained in the rhs of (\ref{g77}).
We now replace $\partial_iu_h(e_i+\nabla\phi_h)$ by $\partial_iu_h(y)(e_i+\nabla\phi_h)$. By
the triangle inequality in $L^2$, we are lead to estimating 
\begin{align*}
\lefteqn{\big(\av_{B_R(y)}|(\partial_iu_h-\partial_iu_h(y))(e_i+\nabla\phi_i)|^2\big)^\frac{1}{2}}\nonumber\\
&\le R\sup_{B_R(y)}|\nabla\partial_iu_h|\big(\av_{B_R(y)}|e_i+\nabla\phi_i|^2\big)^\frac{1}{2}.
\end{align*}
As above, the first factor is estimated by $\frac{R}{|y|}(\frac{\ell}{|y|})^d$. According to
(\ref{L55}) (with $y$ playing the role of the origin and using $R\ge r_*(y)$), the
second factor is estimated by $1$. Because of $\beta\le 1$, 
also this contribution is contained in the rhs of (\ref{g77}). The same argument leads to the
second estimate in (\ref{g78}). The first estimate in (\ref{g78}) follows from the second
one and (\ref{g77}) via the triangle inequality in $L^2$.

\medskip

We now turn to the sensitivity estimate and consider
\begin{align}\label{g75}
w:=u'-u,\quad\mbox{st}\quad-\nabla\cdot a'\nabla w=\nabla\cdot h'\quad\mbox{with}\quad
h':=(a'-a)\nabla u.
\end{align}
We want to apply the localized homogenization error estimate of \cite[Theorem 0.2]{BellaGiuntiOttoPCMINotes}, 
for the medium given by $a'$. This means comparing $w$ to the solution of 
\begin{align}\label{g87}
-\nabla\cdot a_h\nabla w_h=\nabla\cdot (h_i'(e_i+\partial_i\phi')).
\end{align}
Note that by assumption, and by Lemma \ref{L7} we have $r_*'(y),r_*'\lesssim R$.
Since $h'$ is supported in $B_R(y)$, by \cite[Theorem 0.2]{BellaGiuntiOttoPCMINotes} 
(for the medium $a'$, $R$ playing the role of $r_*$ there, the roles of $y$ and the origin exchanged,
and $\beta$ replacing $1-\alpha$ there)  we have
\begin{align*}
\big(\av_{B_R}|\nabla w-\partial_kw_h(e_k+\nabla\phi_k')|^2\big)^\frac{1}{2}
\lesssim(\frac{R}{|y|})^{d+\beta}\big(\av_{B_R(y)}|h'|^2\big)^\frac{1}{2}.
\end{align*}
(As explained in the proof of Lemma \ref{L1}, the logarithm in \cite[Theorem 0.2]{BellaGiuntiOttoPCMINotes}
can be avoided.)
Among other ingredients, this estimate is based on the following estimate of $w_h$,
cf Step 2 of the proof of \cite[Theorem 0.2]{BellaGiuntiOttoPCMINotes},
\begin{align*}
\sup_{B_L^c(y)}(|\nabla w_h|+L|\nabla^2 w_h|)
\lesssim(\frac{R}{L})^{d}\big(\av_{B_R(y)}|h'|^2\big)^\frac{1}{2}.
\end{align*}
In view of the definition of $h'$, cf (\ref{g75}), and (\ref{g78}), these two estimates turn into
\begin{align}
\big(\av_{B_R}|\nabla w-\partial_kw_h(e_k+\nabla\phi_k')|^2\big)^\frac{1}{2}
\lesssim(\frac{\ell}{|y|})^d(\frac{R}{|y|})^{d+\beta},\label{g76}\\
\sup_{B_L^c(y)}(|\nabla w_h|+L|\nabla^2 w_h|)
\lesssim(\frac{\ell}{|y|})^{d}(\frac{R}{L})^{d}.\label{g79}
\end{align}
Like for the passage from (\ref{g70}) to (\ref{g77}), (\ref{g76}) may be post-processed 
with help of (\ref{g79}) to
\begin{align}\label{g80}
\big(\av_{B_R}|\nabla w-\partial_kw_h(0)(e_k+\nabla\phi'_k)|^2\big)^\frac{1}{2}
\lesssim(\frac{\ell}{|y|})^d(\frac{R}{|y|})^{d+\beta}.
\end{align}
In this argument, we just have to replace the medium $a$ by the medium $a'$ and $y$ by the origin.
%
%

\medskip

We continue with post-processing and argue that we may replace $\phi_k'$ by $\phi_k$ in (\ref{g80}):
\begin{align}\label{g82}
\big(\av_{B_R}|\nabla w-\partial_kw_h(0)(e_k+\nabla\phi_k)|^2\big)^\frac{1}{2}
\lesssim(\frac{\ell}{|y|})^d(\frac{R}{|y|})^{d+\beta}.
\end{align}
In fact, we claim that the error term is of (substantially) higher order:
\begin{align*}
|\partial_kw_h(0)|\big(\av_{B_R}|\nabla(\phi_k'-\phi_k)|^2\big)^\frac{1}{2}
\lesssim(\frac{\ell}{|y|})^d(\frac{R}{|y|})^{2d},
\end{align*}
which by (\ref{g79}) reduces to
\begin{align}\label{g84}
\big(\av_{B_R}|\nabla(\phi_k'-\phi_k)|^2\big)^\frac{1}{2}
\lesssim(\frac{R}{|y|})^{d}.
\end{align}
The latter can be seen noting that 
\begin{align}\label{g83}
-\nabla\cdot a\nabla(\phi_k'-\phi_k)=\nabla\cdot(a'-a)(e_k+\nabla\phi_k').
\end{align}
From (\ref{g83}) we learn at first that $\phi_k'-\phi_k$ is $a$-harmonic in $B_R^c(y)$ so
that by the $C^{0,1}$-estimate, cf (\ref{io10}), and since $R\ge r_*$
\begin{align*}
\big(\av_{B_R}|\nabla(\phi_k'-\phi_k)|^2\big)^\frac{1}{2}\lesssim
\big(\av_{B_L}|\nabla(\phi_k'-\phi_k)|^2\big)^\frac{1}{2},
\end{align*}
where we recall our abbreviation $L=\frac{1}{2}|y|$. Moreover, since the rhs of (\ref{g83}) is in
divergence form, we have by the dualized $C^{0,1}$-estimate (see Step 5 in the proof of 
\cite[Theorem 0.2]{BellaGiuntiOttoPCMINotes} for such a duality argument) and since $R\ge r_*(y)$
\begin{align*}
\big(\frac{1}{L^d}\int_{B_L^c(y)}|\nabla(\phi_k'-\phi_k)|^2\big)^\frac{1}{2}
\lesssim(\frac{R}{L})^d\big(\frac{1}{R^d}\int_{B_R^c(y)}|\nabla(\phi_k'-\phi_k)|^2\big)^\frac{1}{2}.
\end{align*}
Finally, by the energy estimate for (\ref{g83}) we have
\begin{align*}
(\frac{1}{R^d}\int|\nabla(\phi_k'-\phi_k)|^2\big)^\frac{1}{2}
\lesssim(\av_{B_R(y)}|e_k+\nabla\phi_k'|^2\big)^\frac{1}{2}\lesssim 1,
\end{align*}
where in the very last estimate, we've used (\ref{L55}) (with $(a,0)$ replaced by $(a',y)$),
which we may since by Lemma \ref{L7} we have $r_*'(y)\lesssim R$. Since by definition of $L$,
$B_L\subset B_L^c(y)$ the last three estimates combine to (\ref{g84}).

\medskip

In the remainder of the proof we argue that we may replace $w_h$ by a more
explicit expression. More precisely, in order to pass from (\ref{g82}) to (\ref{g85}) we have to
replace $\partial_kw_h(0)$ by
\begin{align}\label{g88}
\partial_j\partial_kG_h(-y)\partial_iu_h(y)\int(e_j+\nabla\phi_j')\cdot(a'-a)(e_i+\nabla\phi_i),
\end{align}
and then appeal to the definition (\ref{g55}) of $\delta a_{ij}$. The basis for this is the representation
\begin{align}\label{g89}
\partial_kw_h(0)=\int\partial_j\partial_kG_h(-x)\big((e_j+\nabla\phi_j')\cdot(a'-a)\nabla u\big)(x)dx,
\end{align}
which is a consequence of the definition (\ref{g87}) of $w_h$, yielding the representation
\begin{align*}
\partial_kw_h(0)&=\int\nabla\partial_kG_h(-x)\cdot\big(h_i'(e_i+\partial_i\phi')))(x)dx\\
&=\int\partial_j\partial_kG_h(-x)\big(h'\cdot(e_j+\nabla\phi_j')))(x)dx,
\end{align*}
into which we insert the definition (\ref{g75}) of $h'$. We split the passage from (\ref{g89}) to
(\ref{g88}) into two steps. In view of $(\av_{B_R}|e_k+\nabla\phi_k|^2)^\frac{1}{2}\lesssim 1$,
which is a consequence of $R\ge r_*$, cf (\ref{L55}), it is enough to show
\begin{align*}
R^d\sup_{x\in B_R(y)}|\partial_j\partial_kG_h(-x)-\partial_j\partial_kG_h(-y)|
\av_{B_R}|e_j+\nabla\phi_j'||\nabla u|&\lesssim(\frac{\ell}{|y|})^d(\frac{R}{|y|})^{d+1},\\
R^d|\partial_j\partial_kG_h(-y)|
\av_{B_R}|e_j+\nabla\phi_j'||\nabla u-\partial_iu_h(y)(e_i+\nabla\phi_i)|
&\lesssim(\frac{\ell}{|y|})^d(\frac{R}{|y|})^{d+\beta},
\end{align*}
which by the obvious decay properties of $G_h$ reduces to
\begin{align*}
\av_{B_R}|e_j+\nabla\phi_j'||\nabla u|&\lesssim(\frac{\ell}{|y|})^d,\\
\av_{B_R}|e_j+\nabla\phi_j'||\nabla u-\partial_iu_h(y)(e_i+\nabla\phi_i)|
&\lesssim(\frac{\ell}{|y|})^d(\frac{R}{|y|})^{\beta}.
\end{align*}
By the Cauchy-Schwarz inequality, because of $(\av_{B_R}|e_j+\nabla\phi_j'|^2)^\frac{1}{2}\lesssim 1$,
which is a consequence of $R\gtrsim r_*'$, this reduces to (\ref{g77}) and (\ref{g78}).
\end{proof}

\medskip


\begin{proof}[\sc Proof of Proposition \ref{P2}]\

We may apply Lemma \ref{L9} to the heterogeneous medium $a$ replaced by the homogeneous $a_h$;
in which case we may choose $(\phi,\sigma)\equiv 0$ and thus in particular $r_*=r_*(y)=0$, whereas in view of
(\ref{g16}) and (\ref{g17}), $\phi'$ is being replaced by $\phi_h'$ and thus $\delta a$ by $\delta a_h$, cf (\ref{g96}). 
An inspection of the proof of Lemma \ref{L9} shows that me may replace $w=u'-u$ by $w_h:=u_h'-u_h$,
since the crucial property of that difference was, when translated to the homogeneous medium,
\begin{align*}
w_h:=u'_h-u_h,\quad\mbox{st}\quad-\nabla\cdot a'_h\nabla w_h=\nabla\cdot h'_h\quad\mbox{with}\quad
h'_h:=(a'_h-a_h)\nabla u_h,
\end{align*}
cf (\ref{g75}) for $w$, which for $w_h$ follows immediately from the 
definition (\ref{g12}) of $u_h'$. Hence (\ref{g85}) assumes the form
\begin{align*}
\big(\av_{B_R}|\nabla(u'_h-u_h)
-\partial_iu_h(y)\delta a_{hij}\partial_j\partial_kG_h(-y)e_k|^2\big)^\frac{1}{2}
\lesssim(\frac{\ell}{|y|})^d(\frac{R}{|y|})^{d+\beta}.
\end{align*}
Into this estimate, we insert (\ref{g41}) in form of $\frac{1}{R^d}|\delta a_{ij}-\delta a_{hij}|$
$\lesssim(\frac{r_*(y)}{R})^\alpha$, combined with (\ref{g96bis}) below to
\begin{align*}
\big(\av_{B_R}|\partial_k(u'_h-u_h)&
-\partial_iu_h(y)\delta a_{ij}\partial_j\partial_kG_h(-y)|^2\big)^\frac{1}{2}\nonumber\\
&\lesssim(\frac{\ell}{|y|})^d(\frac{R}{|y|})^{d}\big((\frac{R}{|y|})^{\beta}+(\frac{r_*(y)}{R})^\alpha\big).
\end{align*}
In combination with (\ref{L55}) in form of $\big(\av_{B_R}|e_k+\nabla\phi_k|^2\big)^\frac{1}{2}\lesssim 1$,
this assumes the form
\begin{align*}
\big(\av_{B_R}|\partial_k(u'_h-u_h)(e_k+\nabla\phi_k)&
-\partial_iu_h(y)\delta a_{ij}\partial_j\partial_kG_h(-y)(e_k+\nabla\phi_k)|^2\big)^\frac{1}{2}\nonumber\\
&\lesssim(\frac{\ell}{|y|})^d(\frac{R}{|y|})^{d}\big((\frac{R}{|y|}^{\beta}+(\frac{r_*(y)}{R})^\alpha\big).
\end{align*}
By the triangle inequality in $L^2$ with (\ref{g85}) in its original form, we obtain (\ref{g15}).
\end{proof}

\medskip

\begin{proof}[\sc Proof of Lemma \ref{L7}]\

Clearly, it suffices to establish (\ref{g01}) at the origin (which here plays the role of the general point), that is,
\begin{align*}
r_*'\lesssim r_*+r_*(y)+R.
\end{align*}
According to the definition of $r_*$ as the minimal radius with the property (\ref{L04}), 
by the triangle inequality in $L^2$, it is enough to show
\begin{align*}
\frac{1}{r}\big(\av_{B_r}|(\phi'-\phi,\sigma'-\sigma)&-\av_{B_r}(\phi'-\phi,\sigma'-\sigma)|^2\big)^\frac{1}{2}\nonumber\\
&\lesssim(\frac{r_*(y)+R}{r})^\beta\quad\mbox{for}\;r\ge r_*(y)+R,
\end{align*}
which by Poincar\'e's inequality follows from
\begin{align*}
\big(\av_{B_r}|\nabla(\phi'-\phi,\sigma'-\sigma)|^2\big)^\frac{1}{2}
\lesssim(\frac{r_*(y)+R}{r})^\beta\quad\mbox{for}\;r\ge r_*(y)+R.
\end{align*}
In view of $\beta<1\le\frac{d}{2}$ and of $(\av_{B_r}f^2)^\frac{1}{2}$ 
$\lesssim(\frac{R}{r})^\frac{d}{2}(\frac{1}{R^d}\int f^2)^\frac{1}{2}$
this follow from
\begin{align}\label{g03}
\big(\frac{1}{R^d}\int|\nabla(\phi'-\phi,\sigma'-\sigma)|^2\big)^\frac{1}{2}
\lesssim (1+\frac{r_*(y)}{R})^\frac{d}{2},
\end{align}
which we shall establish now.

\medskip

In order to tackle (\ref{g03}), we fix a coordinate direction $i=1,\cdots,d$ and suppress the index $i$ in $(\phi_i,\sigma_i)$
and $(\phi_i',\sigma_i')$. Since $-\nabla\cdot a'\nabla(\phi'-\phi)$ $=\nabla\cdot(a'-a)(e+\nabla\phi)$, we have
by the energy inequality
\begin{align}\label{g04}
\big(\frac{1}{R^d}\int|\nabla(\phi'-\phi)|^2\big)^\frac{1}{2}
\lesssim\big(\av_{B_R(y)}|e+\nabla\phi|^2\big)^\frac{1}{2}.
\end{align}
We now note that Caccioppoli's estimate (\ref{L55}) implies
\begin{align}\label{l01}
\big(\av_{B_R(y)}|e+\nabla\phi|^2\big)^\frac{1}{2}\lesssim (1+\frac{r_*(y)}{R})^\frac{d}{2};
\end{align}
This is immediate in case of $R\ge r_*(y)$ and follows from 
$(\av_{B_R(y)}f^2)^\frac{1}{2}$ $\le(\frac{R}{r_*(y)})^\frac{d}{2}$ $(\av_{B_{r_*(y)}}f^2)^\frac{1}{2}$
in the other case. The combination of (\ref{g04}) and (\ref{l01}) yields (\ref{g03})
for the $\phi$-part.

\medskip

We now turn to the $\sigma$-part of (\ref{g03}) and note that $-\Delta(\sigma_{jk}'-\sigma_{jk})$ $=\partial_j(q'-q)_k$
$-\partial_k(q'-q)_j$ where $q'-q$ $=a'\nabla(\phi'-\phi)$ $+(a'-a)(e+\nabla\phi)$. Hence we obtain by
the energy estimate
\begin{align*}
\big(\frac{1}{R^d}\int|\nabla(\sigma'-\sigma)|^2\big)^\frac{1}{2}
\lesssim\big(\frac{1}{R^d}\int|\nabla(\phi'-\phi)|^2\big)^\frac{1}{2}+\big(\av_{B_R(y)}|e+\nabla\phi|^2\big)^\frac{1}{2},
\end{align*}
so that it remains to appeal to (\ref{g04}) and (\ref{l01}) once more.
\end{proof}

\medskip


\begin{proof}[\sc Proof of Corollary \ref{Le11}]\

\ignore{
Motivated by (\ref{g41}), let us introduce for abbreviation
\begin{align}\label{g96}
\delta a_{hij}:=\int(e_j+\nabla\phi_{jh}')\cdot(a_h'-a_h)e_i
\end{align}
so that the tensor $\delta a_{h}=\{\delta a_{hij}\}_{ij=1,\cdots,d}$ satisfies
\begin{align}\label{g91}
\tilde\xi\cdot\delta a_h\xi=\int(\xi+\nabla\phi_{h\xi}')\cdot(a_h'-a_h)\tilde\xi\quad\mbox{where}\quad
\phi_{h\xi}':=\xi_j\phi_{jh}'.
\end{align}
We first claim that the non-degeneracy condition (\ref{g90}) translates into
\begin{align}\label{g92}
-\xi\cdot\delta a_h\xi\ge\lambda_0|B_R||\xi|^2\quad\mbox{for all}\;\xi\in\mathbb{R}^d.
\end{align}
Indeed, since by definition (\ref{g91}) of $\phi_{h\xi}'$ and (\ref{g17}) we have
$-\nabla\cdot a_h'(\xi+\nabla\phi_{h\xi}')=0$, which we rewrite as
$-\nabla\cdot a_h'\nabla\phi_{h\xi}'=\nabla\cdot(a_h'-a_h)\xi$, we obtain from (\ref{g91}) the representation
\begin{align*}
\tilde\xi\cdot\delta a_h\xi=\int\xi\cdot(a_h'-a_h)\tilde\xi-\int\nabla\phi_{h\xi}'\cdot a_h'\nabla\phi_{h\tilde\xi}'
\end{align*}
and thus the inequality for the quadratic part
\begin{align*}
\xi\cdot\delta a_h\xi\le\int_{B_R(y)}\xi\cdot(a_0-a_h)\xi.
\end{align*}
Hence (\ref{g90}) indeed translates into (\ref{g92}).
}

We start by post-processing Lemma \ref{L9} by getting rid of the constraint that $R\ge r_*,r_*(y)$.
Indeed, in case of $R\le\max\{r_*,r_*(y)\}$ we apply Lemma \ref{L9} with $\max\{r_*,r_*(y)\}$ playing
the role of $R$ so that (\ref{g85}) turns into 
\begin{align*}
\big(\av_{B_{\max\{r_*,r_*(y)\}}}|\nabla(u'-u)
&-\partial_iu_h(y)\delta a_{ij}\partial_j\partial_kG_h(-y)(e_k+\nabla\phi_k)|^2\big)^\frac{1}{2}\nonumber\\
&\lesssim (\frac{\ell}{|y|})^d(\frac{\max\{r_*,r_*(y)\}}{|y|})^{d+\beta},
\end{align*}
which because of $\av_{B_R}$ $\le(\frac{\max\{r_*,r_*(y)\}}{R})^d\av_{B_{\max\{r_*,r_*(y)\}}}$ yields
\begin{align*}
\big(\av_{B_{R}}|\nabla(u'-u)
&-\partial_iu_h(y)\delta a_{ij}\partial_j\partial_kG_h(-y)(e_k+\nabla\phi_k)|^2\big)^\frac{1}{2}\nonumber\\
&\lesssim (\frac{\ell}{|y|})^d
(\frac{R}{|y|})^{d+\beta}
(\frac{\max\{r_*,r_*(y)\}}{R})^{d+\beta+\frac{d}{2}}.
\end{align*}
Hence we obtain in either case 
\begin{align}\label{g97}
\big(\av_{B_{R}}|\nabla(u'-u)
&-\partial_iu_h(y)\delta a_{ij}\partial_j\partial_kG_h(-y)(e_k+\nabla\phi_k)|^2\big)^\frac{1}{2}\nonumber\\
&\lesssim (\frac{\ell}{|y|})^d
(\frac{R}{|y|})^{d+\beta}
\big(1+\frac{r_*+r_*(y)}{R}\big)^{d+\beta+\frac{d}{2}}.
\end{align}

\medskip

We now pass from a strong-$L^2$-estimate to an estimate of averages, which allows us to get rid
of the corrector in (\ref{g97}):
\begin{align}\label{g95}
\lefteqn{\big|\int\omega\big(\nabla(u'-u)
-\partial_iu_h(y)\delta a_{ij}\partial_j\partial_kG_h(-y)e_k\big)\big|}\nonumber\\
&\lesssim(\frac{\ell}{|y|})^d(\frac{R}{|y|})^{d}
\Big((\frac{R}{|y|})^{\beta}+(\frac{r_*}{R})^\beta\Big)
\big(1+\frac{r_*+r_*(y)}{R}\big)^{\frac{3}{2}d+\beta}.
\end{align}
Indeed, this follows from the obvious estimate $|\int\omega g|$ $\lesssim(\av_{B_R}|g|^2)^\frac{1}{2}$
applied to the field $g$ $:=\nabla(u'-u)-\partial_iu_h(y)\delta a_{ij}\partial_j\partial_kG_h(-y)(e_k+\nabla\phi_k)$,
combined with an argument that the contribution from $\int\omega\nabla\phi_k$ is negligible. For the latter we first
note that
\begin{align}\label{g96bis}
R^d|\nabla^2 G_h(-y)|\lesssim (\frac{R}{|y|})^{d},\quad
|\nabla u_h(y)|\stackrel{(\ref{L06})}{\lesssim}(\frac{\ell}{|y|})^d
\end{align}
and that
\begin{align}\label{g99}
\frac{1}{R^d}|\delta a_{ij}|\lesssim (1+\frac{r_*(y)}{R})^\frac{d}{2}.
\end{align}
In case of $R$ $\ge$ $\max\{r_*(y),r_*'(y)\}$, (\ref{g99}) follows immediately from the
definition (\ref{g55}), the fact that the integral is restricted to $B_R(y)$ by (\ref{g10}),
and the Caccioppoli estimate (\ref{L55}) with $y$ playing the role of the origin and
also applied to the perturbed medium $a'$. In case of $R$ $\le$ $\max\{r_*(y),r_*'(y)\}$
we argue as in the previous paragraph, that is, we apply the above estimate with 
$\max\{r_*(y),r_*'(y)\}$ playing the role of $R$, which is legitimate in view of the only
constraint (\ref{g10}) on $R$. This yields (\ref{g99}) with a rhs given by
$(1+\frac{\max\{r_*(y),r_*'(y)\}}{R})^\frac{d}{2}$, so that it remains to appeal to Lemma \ref{L7}.
We now may turn to $\int\omega\nabla\phi_k$ itself: By integration by parts we obtain
\begin{align}\label{l03}
|\int\omega\nabla\phi_k|=|\int(\phi_k-\av_{B_R}\phi_k)\nabla\omega|
\stackrel{(\ref{g18})}{\lesssim}\frac{1}{R}\big(\av_{B_R}(\phi_k-\av_{B_R}\phi_k)^2\big)^\frac{1}{2}.
\end{align}
In case of $R\ge r_*$, the rhs is controlled by $(\frac{r_*}{R})^\beta$ according to (\ref{L04}).
In the other case, we argue as before to obtain an estimate by $(\frac{r_*}{R})^{\frac{d}{2}+1}$. 
Hence we obtain in either case
\begin{align}\label{l02}
\frac{1}{R}\big(\av_{B_R}(\phi_k-\av_{B_R}\phi_k)^2\big)^\frac{1}{2}
\lesssim(\frac{r_*}{R})^\beta\big(1+\frac{r_*}{R}\big)^{\frac{d}{2}+1-\beta}.
\end{align}
Inserting (\ref{l02}) into (\ref{l03}) and combining with (\ref{g96bis}) and (\ref{g99}) yields
\begin{align*}
\lefteqn{\big|\partial_iu_h(y)\delta a_{ij}\partial_j\partial_kG_h(-y)\int\omega\nabla\phi_k\big|}\nonumber\\
&\lesssim
(\frac{\ell}{|y|})^d(\frac{R}{|y|})^d\big(1+\frac{r_*(y)}{R}\big)^\frac{d}{2}(\frac{r_*}{R})^\beta
\big(1+\frac{r_*}{R}\big)^{\frac{d}{2}+1-\beta},
\end{align*}
which completes the argument for (\ref{g95}) since $(1+\frac{r_*(y)}{R})^\frac{d}{2}(1+\frac{r_*}{R})^{\frac{d}{2}+1-\beta}$
$\lesssim$ $(1+\frac{r_*+r_*(y)}{R})^{\frac{3}{2}d+\beta}$.

\medskip

We now bring Corollary \ref{Le10} into play, which in our abbreviations (\ref{g96}) and (\ref{g55}) reads
$\frac{1}{R^d}|\delta a_{ij}-\delta a_{hij}|\lesssim(\frac{r_*(y)}{R})^\alpha$ in case of $R\ge r_*(y)$.
In the other case, we use the triangle inequality $\frac{1}{R^d}|\delta a_{ij}-\delta a_{hij}|$ $\le
\frac{1}{R^d}|\delta a_{ij}|$ $+\frac{1}{R^d}|\delta a_{hij}|$ recall the argument for (\ref{g99}) that gave
$\frac{1}{R^d}|\delta a_{ij}|$ $\lesssim(\frac{\max\{r_*(y),r_*'(y)\}}{R})^\frac{d}{2}$
$\lesssim(1+\frac{r_*(y)}{R})^\frac{d}{2}$, while we easily get
$\frac{1}{R^d}|\delta a_{hij}|\lesssim 1$ from (\ref{g22}). Hence in either case, we have
\begin{align*}
\frac{1}{R^d}|\delta a_{ij}-\delta a_{hij}|\lesssim(\frac{r_*(y)}{R})^\alpha\big(1+\frac{r_*(y)}{R})^{\frac{d}{2}-\alpha}.
\end{align*}
This allows us to upgrade (\ref{g95}) to
\begin{align}\label{g98}
\lefteqn{\big|\int\omega\nabla(u'-u)-\partial_iu_h(y)\delta a_{hij}\partial_j\partial_kG_h(-y)e_k\big|}\\
&\lesssim(\frac{\ell}{|y|})^d(\frac{R}{|y|})^{d}
\Big((\frac{R}{|y|})^{\beta}+(\frac{r_*}{R})^\beta+(\frac{r_*(y)}{R})^\alpha\Big)
\big(1+\frac{r_*+r_*(y)}{R}\big)^{\frac{3}{2}d+\beta},\nonumber
\end{align}
where we absorbed both $(\frac{r_*}{R})^\beta$ and $(\frac{r_*(y)}{R})^\alpha$ into
$(\frac{\max\{r_*,r_*(y)\}}{R})^\alpha$.

\medskip

Finally, we will substitute $u_h$ by $(\int g)\cdot \nabla G_h$ in (\ref{g98}).
%
%
Starting point is the following representation of $u_h$ which follows from its definition in (\ref{L24}) and (\ref{L18}):
\begin{align*}
u_h(x)&=\int\nabla G_h(x-x')\cdot g(x')dx'+(\int\nabla\phi_m\cdot g)\partial_mG_h(x)\quad\mbox{and thus}\\
\partial_i u_h(y)&=\int\partial_i\partial_mG_h(y-x') g_m(x')dx'+(\int\nabla\phi_m\cdot g)\partial_i\partial_mG_h(y).
\end{align*}
Hence we have
\begin{align}\label{j01}
\lefteqn{|\partial_i u_h(y)-\big(\int(e_m+\nabla\phi_m)\cdot g\big)\partial_i\partial_mG_h(y)|}\nonumber\\
&\stackrel{(\ref{L97})}{\le}\ell^d\sup_{x'\in B_\ell}|\partial_i\partial_mG_h(y-x')-\partial_i\partial_mG_h(y)|
\lesssim(\frac{\ell}{|y|})^{d+1}.
\end{align}
In order to argue that the contribution from $\int\nabla\phi_m\cdot g$ is negligible we argue as for
(\ref{l03}) and (\ref{l02}) to obtain
\begin{align*}
|\int\nabla\phi_m\cdot g|
\lesssim\ell^d(\frac{r_*}{\ell})^\beta\big(1+\frac{r_*}{R})^{\frac{d}{2}+1-\beta}.
\end{align*}
Hence with help of $|\partial_i\partial_mG_h(y)|\lesssim\frac{1}{|y|^d}$ we may upgrade (\ref{j01}) to
\begin{align}
|\partial_i u_h(y)-(\int g_m)\partial_i\partial_mG_h(y)|
\lesssim(\frac{\ell}{|y|})^{d}\Big(\frac{\ell}{|y|}+(\frac{r_*}{\ell})^\beta\big(1+\frac{r_*}{R})^{\frac{d}{2}+1-\beta}\Big).
\end{align}
Since $|\delta a_{hij}|\lesssim R^d$, see above, and $|\partial_j\partial_kG_h(-y)|\lesssim\frac{1}{|y|^d}$
this allows to pass from (\ref{g98}) to (\ref{j05bis}).
\end{proof}

\medskip


\begin{proof}[\sc Proof of Lemma \ref{Le12}]\

In view of the independence properties of the Poisson point process, 
it follows from elementary probability theory that for any random variable $F$
(ie a function $F=F(X)$ of the point configuration $X$) we have
\begin{align*}
\big\langle(\langle F|B_{R-1}(y)\rangle-\langle F\rangle)^2\big\rangle
=\frac{1}{2}\big\langle\big\langle(\langle F'\rangle_{out}-\langle F\rangle_{out})^2
\big\rangle_{in}\big\rangle_{in}',
\end{align*}
where $\langle\cdot\rangle_{out}:=\langle\cdot|B_{R-1}(y)\rangle$,
where $\langle\cdot\rangle_{in}$ and $\langle\cdot\rangle_{in}'$ denote two independent
copies of the Poisson process restricted to $B_{R-1}(y)$, 
and where $F':=F(X')$ with $X'$ denoting the realization of
the Poisson point process that arises from concatenating $X_{|B_{R-1}^c(y)}$ and $X'_{|B_{R-1}(y)}$.
We apply this to $F=\int\omega\nabla u$ and note that $F'=\int\omega\nabla u'$,
provided $a':=A(\cdot,X')$, cf Definition \ref{D1}. By definition of $X$ and the locality assumption (\ref{l06}) on $A$, 
this is consistent with the second condition in (\ref{g10}); 
the first condition in (\ref{g10}) comes for free by setting $a_0:=a'_{|B_R(y)}$.
Hence in order to establish (\ref{l05}) (with $R$ replaced by $R-1$), we have to show for $R\gg 1$ that
\begin{align}\label{l07}
\big\langle\big\langle(\langle\int\omega\nabla u'\rangle_{out}-\langle\int\omega\nabla u\rangle_{out})^2
\big\rangle_{in}\big\rangle_{in}'
\gtrsim\big((\frac{\ell}{|y|})^d(\frac{1}{|y|})^d\big)^2,
\end{align}
where by homogeneity and rotational invariance, we have assumed wlog
\begin{align}\label{j09}
\int\hat g=e_1.
\end{align}

\medskip

Since by the above remark, we have that $a'=a$ outside of $B_{R}(y)$, we may apply Corollary \ref{Le11}.
More precisely, we apply $\langle\cdot\rangle_{out}$ to (\ref{j05bis}) and obtain by Jensen's inequality
in probability
\begin{align}\label{l08}
\lefteqn{\big|\langle\int\omega\nabla u'\rangle_{out}-\langle\int\omega\nabla u\rangle_{out}}\nonumber\\
&-(\int g_m)\partial_m\partial_iG_h(y)\langle\delta a_{hij}\rangle_{out}\partial_j\partial_kG_h(-y)e_k\big|\lesssim
\langle {\rm err}\rangle_{out},
\end{align}
where we have set for abbreviation
\begin{align*}
\lefteqn{{\rm err}:=(\frac{\ell}{|y|})^d(\frac{R}{|y|})^{d}}\nonumber\\
&\times\Big((\frac{R}{|y|})^{\beta}+(\frac{r_*}{R})^\beta+(\frac{r_*(y)}{R})^\alpha
+\frac{\ell}{|y|}+(\frac{r_*}{\ell})^\beta\Big)\big(1+\frac{r_*+r_*(y)}{R}\big)^{\frac{3}{2}d+\beta}.
\end{align*}
Thanks to the normalization (\ref{j09}), which by (\ref{L97}) turns into
$\int g=\ell^de_1$, we obtain from just considering the first component, and using the symmetry and homogeneity
of $\nabla^2G_h$
\begin{align*}
\lefteqn{\big|(\int g_m)\partial_m\partial_iG_h(y)\langle\delta a_{hij}\rangle_{out}\partial_j\partial_kG_h(-y)e_k\big|}\nonumber\\
&\ge (\frac{\ell}{|y|})^d(\frac{R}{|y|})^d
\big|\partial_i\partial_1G_h(\hat y)\langle\frac{1}{R^d}\delta a_{hij}\rangle_{out}\partial_j\partial_1G_h(\hat y)\big|.
\end{align*}
Hence under the proviso
\begin{align}\label{l10}
-\xi\cdot\langle\frac{1}{R^d}\delta a_{hij}\rangle_{out}\xi\ge\lambda_0|\xi|^2\quad\mbox{for all}\;\xi\in\mathbb{R}^d
\end{align}
for some $\lambda_0>0$ to be chosen later, we obtain
\begin{align*}
\big|(\int g_m)\partial_m\partial_iG_h(y)\langle\delta a_{hij}\rangle_{out}\partial_j\partial_kG_h(-y)e_k\big|
\ge \lambda_0(\frac{\ell}{|y|})^d(\frac{R}{|y|})^d|\nabla^2 G_h(\hat y)e_1|^2,
\end{align*}
which by the invertibility of the Hessian matrix $\nabla^2G_h(\hat y)$ and its continuity on $\{|\hat y|=1\}$ implies
\begin{align*}
\big|(\int g_m)\partial_m\partial_iG_h(y)\langle\delta a_{hij}\rangle_{out}\partial_j\partial_kG_h(-y)e_k\big|
\gtrsim \lambda_0(\frac{\ell}{|y|})^d(\frac{R}{|y|})^d.
\end{align*}
Therefore, always under the proviso (\ref{l10}), we obtain the following lower bound from (\ref{l08})
\begin{align}\label{l11}
\lambda_0(\frac{\ell}{|y|})^d(\frac{R}{|y|})^d
\lesssim \big|\langle\int\omega\nabla u'\rangle_{out}-\langle\int\omega\nabla u\rangle_{out}\big|
+\langle {\rm err}\rangle_{out}.
\end{align}

\medskip

We note that by definition, cf (\ref{g14}), (\ref{g17}) and (\ref{g96}), 
$\delta a_{hij}$ only depends on ${a_0}_{|B_R(y)}=a'_{|B_R(y)}=A(\cdot,X)_{|B_R(y)}$, which by the locality of $A$,
cf (\ref{l06}), implies that $\delta a_{hij}$ and $\langle\delta a_{hij}\rangle_{out}$ are
in particular independent of $X_{|B_{R-1}(y)}$. Hence the proviso (\ref{l10}) is not affected by
(squaring and) applying $\langle\cdot\rangle_{in}$ to (\ref{l11}): Under the proviso (\ref{l10}), we have
\begin{align*}
\big(\lambda_0(\frac{\ell}{|y|})^d(\frac{R}{|y|})^d\big)^2
\lesssim \big\langle\big|\langle\int\omega\nabla u'\rangle_{out}-\langle\int\omega\nabla u\rangle_{out}\big|^2\big\rangle_{in}
+\langle {\rm err}^2\rangle,
\end{align*}
where we've used Jensen's inequality in probability on the last term. In particular, we may multiply this estimate with
the characteristic function $I(-\frac{1}{R^d}\langle\delta a_h\rangle_{out}\ge\lambda_0{\rm id})$ of the event
(\ref{l10}) and then apply $\langle\cdot\rangle_{in}'$, which yields
\begin{align*}
{\rm prob}\big(\lambda_0(\frac{\ell}{|y|})^d(\frac{R}{|y|})^d\big)^2
&\lesssim \big\langle\big\langle\big|\langle\int\omega\nabla u'\rangle_{out}
-\langle\int\omega\nabla u\rangle_{out}\big|^2\big\rangle_{in}\big\rangle_{in}'
+{\rm prob}\langle {\rm err}^2\rangle,
\end{align*}
where we have set for abbreviation
\begin{align*}
{\rm prob}:=\langle I(-\frac{1}{R^d}\langle\delta a_h\rangle_{out}\ge\lambda_0{\rm id})\rangle_{in}'.
\end{align*}
Since $r_*$ has finite moments, cf (\ref{l13}), (and thus also $r_*(y)$ in view of the shift-covariance
(\ref{l22}), which translates into shift covariance of $\nabla(\phi,\sigma)$ and thus of $r_*$, 
in conjunction with the stationarity of $\langle\cdot\rangle$) we have
\begin{align*}
\langle {\rm err}^2\rangle\lesssim 
\Big((\frac{\ell}{|y|})^d(\frac{R}{|y|})^{d}\big((\frac{R}{|y|})^{\beta}+(\frac{1}{R})^\beta+(\frac{1}{R})^\alpha
+\frac{\ell}{|y|}+(\frac{1}{\ell})^\beta\big)\Big)^2.
\end{align*}
Hence provided
\begin{align}\label{l15}
|y|\gg\lambda_0^{-\frac{1}{\beta}}R,\quad
R\gg\lambda_0^{-\frac{1}{\alpha}},\quad
|y|\gg\lambda_0^{-1}\ell,\quad
\ell\gg\lambda_0^{-\frac{1}{\beta}},
\end{align}
we may absorb the second rhs term to obtain
\begin{align}\label{l20}
{\rm prob}\big(\lambda_0(\frac{\ell}{|y|})^d(\frac{R}{|y|})^d\big)^2
\lesssim \big\langle\big\langle\big|\langle\int\omega\nabla u'\rangle_{out}
-\langle\int\omega\nabla u\rangle_{out}\big|^2\big\rangle_{in}\big\rangle_{in}'.
\end{align}

\medskip

In order to derive (\ref{l05}) for a fixed, but sufficiently large $R$, from (\ref{l20}), 
it remains to argue that there exists a $\lambda_0>0$ only depending on the ensemble such that
\begin{align}\label{l17}
{\rm prob}=\langle I(-\frac{1}{R^d}\langle\delta a_h\rangle_{out}\ge\lambda_0{\rm id})\rangle_{in}'\ge\exp(-CR^d)
\end{align}
for some constant $C=C(d)$, where it only matters that the rhs of (\ref{l17}) only depends on $R$
and is positive for every finite $R$.
To this purpose, we first claim that
\begin{align}\label{l14}
\frac{1}{(R-2)^d}\int_{B_{R-2}(y)}(a_h-a')\ge2\lambda_0{\rm id}
\quad\Longrightarrow\quad
-\frac{1}{R^d}\delta a_h\ge\lambda_0{\rm id}.
\end{align}
Here comes the argument: by definition (\ref{g96}) we have
\begin{align}\label{g91}
\tilde\xi\cdot\delta a_h\xi=\int(\xi+\nabla\phi_{h\xi}')\cdot(a_h'-a_h)\tilde\xi\quad\mbox{where}\quad
\phi_{h\xi}':=\xi_j\phi_{jh}'.
\end{align}
Indeed, since by definition (\ref{g91}) of $\phi_{h\xi}'$ and (\ref{g17}) we have
$-\nabla\cdot a_h'(\xi+\nabla\phi_{h\xi}')=0$, which we rewrite as
$-\nabla\cdot a_h'\nabla\phi_{h\xi}'=\nabla\cdot(a_h'-a_h)\xi$, we obtain from (\ref{g91}) and (\ref{g14}) the representation
\begin{align*}
\tilde\xi\cdot\delta a_h\xi=\int_{B_R(y)}\xi\cdot(a_0-a_h)\tilde\xi-\int\nabla\phi_{h\xi}'\cdot a_h'\nabla\phi_{h\tilde\xi}'.
\end{align*}
Hence for the quadratic part we have the inequality
\begin{align*}
\xi\cdot\delta a_h\xi\le\int_{B_R(y)}\xi\cdot(a_0-a_h)\xi,
\end{align*}
which in view of $a_0=a'_{|B_R(y)}$ implies
\begin{align*}
\frac{1}{R^d}\int_{B_{R}(y)}(a_h-a')\ge\lambda_0{\rm id}
\quad\Longrightarrow\quad
-\frac{1}{R^d}\delta a_h\ge\lambda_0{\rm id}.
\end{align*}
This yields (\ref{l14}) since 
\begin{align*}
|\frac{1}{(R-2)^d}\int_{B_{R-2}(y)}(a_h-a')-\frac{1}{R^d}\int_{B_{R}(y)}(a_h-a')|\stackrel{(\ref{L03})}{\lesssim}\frac{1}{R}
\end{align*}
and since we are in the regime (\ref{l15}), which in view of $\alpha\le 1$ includes $\frac{1}{R}\ll\lambda_0$.

\medskip

Based on (\ref{l14}) we now argue that
\begin{align}\label{l16}
\langle I(-\frac{1}{R^d}\langle\delta a_h\rangle_{out}\ge\lambda_0{\rm id})\rangle_{in}'
\ge \langle I(\frac{1}{(R-2)^d}\int_{B_{R-2}(y)}(a_h-a)\ge2\lambda_0{\rm id})\rangle.
\end{align}
This follows immediately from (\ref{l14}) and the observation that thanks to the locality assumption
(\ref{l06}), the event
$\frac{1}{(R-2)^d}\int_{B_{R-2}(y)}(a_h-a')\ge2\lambda_0{\rm id}$ only depends on $X'$ via its restriction to 
$B_{R-1}(y)$.
Hence by the independence property of the Poisson point process, on the one hand, this event implies
$-\frac{1}{R^d}\langle \delta a_h\rangle_{out}\ge\lambda_0{\rm id}$, 
and on the other hand, the conditional probability $\langle\cdot\rangle_{in}'$ of this event agrees with the unconditional one
$\langle\cdot\rangle$.

\medskip

In view of (\ref{l16}), in order to establish (\ref{l17}), it is enough to argue that
there exists a $\lambda_0>0$ such that
\begin{align}\label{l19}
\langle I(\frac{1}{(R-2)^d}\int_{B_{R-2}(y)}(a_h-a)\ge2\lambda_0{\rm id})\rangle\ge\exp(-CR^d).
\end{align}
Here comes the argument for (\ref{l19}):
Letting $\lambda_0>0$ be such that $a_h-A(0,\emptyset)\ge \frac{2\lambda_0}{|B_1|}{\rm id}$, 
cf the monotonicity assumption (\ref{l18}), we have by the locality assumption (\ref{l06}) and shift-invariance assumption
(\ref{l22}) that for any realization $X$ of the Poisson point process:
\begin{align*}
X\cap B_{R-1}(y)=\emptyset&\Longrightarrow& a_h-a\ge \frac{2\lambda_0}{|B_1|}{\rm id}\;\mbox{in}\;B_{R-2}(y)\nonumber\\
&\Longrightarrow&\frac{1}{(R-2)^d}\int_{B_{R-2}(y)}(a_h-a)\ge2\lambda_0{\rm id}.
\end{align*}
Hence (\ref{l19}) follows from the defining property of the Poisson point process:
\begin{align*}
\langle I(X\cap B_{R-1}(y)=\emptyset)\rangle=\exp(-|B_{R-1}(y)|).
\end{align*}
\end{proof}

\medskip


\begin{proof}[\sc Proof of Theorem \ref{T2}]\

As a consequence of the independence property of the ensemble $\langle\cdot\rangle$ of the Poisson point process,
we have for any random variable, in particular $F:=\int\omega\nabla u$, that
\begin{align}\label{l27}
\langle|F-\langle F|B_L\rangle|^2\rangle\ge\sum_{i=1}^N\langle|\langle F|B_R(y_n)\rangle-\langle F\rangle|^2\rangle,
\end{align}
provided we have for the sets
\begin{align}\label{l25}
B_L,\;B_R(y_1),\;\cdots\;,B_R(y_N)\quad\mbox{are pairwise disjoint}.
\end{align}
Hence Theorem \ref{T2} follows immediately from Lemma \ref{Le12}, since under the assumptions of the theorem
there exists a family $\{y_n\}_{n=1,\cdots,N}$ of points such that (\ref{l25}) holds while
\begin{align}\label{l26}
|y_n|\le 2L\quad\mbox{and}\quad N\gtrsim L^d,
\end{align}
where $R$ is the order-one radius given by the lemma.
More precisely, we use (\ref{l05}) with $y_n$ playing the role of $y$ and which by the
first property in (\ref{l26}) assumes the form of
\begin{align*}
\big\langle\big|\langle F|B_R(y_n)\rangle-\langle F\rangle\big|^2\big\rangle
\gtrsim\big((\frac{\ell}{L})^d(\frac{1}{L})^d|\int\hat g|\big)^2.
\end{align*}
By the second property in (\ref{l26}) we obtain for the sum
\begin{align*}
\sum_{n=1}^N\big\langle\big|\langle F|B_R(y_n)\rangle-\langle F\rangle\big|^2\big\rangle
\gtrsim\frac{1}{L^d}\big((\frac{\ell}{L})^d|\int\hat g|\big)^2,
\end{align*}
so that (\ref{j17}) follows via (\ref{l27}).
\end{proof}

\bigskip 

\noindent
\textbf{Acknowledgments.}
The work of JL is supported in part by the National Science Foundation under award DMS-1454939. FO learned from Jim Nolen the general strategy (\ref{l27}) for a lower bound on the variance.


\bibliographystyle{plain}
\bibliography{stochhom}

\begin{bibdiv}
\begin{biblist}

\bib{AKM}{article}{
      author={Armstrong, S.~N.},
      author={Kuusi, T.},
      author={Mourrat, J.-C.},
       title={The additive structure of elliptic homogenization},
        date={2017},
     journal={Invent.~Math},
      volume={208},
       pages={999\ndash 1154},
}

\bib{AS}{article}{
      author={Armstrong, S.~N.},
      author={Smart, C.~K.},
       title={Quantitative stochastic homogenization of convex integral
  functionals},
        date={2016},
     journal={Ann. Sci. Ecole. Norm. S.},
      volume={49},
       pages={423\ndash 481},
}

\bib{AL}{article}{
      author={Avellaneda, M.},
      author={Lin, F.-H.},
       title={Compactness methods in the theory of homogenization},
        date={1987},
     journal={Comm. Pure Appl. Math.},
      volume={40},
       pages={803\ndash 847},
}

\bib{BellaGiuntiOttoarXiv}{misc}{
      author={Bella, P.},
      author={Giunti, A.},
      author={Otto, F.},
       title={Effective multipoles in random media},
        date={2017},
        note={preprint, arXiv:1708.07672},
}

\bib{BellaGiuntiOttoPCMINotes}{inproceedings}{
      author={Bella, P.},
      author={Giunti, A.},
      author={Otto, F.},
       title={Quantitative stochastic homogenization: local control of
  homonization error through corrector},
        date={in press},
   booktitle={{Mathematics and Materials, IAS/Park City Mathematics Series,
  Vol. 23}},
      editor={Bowick, M.~J.},
      editor={Kinderlehrer, D.},
      editor={Menon, G.},
      editor={Radin, C.},
   publisher={AMS, IAS/Park City Mathematics Institute, and SIAM},
}

\bib{FischerOttoCPDE}{article}{
      author={Fischer, J.},
      author={Otto, F.},
       title={A higher-order large-scale regularity theory for random elliptic
  operators},
        date={2016},
     journal={Comm. PDE},
      volume={41},
       pages={1108\ndash 1148},
}

\bib{FischerOttoSPDE}{article}{
      author={Fischer, J.},
      author={Otto, F.},
       title={Sublinear growth of the corrector in stochastic homogenization:
  optimal stochastic estimates for slowly decaying correlations},
        date={2017},
     journal={Stoch. Partial Differ. Equ. Anal. Comput.},
      volume={5},
       pages={220\ndash 255},
}

\bib{GNOarXiv}{misc}{
      author={Gloria, A.},
      author={Neukamm, S.},
      author={Otto, F.},
       title={A regularity theory for random elliptic operators},
        date={2014},
        note={preprint, arXiv:1409.2678},
}

\bib{GloriaOttoAP}{article}{
      author={Gloria, A.},
      author={Otto, F.},
       title={An optimal variance estimate in stochastic homogenization of
  discrete elliptic equations},
        date={2011},
     journal={Ann.~Probab.},
      volume={39},
       pages={779\ndash 856},
}

\bib{GloriaOttofiniterange}{misc}{
      author={Gloria, A.},
      author={Otto, F.},
       title={The corrector in stochastic homogenization: optimal rates,
  stochastic integrability, and fluctuations},
        date={2016},
        note={preprint, arXiv:1510.08290v3},
}

\bib{GloriaOttoJEMS}{article}{
      author={Gloria, A.},
      author={Otto, F.},
       title={Quantitative results on the corrector equation in stochastic
  homogenization},
        date={2017},
     journal={J.~Eur.~Math.~Soc.},
      volume={19},
       pages={3489\ndash 3548},
}

\bib{OdenVemaganti}{article}{
      author={Oden, J.~T.},
      author={Vemaganti, K.~S.},
       title={Estimation of local modeling error and goal-oriented adaptive
  modeling of heterogeneous materials: {I}. {E}rror estimates and adaptive
  algorithms},
        date={2000},
     journal={J. Comput. Phys.},
      volume={164},
       pages={22\ndash 47},
}

\end{biblist}
\end{bibdiv}

\end{document}